\documentclass[12pt]{article}   
\usepackage{e-jc-arxiv}

\dateline{Feb 28, 2020}{Sep 16, 2020}{TBD}
                     

\Copyright{The authors. Released under the CC BY-ND license (International 4.0).}         

\usepackage{graphicx}
\usepackage{url}

\expandafter\let\csname equation*\endcsname\relax
\expandafter\let\csname endequation*\endcsname\relax

\usepackage{amssymb}

\usepackage{amsthm}
\usepackage{amsfonts}
\usepackage{url}
\usepackage{color}
\usepackage{xcolor}

\usepackage{breakurl}
\usepackage{rotating}

\newcommand{\seq}{{\bf a}}
\newcommand{\C}{{\mathbb C}}
\newcommand{\E}{{\mathbb E}}
\newcommand{\R}{{\mathbb R}}
\newcommand{\N}{{\mathbb N}}

\newcommand{\Z}{{\mathbb Z}}
\def\Pr{\mathbb{P}}

\theoremstyle{plain}

\usepackage[latin1]{inputenc}
\usepackage{amsmath}
\usepackage{epsfig}
\usepackage{xcolor}
\newcounter{fig}
\def\TwoDstepset#1#2#3#4#5#6#7#8{%
  \begin{picture}(20,20)(-10,-10)
    \put(0,0){\ifx1#1\thicklines\let\x\vector\else\thinlines\let\x\line\fi\x(-1,-1){10}}
    \put(0,0){\ifx1#2\thicklines\let\x\vector\else\thinlines\let\x\line\fi\x(-1,0){10}}
    \put(0,0){\ifx1#3\thicklines\let\x\vector\else\thinlines\let\x\line\fi\x(-1,1){10}}
    \put(0,0){\ifx1#4\thicklines\let\x\vector\else\thinlines\let\x\line\fi\x(0,-1){10}}
    \put(0,0){\ifx1#5\thicklines\let\x\vector\else\thinlines\let\x\line\fi\x(0,1){10}}
    \put(0,0){\ifx1#6\thicklines\let\x\vector\else\thinlines\let\x\line\fi\x(1,-1){10}}
    \put(0,0){\ifx1#7\thicklines\let\x\vector\else\thinlines\let\x\line\fi\x(1,0){10}}
    \put(0,0){\ifx1#8\thicklines\let\x\vector\else\thinlines\let\x\line\fi\x(1,1){10}}
  \end{picture}}

\newcommand{\oeis}[1]{{\href{https://oeis.org/#1}{\footnotesize\bf #1}}}

\title{Stieltjes moment sequences \\ for pattern-avoiding permutations} 

\author{
Alin Bostan 
\\
\small Inria and Universit\'e Paris-Saclay\\[-0.8ex]             
\small 1 rue Honor\'e d'Estienne d'Orves \\[-0.8ex] 
\small 91120 Palaiseau, France\\[-0.8ex]   
\small\tt \href{mailto:alin.bostan@inria.fr}{\tt alin.bostan@inria.fr}         
\and
Andrew Elvey Price\\
\small LaBRI,  Universit\'e Bordeaux I \\[-0.8ex]       
\small 351 cours de la Lib\'eration    \\[-0.8ex]
\small 33405 Talence Cedex, France    \\[-0.8ex]
\small\tt \href{mailto:andrewelveyprice@gmail.com}{\tt andrewelveyprice@gmail.com}  
\and                                                                    
Anthony~John~Guttmann\\
\small School of Mathematics and Statistics  \\[-0.8ex] 
\small The University of Melbourne \\[-0.8ex] 
\small Vic. 3010, Australia\\[-0.8ex]
\small\tt \href{mailto:guttmann@unimelb.edu.au}{\tt guttmann@unimelb.edu.au}         \and
Jean-Marie Maillard\\
\small LPTMC, 
CNRS, Sorbonne Universit\'e \\[-0.8ex]  
\small 4 Place Jussieu, Tour 23, case 121 \\[-0.8ex]
\small 75252 Paris Cedex 05, France
\\[-0.8ex]
\small\tt \href{mailto:maillard@lptmc.jussieu.fr}{\tt maillard@lptmc.jussieu.fr}                                                                                                                                 
}

\begin{document}

\maketitle

\begin{abstract}

A small {set} of combinatorial sequences have coefficients that can be
represented as moments of a nonnegative measure on $[0, \infty)$. Such
sequences are known as \emph{Stieltjes moment sequences}. {They have a
number of nice properties, such as log-convexity, which are useful to
rigorously bound their growth constant from below.}

This article focuses on some classical sequences in enumerative combinatorics,
denoted $Av(\mathcal{P})$, and counting permutations of $\{1, 2, \ldots, n \}$
that avoid some given pattern $\mathcal{P}$. For increasing patterns
$\mathcal{P}=(12\ldots k)$, we recall that the corresponding sequences,
$Av(123\ldots k)$, are Stieltjes moment sequences, and we explicitly find the
underlying density function, either exactly or numerically, by using the
Stieltjes inversion formula as a fundamental tool.

We first illustrate our approach on two basic examples, $Av(123)$ and
$Av(1342)$, whose generating functions are algebraic. We next investigate the
general (transcendental) case of $Av(123\ldots k)$, which counts permutations
whose longest increasing subsequences have length~at most~$k-1$. We show that
the generating functions of the sequences $\, Av(1234)$ and $\, Av(12345)$
correspond, up to simple rational functions, to an order-one linear
differential operator acting on a classical modular form given as a pullback
of a Gaussian $\, _2F_1$ hypergeometric function, respectively to an order-two
linear differential operator acting on the square of a classical modular form
given as a pullback of a $\, _2F_1$ hypergeometric function.

We demonstrate that the density function for the Stieltjes moment sequence
$Av(123\ldots k)$ is closely, but non-trivially, related to the density
attached to the distance traveled by a walk in the plane with $k-1$ unit steps
in random directions.

Finally, we study the challenging case of the $Av(1324)$ sequence and give
compelling numerical evidence that this too is a Stieltjes moment sequence.
Accepting this, we show how rigorous lower bounds on the growth constant of
this sequence can be constructed, which are stronger than existing bounds. A
further unproven assumption leads to even better bounds, which can be
extrapolated to give an estimate of the (unknown) growth constant.

\end{abstract}

\vskip .1cm


\vskip .1cm

\section{Introduction}
\label{introduction}

\subsection{Context and motivation}

Many important distributions in probability theory have moments that
happen to coincide with some classical counting sequences in combinatorics.
For instance, the moments of the standard exponential distribution
$\mathcal{E}(1)$ with rate parameter 1, 
\begin{equation}\label{ex:exp}
a_n:=\int_{0}^\infty  x^n\, \mu(x) \, dx, \qquad
\hbox{with}\quad \mu(x):= e^{-x},
\end{equation}
are equal to $a_{n}= n!$, and thus count the number of permutations of $\mathfrak{S}_n = \{1,
\ldots, n\}$. Similarly, the moments of the standard Gaussian (normal)
distribution $\mathcal{N}(0,1)$,
\begin{equation}\label{ex:normal}
a_n:=\int_{0}^\infty  x^n\, \mu(x) \, dx, \qquad
\hbox{with}\quad \mu(x):=\frac{1}{\sqrt{2\pi}}e^{-x^2/2}
\end{equation}
satisfy, for all $n\ge0$,
\begin{equation}\label{e2}
a_{2n}=1\cdot3\cdot5\cdots(2n-1), \qquad
a_{2n+1}=0,
\end{equation}
so that $a_{2n}$ counts the number of pairings of~$2n$ elements, or,
equivalently, the number of perfect matchings of the complete graph~$K_{2n}$.

Such connections between probability and combinatorics turn out to be useful
in a number of contexts. For instance, Billingsley's proof of the Central
Limit Theorem by the ``method of moments''~\cite[p.~408--410]{Billingsley86}
arrives at a limiting Gaussian distribution~\eqref{ex:normal}, for sums of
independent random variables, precisely thanks to the property of the
corresponding moments to enumerate pairings~\eqref{e2}.

In random graph theory, the Poisson distribution plays an important role.
When the mean is $1$, its moments
\begin{equation}\label{e3}
a_n:=\E(X^n),\qquad\hbox{where}\quad \Pr(X=k)=\frac{e^{-1}}{k!},
\end{equation}
count the  number of set  partitions of a collection  of~$n$ elements,
giving rise to the sequence of Bell numbers~\cite[p.~109--110]{FlSe09}:
\begin{equation}\label{Bell}
a_n = n!\times \left(\text{coefficient of } \; x^n \, \text{in} \; e^{e^x-1} \right).
\end{equation}
It could be argued that the Poisson distribution arises, in many cases, from
the emergence of set partitions, in accordance with Equations~\eqref{e3}
and~\eqref{Bell}. 

Random matrix theory also abounds with examples, the most famous case
being Wigner's \emph{semi--circle} law. Here the Catalan numbers, $(C_n)_{n \geq 0}=(1,1,2,5,14,42,\ldots)$, known to
enumerate a variety of trees and lattice paths and given by the formula
\begin{equation}\label{catwig}
C_n=\frac{1}{n+1}\binom{2n}{n}
\end{equation}
appear, in the asymptotic limit,  as (renormalized) expectations of powers of traces of random matrices in
various ensembles. Their moment representation
\[
C_{n+1}=\frac{1}{2\pi} \int_0^4 x^n\, \sqrt{x(4-x)}\,  dx
\]
then implies convergence of the spectrum of the random matrices under
consideration to the semicircle law with density~$\frac{1}{2\pi}
\sqrt{x(4-x)}$.

\subsection{The main problem}

There are thus good reasons to develop methods aimed at explicitly
solving the so-called {\em moment problem}~\cite{StTa43,Akhiezer65}, which
considers a sequence $\seq = (a_n)_{n \ge 0}$ of real nonnegative numbers, and
searches to express its general term $a_n$ as the integral
	\[a_n = \int_\Gamma x^n d\rho(x) \qquad  \text{for all} \;\; n \ge 0,\] 
for some support $\Gamma \subseteq \mathbb{R},$ and for some probability
measure~$\rho$. In most cases one can write $d\rho(x) = \, \mu(x) \, dx,$
where $\mu(x)$ is a nonnegative function known as the {\em probability density
function} of $\rho$ (in short, {\em density}), in which case the above
equation becomes
\begin{equation}\label{eq:nth-as-mom}
a_n = \, \int_\Gamma x^n \, \mu(x) \, dx.
\end{equation}

When in addition the support $\Gamma$ is a subset of the
half-line~$[0,\infty)$, the problem~\eqref{eq:nth-as-mom} is classically
called the \emph{Stieltjes moment problem}, and the sequence $\seq \equiv
(a_n)_n$ a \emph{Stieltjes moment sequence}\footnote{Equivalently: for all
$n$, the term $a_n$ is the expectation $\mathbb{E}[X^n]$ of a random variable
$X$ on $\Gamma$, taking real nonnegative values.}. The same
problem is called {the} \emph{Hausdorff moment problem} when $\Gamma=[0,1]$ and {the}
\emph{Hamburger moment problem} when $\Gamma = \mathbb{R}$, see~\cite{Kjeldsen93}.

In this article, we will not consider the Hamburger moment problem, but only
the Stieltjes and the Hausdorff problems. Note that a probability measure is
entirely determined by its moments
(i) in the Stieltjes case if $|a_n| \leq n!$ for
all~$n$~\cite[Thm.~30.1]{Billingsley86},
(ii) in the Hausdorff case, and more generally in the Stieltjes case with
compact support~\cite{BD04}, so that the problem~\eqref{eq:nth-as-mom} has a
unique solution in those cases.

There are several necessary and sufficient conditions that the sequence $\seq
$ must satisfy in order to be a Stieltjes moment sequence, or equivalently,
for a density function $\mu(x)$ to solve~\eqref{eq:nth-as-mom} for some
support $\Gamma \subseteq [0,\infty]$. To express these conditions, one
classically attaches to the sequence $\seq $ the (infinite) Hankel matrices
\[
H_n^\infty(\seq
)=
  \begin{bmatrix}
    a_n & a_{n+1} & a_{n+2} & \ldots \\
    a_{n+1} & a_{n+2}& a_{n+3} & \ldots\\
   a_{n+2} & a_{n+3}& a_{n+4} & \ldots\\
\vdots & \vdots & \vdots & \ddots \\
  \end{bmatrix}. 
\]
The following classical result (see also~\cite[Thm.~2.8]{FJS17} and
\cite[Thm~2.2]{Sokal}), characterizes the property of being a Stieltjes moment
sequence by using the Hankel matrices $H_n^\infty$. It was proved partly in
1894 by Stieltjes~\cite{St} and partly in 1937 by Gantmakher and
Krein~\cite{GK}. In particular, properties (a) and (d) were shown to be
equivalent in~\cite{St}, while these were later shown to be equivalent to (b)
and~(c) in~\cite[Th.~9]{GK}, see also: (i) for
(a)$\Longleftrightarrow$(b),
\cite[Thm.~1.3]{StTa43}, % a <=> b
\cite[Thm.~3.12, p.~65]{Schmudgen}, % a <=> b
\cite[Thm.~87.1, p.~327]{Wall}, % a <=> b
\cite[Thm.~3.12, p.~65]{Schmudgen}, % a <=> b
and \cite[Thm.~1, p.~86]{Simon}; % a <=> b
(ii) for (b)$\Longleftrightarrow$(c),
\cite[Thm.~9]{GK}, % b <=> c
and~\cite[Thm.~4.4]{Pinkus}; % b <=> c
(iii) for (a)$\Longleftrightarrow$(d),
\cite[Appendix, Theorem~0.4, p.~237]{Akhiezer65}, % a <=> d
and~\cite[Satz~4.14, p.~230; Satz~3.11, p.~120]{Perron}. % a <=> b <=> d

{\theorem \label{thm:Stieltjes} For a sequence $\seq
 \equiv (a_n)_{n \ge 0}$,
  the following are equivalent:
\begin{itemize}
\item[\emph{(a)}] There exist $\Gamma \subseteq [0,\infty)$
and a nonnegative measure $\rho$ on $\Gamma$ such that
\[a_n = \, \int_\Gamma \, x^n \, d\rho(x).\]
\item[\emph{(b)}] The matrices $H_0^\infty(\seq)$ and $H_1^\infty(\seq
)$  are both positive semidefinite (i.e., all their leading principal minors, {called Hankel determinants}, are nonnegative).\\
\item[\emph{(c)}] The matrix $H_0^\infty(\seq)$ is totally 
nonnegative
    (i.e., all of its minors are nonnegative). \\
\item[\emph{(d)}] There exists a sequence of nonnegative real numbers
    $(\alpha_n)_{n \ge 0}$,  such that the generating
    function $A(x)\, =  \,\sum_{n=0}^\infty a_n x^n $ of the sequence $\seq
$  satisfies
    \[A(x) \, = \,  \,
    \cfrac{\alpha_0}{1 \, -\cfrac{\alpha_1 x}{1- \, \cfrac{\alpha_2x}{
      \begin{array}{@{}c@{}c@{}c@{}}
        1 \, - \cdots 
      \end{array}
    }}}. \]
\end{itemize}
}
\noindent 
Moreover, the leading principal $n\times n$ minors $\Delta_0^n(\seq)$ 
of $H_0^\infty(\seq)$ and $\Delta_1^n(\seq)$ of $H_1^\infty(\seq )$ occurring in (b) are related to the coefficients $\alpha_i$ in the continued fraction~(d) by: 
\begin{align*}\Delta_0^n(\seq) = \alpha_0 \cdot (\alpha_1 \alpha_2)^{n-1}\cdot (\alpha_3 \alpha_4)^{n-2} \cdots (\alpha_{2n-3} \alpha_{2n-2}),\\
\Delta_1^n(\seq) = \alpha_0 \cdot \alpha_1^n \cdot (\alpha_2 \alpha_3)^{n-1}\cdot (\alpha_4 \alpha_5)^{n-2} \cdots (\alpha_{2n-2} \alpha_{2n-1}).
\end{align*} 
One calls $\seq $ a {\em Stieltjes moment sequence} if it satisfies the
conditions of Theorem~\ref{thm:Stieltjes}.

\subsection{Examples}

\medskip Geometric sequences $(\tau^n)_{n\geq 0}$ with growth rate~$\tau>0$
are Stieltjes moment sequences, as~$\tau^n$ is equal to the $n$-th moment of
the distribution with support $\Gamma={\{\tau\}}$ which is~$\tau$ with
probability 1; equivalently, one can take $\alpha_0=1, \alpha_1=\tau$ and
$\alpha_k=0$ for $k\geq 2$ in~(d). However, other very simple sequences with
rational generating function are not Stieltjes; for instance, the Fibonacci
sequence $\seq = (1,1,2,3,5,8,\ldots)$ with generating function $1/(1-x-x^2)$,
is such that the second leading principal minor $\Delta_1^2(\seq)$ of
$H_1^\infty(\seq )$ equals $1\times 3 - 2 \times 2 = -1<0$, and thus by (b) it
cannot be Stieltjes. It is trivially a moment sequence for a linear
combination of Dirac measures, hence it is a Hamburger moment sequence (with
support $\Gamma$ included in $\mathbb{R}$ but not in $[0,\infty)$). More
generally, Stieltjes moment sequences with rational generating functions are
well-understood: for instance, if $(a_n)_{n\geq 0}$ has a rational generating
function with simple poles only, then $a_n$ is a finite sum of the form
$\sum_j c_j p_j^n$ for some complex $c_j, p_j$, the corresponding measure
$\rho(x)$ is a weighted sum of Dirac measures $\sum_{j} c_j \delta(x-p_j)$,
and one has that $(a_n)_{n\geq 0}$ is Stieltjes if and only if all $c_j$ and
$p_j$ are \emph{nonnegative real numbers}.

\smallskip A useful refinement of Theorem~\ref{thm:Stieltjes} (see e.g., the
references cited before its statement), which essentially excludes Stieltjes
moment sequences whose generating functions are rational, is that: $(a')$
$\seq \equiv (a_n)_{n \ge 0}$ is a Stieltjes moment sequence with a
representing measure $\rho$ \emph{having infinite support}~$\Gamma$
(containing a dense subset of $\mathbb{R}$) $\Longleftrightarrow$ $(b')$
$H_0^\infty(\seq)$ and $H_1^\infty(\seq )$ are \emph{both positive definite}
(i.e., all their \emph{leading principal minors are positive})
$\Longleftrightarrow$ $(c')$ $H_0^\infty(\seq)$ is \emph{totally positive}
(i.e., all of \emph{its minors are positive}) $\Longleftrightarrow$ $(d')$ all
the \emph{$\alpha_i$'s are positive}.

\smallskip Perhaps the simplest example of a Stieltjes moment sequence with
irrational generating function is $a_n=\tau^{n+1}/(n+1)$ with $\tau>0$, which
corresponds to the uniform distribution $\mathcal{U}(0,\tau)$, with
$\Gamma=[0,\tau]$ and $\mu(x)\equiv 1$. Observe that the generating function
$\sum_{n=0}^\infty a_n x^n = -{ {\ln \left( 1-\tau\,x \right) }/{x}}$ is not
only irrational but already transcendental (i.e., non-algebraic). However, it
is D-finite, i.e., it satisfies a linear differential equation with polynomial
coefficients. Another basic example of Stieltjes sequence is $n!$,
corresponding to the standard exponential distribution~$\mathcal{E}(1)$ for
which, according to~\eqref{ex:exp}, one can take $\Gamma = [0,\infty)$ and
$\mu(x) = \exp(-x)$. This can also be seen by using $(d)$ and Euler's 1746
continued fraction expansion~\cite[\S21]{Euler}, with $\alpha_0=1$ and
$\alpha_{2k-1}=\alpha_{2k}=k$ for $k\geq 1$. Once again, the generating
function $\sum_{n \geq 0} n! x^n$ is D-finite, but not algebraic. More
generally, Euler~\cite[\S26]{Euler} showed that for the sequence with general
term $c\, (c+1)\cdots(c+n-1)$ one can take $\alpha_{2n-1}=c+n-1$ and
$\alpha_{2n}=n$, so that the sequence is again Stieltjes for any real $c>0$,
and {its} generating function is transcendental D-finite and satisfies the
linear differential {equation} {$L(y(x))=0$, where $L=$} ${x}^{2} D_x^2+ \left( (c+2)\,x-1 \right) D_x+c=
\left({x}^{2}D_x+2\,x \right) \left(D_x+{\frac {c\,x-1}{{x}^{2}}}\right)$.
Here, and in all that follows, $D_x$ denotes the usual derivation operator
$\frac{\partial}{\partial x}$.

\smallskip A slightly more involved example is ${4^n \, (n+1)! \,
n!}/{(2n+3)!}$ which is Stieltjes using the support $\Gamma=[0,1]$ and the
algebraic density $\mu(x)=\sqrt{1-x}/4$ and whose generating function is
D-finite transcendental, equal to \[\frac16 \, {_2F_1} \left(1, 1;
 \frac52;x \right).\] Here, and in all that follows
${_2F_1} \left(a, b;  c; x \right)$ denotes the
classical Gauss hypergeometric function $\sum_{n=0}^\infty
\frac{(a)_n(b)_n}{(c)_n} \, \frac{x^n}{n!}$ defined for any $a,b,c \in
\mathbb{Q}$, $-c \notin \mathbb{N}$, where $(x)_n$ denotes the Pochhammer
symbol $(x)_n=x\, (x+1)\cdots(x+n-1)$ for $n\in\mathbb{N}$.

\smallskip One of the simplest Stieltjes moment sequences with an algebraic
but irrational generating function is the aforementioned Catalan
sequence~\eqref{catwig}, which has a plethora of combinatorial
interpretations~\cite{Stanley-Catalan}. As we will see in
Section~\ref{sec:Catalan}, in order to show that the Catalan sequence is
Stieltjes, one can take in~$(a)$ the support $\Gamma=[0,4]$ and the density
$\mu(x)=\sqrt{\frac{4-x}{x}}$. Equivalently, one can show that in $(b)$, all
the leading principal minors $\Delta_0^n(\seq)$ and $\Delta_1^n(\seq)$ of
$H_0^\infty(\seq)$ and $H_1^\infty(\seq )$ are equal to~1, a property that
uniquely characterizes the Catalan
sequence~\cite[Pb.~A35(b)]{Stanley-Catalan}. Equivalently, all the
$\alpha_i$'s in the continued fraction expansion $(d)$ are equal to 1.

\subsection{Main objectives}

\medskip In this article we are mainly interested in identifying combinatorial
sequences that are Stieltjes moment sequences, and in finding \emph{an
explicit description of the corresponding densities}: either a closed formula,
or a defining algebraic or differential (potentially nonlinear) equation, or a
numerical approximation if none of those are available. We will mainly focus
on sequences that occur in connection with pattern-avoiding permutations. The
corresponding sequences, denoted $Av(\mathcal{P})$, count permutations of
$\mathfrak{S}_n=\{1, 2, \ldots, n \}$ that avoid some given pattern
$\mathcal{P}$. In that context, the Catalan sequence emerges as the simplest
case, corresponding to the enumeration of permutations avoiding any fixed
pattern of length 3, such as $(213)$ or $(312)$~\cite{SiSc85,K01}.

\smallskip One reason for attempting to identify combinatorial sequences as
Stieltjes (or Hausdorff) moment sequences is that such sequences are
log-convex\footnote{This is not true for Hamburger moment sequences, as
illustrated by the Fibonacci sequence $\seq = (1,1,2,3,5,8,\ldots)$, which
satisfies Cassini's identity $ \, a_{n+1}a_{n-1} \, -a_n^2 = (-1)^{n+1}$ for
all $n\geq 1$.}. To see that a Stieltjes moment sequence is log-convex, it
suffices to observe that for each $n \ge 1$, the expression $ \,
a_{n+1}a_{n-1} \, -a_n^2 \ge \, 0$ is a minor of $H_0^\infty(\seq )$, so it is
nonnegative by condition~$(c)$. Alternatively, one can deduce log-convexity
directly from condition $(a)$ as follows: observe that for any real numbers
$\beta$ and $\gamma$ the expression 
\[ \beta^2\, a_{n-1} \, +2\beta \gamma\,
a_n \, +\gamma^2 \, a_{n+1} = \, \int \, x^{n-1} \, ( \beta + \gamma x)^2 \,
d\rho(x) \] 
is nonnegative. Then setting $\beta =\, a_n$ and $\gamma = \,
-a_{n-1}$ yields $ \, a_{n+1}a_{n-1} \, -a_n^2 \ge \, 0$. If the sequence
$\seq $ is positive, then its log-convexity implies that the ratios $\frac
{a_n}{a_{n-1}}$ are lower bounds on the growth rate $\tau$ of the
sequence\footnote{Actually, Stieltjes moment sequences  \emph{are
infinitely log-convex},~\cite{WZ16}.}.

\smallskip In this article, we will only deal with sequences having at most
exponential growth. An additional characterization of the Stieltjes
property~\cite{LiuPego16,Price18} is available to study problems on such
sequences, and is given in Theorem~\ref{thm:pick} below~(taken
from~\cite[Theorem~6.6]{Price18}; the proof is based
on~\cite[Corollary~1]{LiuPego16}). To state it, let us recall that a
\emph{Nevanlinna function}, (sometimes called a \emph{Pick function}, an
\emph{R function} or a \emph{Herglotz function}) is a complex analytic
function on the open upper half-plane $\mathbb{C}^{+}$ and has nonnegative
imaginary part. A Nevanlinna function maps $\mathbb{C}^{+}$ into itself, but
is not necessarily injective or surjective.

{\theorem \label{thm:pick}
For a sequence of real numbers $\seq
 = (a_n)_{n \ge 0}$, and a real number $\tau > 0$, the following assertions are equivalent:
\begin{itemize}
\item[\emph{(a)}] $\seq
$ is a Stieltjes moment sequence with exponential growth rate at 
most~$\tau$,
\item[\emph{(b)}] 
There exists a positive measure $\rho$ on $[0,\tau]$ such that, for each
$n\geq 0$, \[ a_n = \int_0^\tau x^n d\rho(x),\]
\item[\emph{(c)}] The analytic continuation of the generating function $f(z) = \, \sum_{n=0}^\infty \, a_n z^n$ of $\seq$
is a Nevanlinna function which is 
analytic and nonnegative on $(-\infty,1/\tau).$
\end{itemize}
}

\subsection{Contributions and structure of the article}

The paper is organized as follows. In Section~\ref{basic} we begin by
recalling a basic but extremely useful tool, the \emph{Stieltjes inversion
formula}. We review in~\S\ref{ssec:Poisson} a first proof which is essentially
the one given by Stieltjes himself in~\cite{St}, and then give
in~\S\ref{ssec:Hankel} a second instructive proof based on a complex analytic
approach, using the Hankel contour technique.

In \S\ref{ssec:algalg} we discuss some structural consequences of the
inversion formula, ensuring that algebraic properties of the generating
function of moments are nicely (and algorithmically!) preserved by the density
function. We illustrate this in several ways in~\S\ref{sec:Catalan} on the toy
example of the Catalan sequence, and in~\S\ref{1342} on a first example coming
from the world of pattern-avoiding permutations of length greater than~3.

In Section~\ref{Increasing} we study our basic objects, permutations in
$\mathfrak{S}_n$ avoiding the pattern $(12\ldots k)$, and effectively compute
the corresponding density functions. The case $k=3$ corresponds to the Catalan
sequence, already treated in~\S\ref{sec:Catalan}. The cases $k=4,5,6,7,8$ are
respectively considered in \S\ref{Av1234}--\ref{Av1234as};
\S\ref{Av12345}--\ref{Av12345as}; \S\ref{Av123456}; \S\ref{Av1234567};
\S\ref{Av12345678}.

In Sections~\S\ref{Av12345as:rel} and \S\ref{Av123456as:rel} we show that the
density function for the Stieltjes moment sequence counting permutations that
avoid the pattern $(12\ldots k)$ is closely, but non-trivially, related to the
density attached to the distance traveled by a walk in the plane with $k-1$
unit steps in random directions.

In Section~\ref{1324} we consider the much more challenging (and still
unsolved!) case of permutations avoiding the pattern $(1324)$. We provide
compelling numerical evidence (but not a proof) that the corresponding
counting sequence is a Stieltjes sequence. This is based on extensive
enumerations given in~\cite{CGZ18}, in which the first 50 terms of the
generating function are found. Assuming that the sequence can indeed be
expressed as a Stieltjes moment sequence, we numerically construct the density
function, and obtain lower bounds to the growth constant that are better than
existing rigorous bounds. By extrapolation, we make a rather precise
conjecture as to the actual growth rate.

\medskip 
\subsection{Related work}

\paragraph{Stieltjes moment sequences.} 

A basic tool is the so-called \emph{Stieltjes inversion formula}, invented by
Stieltjes in his pioneering work~\cite{St}, see also the references cited
in~\S\ref{basic}. It has been notably used to study the moment problem in many
\emph{algebraic cases}, that is for sequences whose generating function is an
algebraic function, such as the Catalan sequence~\eqref{catwig}, and many
variants, see the references cited at the end of \S\ref{ssec:algalg}. Various
generalizations of the Catalan numbers are proved to be Stieltjes sequences
e.g., in~\cite{PS01b,LSW11,MP11,Mlotkowski,GP13,LMW16,Katz,LWZ18,Lin19} to
name just a few, using the inverse Mellin transform (a tool similar to the
Stieltjes inversion) and Meijer G-functions.

At a higher level of complexity lie sequences whose generating function is
D-finite but \emph{transcendental} (that is, not algebraic). For instance, the
sequence $n!^3$ is proved to be Stieltjes in~\cite[Eq.~(64)]{KPS01}, with an
explicit density function involving Bessel functions~\cite[Eq.~(65)]{KPS01}.
Various interesting combinatorial sequences arise as binomial sums, and a
natural question is to study if they are Stieltjes or not. For instance, the
sequence $g_n=\sum _{k=0}^{n} {2\,k\choose k} {n\choose k} ^{2}$ is proved to
be Stieltjes in~\cite[p.~7]{ZhuSun}. The same paper proves that the Franel
sequence $f_n = \sum_{k=0}^n \binom{n}{k}^3$ is Hamburger and the
paper~\cite{LWZ18} conjectures that it is Stieltjes. In the same vein, the
famous Ap\'ery sequence $A_n = \sum_{k=0} ^n \binom{n}{k}^2
\,\binom{n+k}{k}^2$ has been conjectured by Sokal to be Stieltjes, and
apparently proved to be so in an unpublished work by Edgar~\cite{Sokal}.
The sequence $g_n$ counts the moments of the distance from the origin of a 3-step
random walk in the plane. More generally, Borwein et
al.~\cite{BNSW11,BSWZ12,BSW13} studied the densities of uniform $n$-step
random walks with unit steps in the plane, corresponding to the moment
sequence~$W_n(k)$ in~\eqref{eq:BNSW1}. Among the aforementioned references,
these papers by Borwein are the closest in spirit to the present article: they
discover features of the corresponding densities, like modularity and
representations via hypergeometric functions (for 3 and 4 steps), and perform
a precise numerical study of their properties. The tools used
in~\cite{BNSW11,BSWZ12,BSW13} are again Meijer G-functions and inverse Mellin
transforms. Note that Bessel functions emerge again, see also~\cite{BBBG08}, a
fact that we will also pop up in our article.

Not all combinatorial sequences have D-finite generating functions; at an even
higher level of complexity lie those {whose generating function is
\emph{D-algebraic} (i.e., the solution of a \emph{nonlinear} differential
equation) or whose exponential generating function is D-algebraic}, but
non-D-finite. 
{
For instance, the Bell sequence $(B_n)_{n \geq 0}$
defined as in~\eqref{Bell} by                                   
$\sum_{n\geq 0} B_n x^n/n! = \exp (e^x-1)$ 
belongs to this
{latter} class, as {do} the Euler and the Springer sequences  
$(E_n)_{n \geq 0}$ and $(S_n)_{n \geq 0}$
defined by $\sum_{n\geq 0} E_n x^n/n! = \sec x + \tan x$ and
$\sum_{n\geq 0} S_n x^n/n! = {1}/(\cos x - \sin x)$.
The five integer sequences 
$(B_n)_{n \geq 0}$, 
$(E_{2n})_{n \geq 0}$,    \label{Euler}
$(E_{2n+1})_{n \geq 0}$, 
$(S_{2n})_{n \geq 0}$
and $(S_{2n+1})_{n \geq 0}$ 
have been proved to be Stieltjes moment sequences~\cite{PBGHSS04,LMW16,Sokal}
and exponential versions of their generating functions are known to be        D-algebraic\footnote{{To be more precise, the five power series $\sum B_{n} x^{n}/n!$, $\sum E_{2n} x^{2n}/(2n)!$, $\sum E_{2n+1} x^{2n+1}/(2n+1)!$, $\sum S_{2n} x^{2n}/(2n)!$, $\sum S_{2n+1} x^{2n+1}/(2n+1)!$ are D-algebraic.}}.
However, the full sequences 
$(E_{n})_{n \geq 0}$ 
and $(S_{n})_{n \geq 0}$
\emph{are not} Stieltjes moment sequences~\cite{Sokal}\footnote{{Interestingly, it is noted  in~\cite{Sokal} that $(S_{n})_{n \geq 0}$ is a Hamburger moment sequence, but not~$(E_{n})_{n \geq 0}$.}}.
}A complete picture, classifying combinatorial
Stieltjes moment sequences according to the (differential) transcendental
nature of their generating functions is however still missing.

\paragraph{Longest increasing subsequences.} The increasing subsequence
problem covers a vast area, difficult to summarize in just one paragraph. It
goes back at least to the 1930s, see for example \cite{ES35}, and was
popularized in a wide-ranging article by Hammersley~\cite{H72}, in response to
a query of Ulam~\cite{Ulam61}. Such subsequences are studied in many
disciplines, including computer science~\cite{K73, Tarjan, Pratt}, random
matrix theory~\cite{DS94, R98,BDJ99}, representation theory of
$\mathfrak{S}_n$~\cite{S13}, combinatorics of Young tableaux~\cite{K73}, and
physics~\cite{Johansson98,Novak,TracyWidom} and~\cite[Chap.~10]{Forrester}. A
recent book~\cite{Romik} is entirely devoted to the surprising and beautiful
mathematics of the longest increasing subsequences. See also Stanley's
comprehensive survey~\cite{Stanley} of the vast literature on increasing
subsequences in permutations.

Counting permutations whose longest increasing subsequences have length less
that~$k$ is in close connection with counting permutations avoiding the
pattern $(12\ldots k)$. For a summary of results on this topic, see Kitaev's
book~\cite[\S6.1.3, \S6.1.4]{Kitaev}. More comprehensively, Romik's book
\cite{Romik} mentioned above. In order to avoid possible confusions, we stress
that counting \emph{consecutive patterns} in permutations is a related but
different area, in which an occurrence of a consecutive pattern in a
permutation corresponds to a contiguous factor of the permutation. The
corresponding generating functions are different from the ones for
permutations avoiding the pattern $(12\ldots k)$, and are quite well
understood~\cite[\S5]{Kitaev} and~\cite{Elizalde}.

\paragraph{General pattern avoidance.} As mentioned before, avoiding the
pattern $(12\ldots k)$ corresponds to the longest increasing subsequences
problem; avoiding more general patterns, or sets of patterns, is an even
vaster topic. Whole books are dedicated to it~\cite{Kitaev,Bona12} and many
problems remain widely open, including in cases when the pattern is simple in
appearance. This is already the case for the innocently looking
pattern~$(1324)$, for which even the nature of the corresponding generating
function remains unknown, {as does} the growth of the
sequence~\cite{JoNa14,CG15,BBEP17,CGZ18}. Until very recently, it was believed
by some experts that for any pattern, the corresponding generating function
for avoiding permutations would be D-finite. This is known under the name of
\emph{the Noonan-Zeilbe{r}ger conjecture}~(1996). The conjecture was disproved by
Garrabrant and Pak~\cite{GaPa15,Pak} using tools from complexity/computability
theory; their striking result is that there exists a family of patterns {which each have length 80}, such that the corresponding generating function is not
D-finite\footnote{See {Igor Pak's very refreshing} post
``\href{https://igorpak.wordpress.com/2015/05/26/the-power-of-negative-thinking-part-i-pattern-avoidance/}{The power of negative thinking, part I. Pattern avoidance}''.}. Note that this existential result is not yet
complemented by any concrete example, ideally of a single avoided permutation,
such as $(1324)$. We might say that the current understanding is at a similar
level as in (good old) times when the existence of transcendental numbers was
proved, without being able to exhibit a single concrete transcendental number.
Coming back to permutations avoiding the particular pattern $(1324)$, there
are compelling (but still empirical) arguments that the generating function
\emph{is not D-finite} and that its coefficients grow as $(11.60 \pm
0.01)^n$~\cite{CGZ18}, the smallest rigorously proved interval containing the
(exponential) growth constant being $(10.27, \, 13.5)$~{\cite{BBEP17}}.

\subsection{Closure properties}

The class of Stieltjes moment sequences $\mathcal{S}(\Gamma)$ with support
$\Gamma$ is a ring and is closed under forward and backward shift, under
dilation and under division by the sequence $(n+1)_{n \geq 0}$. There are a
few basic transformations on sequences that correspond to simple operations on
densities. For instance, if $b_n=a_{n+1}$ is a sequence obtained by
\emph{forward shift}, then the density associated with $b_n$ is $x\cdot
\mu(x)$; similarly, the \emph{backward shift} $b_n=a_{n-1}$ is associated with
density function $(\mu(x)-\mu(0))/x$, but now suitable integrability
conditions of $\mu$ at~$x\sim 0$ are needed. We list a few such
correspondences in Fig.~\ref{opns-fig}. We only stress the formal aspects,
give the basic versions, and refrain from stating detailed validity
conditions, as these are obvious consequences of change-of-variables formulae,
partial integration, and similar elementary techniques. Many closure
properties are proved by Bennett in~\cite[\S2]{Bennett}.

\def\fs{\footnotesize}
\def\ds{\displaystyle}

\begin{figure}\small
\begin{center}\renewcommand{\arraystretch}{1.6}
\begin{tabular}{ll|ll}
\hline\hline
\multicolumn{2}{c|}{\emph{Transformation on sequences}} &
\multicolumn{2}{c}{\emph{Operation on densities}}\\
\hline
forward shift & $b_n=a_{n+1}$ & $\nu(x)= x \,\mu(x)$ &
\\
backward shift & $b_n=a_{n-1}$ & $\nu(x)=(\mu(x)-\mu(0))/x$ & \fs (requires integrability)
\\
differences & $\ds b_n=\sum \binom{r}{k}(-1)^{k}a_{n+k}$ & $\nu(x)=(1-x)^r \, \mu(x)$ & \fs($r\in\Z_{\ge0}$, $\Gamma=[0,1]$)
\\
\hline
``derivative'' & $b_n=na_{n-1}$ & $\nu(x)=-w'(x)$ & \fs (plus boundary terms)\\
``primitive'' & $b_n=\frac{1}{n+1}a_{n+1}$ & $\nu(x)=-\int \mu(x)$ & \fs (plus boundary terms)\\
\hline 
sum & $b_n=a_n+u_n$ & $\nu(x)=\mu(x)+\lambda(x)$ & \\
product & $b_n=a_n u_n$ & $\nu(x)=\int \mu(s) \, \lambda(x/s) \, \frac{ds}{s}$ & \fs (multiplicative convolution)\\
\hline
dilation & $b_n=a_{rn}$ & $\nu(x)=\frac{1}{r} \, \mu(x^{1/r}) \, x^{-1+1/r}$ & \fs($r\in\Z_{\ge2}$, $\Gamma\subseteq\R_{\ge0}$)\\
\hline\hline
\end{tabular}
\end{center}
\caption{\label{opns-fig}\small
The correspondence between basic transformations on sequences and operations
on densities: here $a_n, u_n, b_n$ are respectively associated with densities
$\mu(x), \lambda(x), \nu(x)$.} 
\end{figure}

\section{The moment problem and the Stieltjes inversion formula}\label{basic}

In what follows, we restrict our attention to combinatorial sequences, that
is, we assume that the coefficients $a_n$ of the sequence $\seq $ are
\emph{nonnegative integers}, counting objects of size~$n$ in some fixed
combinatorial family. We assume that the generating function $f(z)= \,
\sum_{n=0}^\infty a_n z^n$ has radius of convergence $z_c\in(0,\infty)$, with
$z_c=1/\tau$, that $\, f(z)$ is analytically continuable to the complex plane
slit along {$[z_{c},\infty)$}, and that $\, f(z) = \, O(z^{-1})$ as $\, z \,
\rightarrow \, \infty$. Much of this section applies more generally, but to
simplify matters, it is more convenient to consider this more restricted
situation.

Assuming that $a_n$ is equal to the $n$-th moment of a probability measure with density~$\mu$, that is, assuming equation~\eqref{eq:nth-as-mom} holds,
our aim is to solve the \emph{Stieltjes moment problem},
that is, to obtain the density $\mu$ in terms of the generating function $f$ of $\seq
$. In other
words, we need to solve the equation that expresses $f$ in
terms of $\mu$, namely 
\begin{equation} \label{eq:mompb}
	f(z)\,
=\,\, \int_0^\tau\, \frac{\mu(x)}{1-xz}\,\, dx.
\end{equation}
The basic tool used to solve the moment problem is the 
\emph{Stieltjes inversion formula}, invented by Stieltjes in his pioneering
work~\cite[\S39, p.~72--75]{St}, then popularized in the first half of the
20th century in influential books by Perron~\cite[Chap. IV, \S32]{Perron},
Stone~\cite[Chap.~5]{Stone1932}, Titchmarsh~\cite[Chap.~XI]{Titchmarsh},
Widder~\cite[Chap.~III and VIII]{Widder}, Shohat and Tamarkin~\cite{StTa43}
and Wall~\cite[Chap.~XIII, \S65]{Wall}. Various generalizations have been
studied in the second half of the 20th century, see the books by
Akhiezer~\cite[Chap.~3]{Akhiezer65},
Chihara~\cite[Chap.~III]{Chihara1978}, Teschl~\cite[Chap.~2]{Teschl}, and also
Masson's article~\cite{Masson1970} and Simon's survey~\cite{Simon}.

\subsection{Via the Poisson kernel}\label{ssec:Poisson}

We follow here the presentation in~\cite[Lecture~2]{NiSp06}.
If $\rho$  is a probability measure with density $\mu$ 
with support $\Gamma \subseteq\mathbb{R}$, then its
\emph{Stieltjes transform} (or, \emph{Cauchy transform}) is the map
$G_\mu: \mathbb{C}^+ \longrightarrow \mathbb{C}^-$
defined by 
\[ G_\mu(z) := \int_\Gamma \,\frac{\mu(x)}{z-x}\, \,dx,\]
where $\mathbb{C}^+$ and $\mathbb{C}^-$ denote respectively the upper and the lower half-planes.

It is easy to show that $G_\mu$ is analytic on $\mathbb{C}^+$, and moreover if the measure $\rho$ is compactly supported (as assumed from the very beginning in our setting), with a support $\Gamma$ (in our case, $\Gamma = [0,\tau]$) contained in the interval $[-R,R]$, then 
\[ G_\mu(z) = \sum_{n \geq 0} \frac{a_n}{z^{n+1}}, \qquad \text{for all} \;\; |z| > R, \]
where $a_n = \int_\Gamma x^n d\rho(x) = \int_\Gamma x^n \mu(x) \,dx$ is the $n$-th moment of $\rho$. (This follows, after multiplication by $\mu(x)$ followed by integration, from the uniformly convergent expansion $1/(z-x) = \sum_{n \geq 0} {x^n}/{z^{n+1}}$ which holds for all $|z|>R\geq |x|$.)

Therefore, Equation~\eqref{eq:mompb} simply states:
Starting from $f$, find $\mu$ such that 
\begin{equation}\label{eq:mompb3}
	G_\mu(z) = \frac{1}{z}
f\left(\frac{1}{z}\right).
\end{equation}

The \emph{Stieltjes inversion formula} is an effective way of solving~\eqref{eq:mompb3}, i.e. of recovering the density $\mu$ of the probability measure $\rho$ from its Stieltjes transform $G_\mu$. Denoting 
\begin{equation}\label{inversemu}
 h_\epsilon(t) := -\frac{1}{\pi} \, \Im \left( G_\mu(t+i \epsilon)\right),
\end{equation}	
this formula reads:
\begin{equation}\label{eq:Stieltjes}
	\mu(t) dt = d\rho(t) = \lim_{\epsilon \to 0^{+}} h_\epsilon(t) dt.
\end{equation}	
(Here, and below, ``$\Im$'' stands for the operation of taking the imaginary part of a complex number.)
The original proof of Stieltjes is based on the so-called
\emph{Poisson kernel on the upper half plane} defined by ${P_\epsilon(t):=\frac{1}{\pi} \frac{\epsilon}{t^2+\epsilon^2}}$. First, by definition,
\[
h_\epsilon(t)  = -\frac{1}{\pi} \, \Im \left( \int_\Gamma \,\frac{1}{t+i \epsilon-x}\, \,d\rho(x)\right)
= -\frac{1}{\pi} \, \Im \left( \int_\Gamma \,\frac{t-i \epsilon-x}{(t-x)^2+\epsilon^2}\, \,d\rho(x)\right)
\]
hence $h_\epsilon(t)$ can be expressed as the convolution integral
\[
h_\epsilon(t)  
= \frac{1}{\pi} \, \left( \int_\Gamma \,\frac{\epsilon}{(t-x)^2+\epsilon^2}\, \,d\rho(x)\right)
=
\int_\Gamma P_\epsilon(t-x) \,d\rho(x)
\]
and the properties of the Poisson kernel permit one to conclude
the proof of~\eqref{eq:Stieltjes}.

In the important case when $G_\mu$ admits a continuous extension to $\C^+
\cup I$,  for some interval $I\subseteq \R$,
the Stieltjes inversion formula simply reads: 
\begin{equation}\label{eq:Stieltjes2}
\mu(x) 
=
-\frac{1}{\pi} \, \Im \left( \lim_{\epsilon \to 0^{+}}  \left(G_\mu(x+i \epsilon) \right) \right), \qquad \text{for all} \; x\in I.
\end{equation}

Putting things together, we conclude that, given a sequence $\seq
 = (a_n)_{n \geq 0}$ and a probability density $\mu(x)$,
the following assertions are equivalent:
\begin{itemize}
\item  $a_n$ is the $n$-th moment of $\mu$, i.e., $\displaystyle{a_n  = \int_\Gamma x^n \mu(x) \,dx}$, for all $n\geq 0$, 
\item  the generating function $f(z)=\sum_{n \geq 0} a_n \, z^n$ is equal to $\displaystyle{\int_\Gamma\, \frac{\mu(x)}{1-xz}\,\, dx}$,
\item $\displaystyle{g(z) = \frac{1}{z} f\left(\frac{1}{z}\right)}$ is equal to 
$\displaystyle{G_\mu(z) = \int_\Gamma \,\frac{\mu(x)}{z-x}\, \, dx}$,
\item $\displaystyle{\mu(x) = -\frac{1}{\pi} \, \lim_{\epsilon \to 0^{+}} \Im \left(g(x+i \epsilon) \right)}$, where $\displaystyle{g(z) = \sum_{n\geq 0} \frac{a_n}{z^{n+1}}}$,
\end{itemize}
in which case, moreover, the following integral representation holds:
\begin{equation}\label{eq:Stieltjes3}
a_n  =  -\frac{1}{\pi} \, \int_\Gamma x^n \, \lim_{\epsilon \to 0^{+}} \Im \left( \frac{1}{x+i \epsilon} f\left(\frac{1}{x+i \epsilon}\right)  \right) \,dx.
\end{equation}

%%%
\subsection{Via the Hankel contour technique}\label{ssec:Hankel}

The previous inversion process~\eqref{eq:Stieltjes3} can also be performed
using Cauchy's coefficient formula in conjunction with special contours of
integration known as \emph{Hankel contours}. These contours are classical in
complex analysis (e.g., to express the inverse of the Gamma function in the
whole complex plane), in the inversion theory of integral
transforms~\cite{Doetsch}, and more recently in \emph{singularity analysis},
see~\cite{FlOd90} and~\cite[Chap.~VI]{FlSe09}. Hankel contours come very close
to the singularities then steer away: by design, they capture essential
asymptotic information contained in the functions' singularities.

The proof which follows is taken from~\cite{BFP11}, and does not appear to be
published elsewhere in the rich literature on the moment problem, which is why
we reproduce it here.
To start with, Cauchy's coefficient formula provides an integral 
representation
\[
a_n=\frac{1}{2 \, \pi \, i } \, \int_{\cal C} f(z)\, \frac{dz}{z^{n+1}},
\]
where the integration path~$\cal C$ should encircle the origin 
and stay within the domain of analyticity of $f(z)$.
(For instance, any circle of radius $<z_c$ is suitable.)

Due to the assumptions made on $\, f(z)$, 
we can extend the contour $\, {\cal C}$
to be part of a large circle centered at the origin and with radius $\, R$,
complemented by a Hankel-like contour
that starts from $\, R \, - \, i \, \epsilon$ to $\, z_c \, -i \, \epsilon$,
then winds around  $\, z_c$ to the left, then continues
from $\, z_c + \, i \, \epsilon$ 
to $\, R \, + \, i \, \epsilon$. As
$\, R \, \rightarrow \, \infty$, only the contribution
of the Hankel part~$\mathcal{H}$ of the contour survives, so that
\[
  a_n \, = \,\, \frac{1}{2 \, \pi \, i }  \, 
  \int_{\mathcal{H}} \, f(z) \, {{d z} \over { z^{n+1}}}.
\]

Since $\, f(z)$ has real Taylor coefficients, it satisfies 
$\, f(\bar{z}) \, = \, \bar{f}(z)$. This was originally true
for $\, |z| \, < z_c$ and it survives in the split
plane, by analytic continuation. Thus, if we let $\, \epsilon$
tend to $\, 0$, we obtain in the limit 
\[  a_n \, = \,\, \frac{1}{2 \, \pi \, i}    \, 
  \int_{\mathcal{H}} \,\left  (f_{-}(z) \, - f_{+}(z)\right ) \,  \, {{d z} \over { z^{n+1}}},
\]
where $\,  f_{-}(z)$,  $\,  f_{+}(z)$ are the lower and upper limits
of $\, f(z)$ as $\, z \, \rightarrow \, x \, \in \,  {\mathbb R}_{+}$ from
below and above, respectively. This is a typical process of contour integration,
used here in a semi-classical way.
By conjugacy, we also have 
\[
   f_{+}(z) \, - f_{-}(z)  \, = \, \, \,  2 \, i \,  \, \, \Im\left(f_{+}(x)\right),
\]
since the contributions 
to the integral 
from the real part of 
$f$ cancel out. By denoting
\[ 
  f_{\Im}(x)  \, = \, \,
  \lim_{\epsilon \, \rightarrow \, 0^{+}} \, \Im(f(x\, + \, i \, \epsilon)),
\]
we obtain the real integral representation
$
\displaystyle{  a_n  \, = \,- \, {{1} \over {\pi}}  \, 
  \int_{z_c}^{\infty} \,  f_{\Im}(x) \,   \, {{dx } \over {x^{n+1}}}.} 
$

Then, under the change of variables $\, t \, = \, 1/x$, one
again obtains the Stieltjes inversion formula~\eqref{eq:Stieltjes2}, under the equivalent form
\[
  a_n  \, = \, - \, {{1} \over { \pi}}  \, 
  \int_{0}^{1/z_c} {t^n}  \cdot \left({1 \over {t}} \,  \, \,  f_{\Im}\Bigl( {{1} \over {t}} \Bigr)  \right) \, dt. 
\]

\subsection{Stieltjes inversion formula for non-convergent series} In the
case where the series $f(z)=\sum_{n=0}^{\infty}a_{n}z^n$ has radius of
convergence $0$, we cannot apply the Stieltjes inversion formula as described,
as $f(z)$ only makes sense as a formal power series, as the sum does not
converge except at $z=0$. Nonetheless, we mention a generalization of this
method due to Hardy \cite{Hardy1}, which applies as long as there exists a
constant $C>0$ satisfying
\begin{equation}\label{Hardy_condition}
	a_n\leq C^n \cdot (2n)! \qquad \text{for all} \; n>0.
\end{equation}
The idea is to introduce a new function
\[F(s)=\sum_{n=0}^{\infty}\frac{{a_{n}} \cdot (-s)^{n}}{(2n)!},\]
which {\em does} converge for small $s$, and which extends to an analytic function on $\mathbb{C}\setminus\mathbb{R}^{+}$. One can then prove that $G_{\mu}(z):=\int_{\Gamma}\frac{\mu(x)}{z-x}dx$ is given by
\[G_{\mu}(z)=\frac{1}{z}\int_{0}^{\infty}e^{-t}F\left(\frac{t^2}{z}\right)dt.\]
Finally the Stieltjes inversion formula can be applied to solve for $\mu$.

\subsection{Algebraic and algorithmic consequences of the Stieltjes inversion}
\label{ssec:algalg} 

An important consequence of the Stieltjes inversion formula is that the
properties of (the generating function of) a sequence of moments are perfectly
mirrored by the properties of the corresponding density function. More 
precisely, the following holds:

{\theorem 
Assume as before that $\seq
 = (a_n)_{n \geq 0}$ is the moment sequence of
the density function~$\mu$ on a bounded domain~$\Gamma$,
i.e., $\displaystyle{a_n = \int_\Gamma x^n \mu(x) \,dx}$
for all $n\geq 0$.  
If the generating function $f(z)=\sum_{n \geq 0} a_n \, z^n$ of $\seq$ belongs 
(piecewise) to one of the following classes
\begin{enumerate}
\item[$(i)$] \emph{algebraic}, i.e., root of a polynomial equation $P(z,f(z)) = 0$, with $P\in\C[x,y]$,
\item[$(ii)$] \emph{D-finite}, i.e., solution of a linear ODE with polynomial
coefficients,
\item[$(iii)$] \emph{D-algebraic}, i.e., solution of a nonlinear ODE with polynomial coefficients,
\end{enumerate}
then the same is true for the density function~$\mu$.
}

\bigskip Moreover, one can effectively compute an (algebraic, resp.
differential) equation satisfied by the density function $\mu$ starting from
an equation for~$f$. Instead of proving the general statement, we illustrate
the proof on an example in the next subsection.

\bigskip
We expect that this theorem generalizes to some cases where the support is not bounded, perhaps to all sequences satisfying Hardy's condition \eqref{Hardy_condition}.
Note that the converse of $(ii)$ holds true~\cite{Batenkov,BJL19}. However, this is not the case for $(i)$ and we do not believe it is the case for $(iii)$.  

For instance, concerning the converse of~$(i)$, 
if $\Gamma=[0,1]$ and $\mu(x) = C \cdot
\sqrt [3]{x \left( 1-x \right) }$, where
\[C = \frac{5 \, \sqrt {3}}{\sqrt[3]{2}}\,{\frac {\Gamma \left( \frac23 \right) \Gamma \left( \frac56 \right) }{{\pi}^{\frac32}}} \; \approx \; 1.88682,\]
then the generating function $F(z)$ of the sequence of moments $a_n = \int_\Gamma x^n \mu(x) dx$,
\begin{align*}
F(z) \, = \, \frac{5}{z} \,  \left( 1 - \, \sqrt [3]{1-z} \,\, {\mbox{$_2$F$_1$}\left(\frac13,\frac23;\, \frac53 ; \,z\right)} \right) \qquad \qquad \\  \qquad \qquad = \, 1+{\frac{1}{2}}z+{\frac{7}{22}}{z}^{2}+{\frac{5}{22}}{z}^{3}+{\frac{65}{374}}{z}^{4}+{\frac{26}{187}}{z}^{5}+ \cdots
\end{align*}
is a transcendental function. See~\cite{PRY04, RY05} for a study of the
algebraicity of integrals of the form 
\[g(z)=\frac{1}{2 \, \pi \, i } \int_{\cal C} \mu(x)\, \frac{dx}{x-z}.\] 

As for the converse of~$(iii)$, we do not have a counterexample for which the support~$\Gamma$ is bounded, but for unbounded support, we have the following counterexample: take the log-normal distribution $\mu(x) =
\exp(-\ln(x)^2/4)/x$, which is clearly D-algebraic. Its moments over
$\Gamma=[0,\infty)$ are equal to $a_n = \exp(n^2)$. The generating function
$f(z)$ of $(a_n)_n$ is not D-algebraic, since by a result of Maillet and
Mahler~\cite[Eq.~(2)]{Rubel}, if $f(z)$ were D-algebraic, then there would exist
two positive constants $K$ and $C$ such that $a_n < K \, n!^C$ for all $n
\in \N$, which is clearly not the case.

\smallskip 
Another interesting example is the D-algebraic density 
{$\mu(x) = 1 / \left( 2\, {\sqrt {x}} \, {\cosh \left( \frac{\pi}{2} \,\sqrt {x} \right)} \right)$}, 
whose moment sequence $e_n =\int_0^\infty x^n \mu(x) dx$ is the sequence
$({1}, 1, 5, 61, 1385, 50521, \ldots)$ of Euler's secant numbers
{(these are precisely the numbers $E_{2n}$ mentioned on page~\pageref{Euler})}.
Although the exponential generating function $\sum_{n \geq 0} e_{n}
\frac{z^{2n}}{(2n)!}$ is $\sec(z) = \frac{1}{\cos(z)}$, hence it is
D-algebraic, the ordinary generating function $f(z) = \sum_{n\geq 0} e_{n}
z^n$ is not D-algebraic~{\cite{BDR20}}.

\smallskip 

Finally, we give a potential counterexample to~$(iii)$ which has bounded
support. The density $\mu(x)=\ln(2\sin(x/2))$ on the interval $[0,2\pi]$ is
clearly D-algebraic. However, the corresponding moments
$a_{n}=\int_{0}^{2\pi}x^n \, \mu(x)dx$ are polynomials in $\pi$ whose coefficients
are linear combinations of odd zeta values $\zeta(2k+1)$:
\[a_{n}=\frac{n}{2} \, (2\pi)^{n}\sum_{j=0}^{\lceil\frac{n}{2}\rceil-1}\frac{(-1)^{j+1}}{4^{j}}(2j)!{n-1\choose
2j}\frac{\zeta(2j+3)}{\pi^{2j+1}}.\] 
If the generating series of $(a_n)_n$ were D-algebraic, this would
contradict a commonly accepted number theoretic conjecture saying that the numbers
$\pi$ and $\zeta(2k+1)$ are algebraically independent.

\subsection{A basic, yet important example: The Catalan case} \label{sec:Catalan}

Here, we illustrate the above results and procedures on the case of $Av(123)$,
the permutations avoiding the pattern $(123)$, which are counted by Catalan
numbers $C_n = \frac{1}{n+1}\binom{2n}{n}$~\cite{SiSc85,K01}.
Let $f(z)$ be the generating function of the Catalan sequence,
\[f(z) = 1+z+2\,{z}^{2}+5\,{z}^{3}+14\,{z}^{4}+42\,{z}^{5}+132\,{z}^{6}+\cdots\]

\subsubsection{The first way: exploiting algebraicity}
It is well-known that $f(z)$ is algebraic, being the root of the polynomial $P(x,y)=xy^2-y+1$. 
It follows that  $g(z) = f(1/z)/z$ is also algebraic, 
{to be precise, it is a} root of $K(z,y) = 1-zy+zy^2$.

Assume that $C_n$ is the $n$-th moment of a probability measure with density
$\mu(x)$, so that
\[ C_n = \int_\Gamma x^n \, \mu(x)  \; dx, \qquad \text{for all} \;\;  n\geq 0.\]
As the exponential growth rate of $C_n$ is $\tau=4$, one may further assume
the support of $\mu$ is $\Gamma=[0,4]$. By the Stieltjes inversion formula,
one has
\[ \mu(x) = -\frac{\psi(x)}{\pi} , \qquad \textrm{where \qquad } \psi(x) = \lim_{y 
\rightarrow 0^{+}} \Im \Big(  \, g(x+iy)\Big).\] 

In this particular case, one can solve $P$ by radicals, and obtain directly
that $f(z) = {\frac {1-\sqrt {1-4\,z}}{2z}}$ and that $\mu(x) = \frac{1}{2
\pi} \, \sqrt{\frac{4-x}{x}}.$ We explain now a method that works more
generally in the case when $P$ cannot be solved by radicals. It yields, by a
resultant computation, a polynomial equation satisfied by the density~$\mu$.

Writing $g(x+iy)$ as $A(x,y) + iB(x,y)$, where $A$ and $B$ are real functions,
and denoting {by} $\varphi(x)$ and $\psi(x)$ the limits in $y=0^{+}$ of $A$ and of
$B$, we deduce
\[ 0 = \, K(z,g(z)) = \, K(x+iy, A(x,y) + i B(x,y)) \]
and taking the limit $y \rightarrow 0^{+}$ yields 
$K(x, \varphi(x)+i\psi(x))=0$. We thus have 
\[ 1 \, - \, x \cdot (\varphi(x)+i\,\psi(x)) + x \cdot  (\varphi(x)+i\,\psi(x))^2 = 0 \] 
and by identification of real and imaginary parts, $\varphi$ and $\phi$
satisfy
\begin{equation}\label{sys}
	1 \, - \, x \, \varphi(x) \, + \,  x \cdot  (\varphi(x)^2 - \psi(x)^2) = 0, \quad 
- \, x \, \psi(x) + 2 \, x \,  \varphi(x) \psi(x) =0.
\end{equation}
In particular, $\varphi(x)$ and of $\psi(x)$ are algebraic; moreover,
$\psi(x)$ is {a} root of the resultant
\[\textrm{Res}_u(1 - x u + x \cdot (u^2 - v^2), \, -x v + 2x uv) \, =
\, x^2v^2 \, (4\, (1-xv^2) \, - x). \]
As $\mu(x)$ is nonnegative, the conclusion is that
\[ \mu(x) = - \frac{1}{\pi} \, \psi(x) = \frac{1}{2 \pi} \, \sqrt{\frac{4-x}{x}}. \]
In other words, the previous procedure finds that
\begin{equation} \label{eq:Catalan}
	C_n = \frac{1}{2\pi} \, \int_0^4 x^n \,  \sqrt{\frac{4-x}{x}}  \; dx.
\end{equation}

\smallskip \noindent {\bf Remark.} Once deduced, this kind of equality can be
proved algorithmically by the method of~\emph{creative telescoping}, see
e.g.,~\cite{AlZe90,Koutschan13,BDS16}; in our case, this method finds and
proves that the integrand $U(n,x) = x^n \sqrt{\frac{4-x}{x}}$ is a solution of
the \emph{telescopic equation}
\[ (n+2) \,  U(n+1,x) \, - \, (4n+2) \,  U(n,x) \, = \, \partial_x \Big(x\,  (x-4) \,  U(n,x) \Big).\]
Integrating the last equality w.r.t.~$x$ between 0 and 4 implies a linear
recurrence satisfied by $f_n = \int_0^4 U(n,x) \, dx$:
\[ (n+2) \,  f_{n+1} \, - \, (4n+2) \,  f_n = 0.\]
Since $f_0 = \frac{1}{2\pi} \,  \int_0^4 \sqrt{\frac{4-x}{x}}=1 = C_0$, the
sequence ${f_n}/(2\pi)$ satisfies the same recurrence, and the same
initial conditions, as the sequence $(C_n)$. Therefore, they coincide, and
this proves~\eqref{eq:Catalan}.

\smallskip \noindent {\bf Remark.} Note that the Stieltjes inversion formula
has allowed us to recover and to prove algorithmically the Marchenko-Pastur
(or, free Poisson) density, arising in connection with the famous Wigner
semicircle law in free probability and in random matrix
theory~\cite{VDN92,Voiculescu}! The same approach allows one to prove
that, for any $s\geq 2$, the Fuss-Catalan numbers $\frac{1}{(s-1)n+1}
\binom{sn}{n}$ admit a moment representation over the interval
$[0,s^s/(s-1)^{s-1}]$, where the density function is a positive
\emph{algebraic} function of degree at most~$s(s-1)/2$~\cite{BFP11}, see
also~\cite{PS01,FussCat,BBCC11,LSW11,MP11,MPZ13,MP14,MP18} for related results
and generalizations; these can be obtained using the Stieltjes inversion
formula. In all these references, the generating functions and the
corresponding densities are algebraic functions.

\smallskip \noindent {\bf Remark.} Catalan numbers count Dyck paths (among
many other combinatorial objects). More generally, the previous method works
\emph{mutatis mutandis} for sequences that count lattice paths in the quarter
plane whose stepset is contained in a half-plane. The generating function of
such walks is known to be algebraic~\cite{BaFl02}, hence the associated
measure is algebraic as well. However, the support is generally not contained
in $[0,\infty)$, therefore the counting sequences are only Hamburger, and not
Stieltjes, moment sequences. The simplest example is that of Motzkin paths,
whose generating function is $\frac{1-x-\sqrt{1-2x-3x^2}}{2x^2} =
1+x+2\,{x}^{2}+4\,{x}^{3}+9\,{x}^{4}+21\,{x}^{5}+51\,{x}^{6}+\cdots$. The
corresponding measure is $\mu(x) = \frac{1}{2\pi}\, \sqrt{(3-x)(1+x)}$ for the
support $\Gamma = [-1,3]$. The principal minor $\Delta_0^n(\seq)$ is equal to
1 for all $n$ but $\Delta_0^n(\seq)$ is $1, 0, -1, -1, 0, 1$ for $n=1,\ldots,
6$, repeating modulo 6 thereafter~\cite[Prop.~2]{Aigner98}, thus confirming
that the Motzkin sequence is not Stieltjes, but only Hamburger. An alternative
way of seeing this is via Flajolet's combinatorial continued
fractions~\cite{Flajolet80}. Indeed, the generating series for Motzkin paths
has a Jacobi continued fraction~\cite[Prop.~5]{Flajolet80} instead of a
Stieltjes continued fraction as required for a Stieltjes moment sequence by
part (d) of Theorem~\ref{thm:Stieltjes}. Some more examples are given in~Appendix~\ref{sec:hamburger}.

\vskip .1cm 

\subsubsection{The second way: exploiting D-finiteness} We can use the
D-finiteness of the generating function $f(z) = \sum_{n \geq 0} C_n z^n$
rather than its algebraicity. We give here the argument, since it will be used
in the subsequent sections for other D-finite (but transcendental) generating
functions. Recall that from the recurrence relation
\[ (n+2) \, C_{n+1} \, - \, (4n+2) \, C_n \, = 0, \qquad \textrm{for 
all } n \geq 0,\]
the generating function $f(z)$ is D-finite and it satisfies the linear differential equation
\[z(4z-1) \, f'(z) + (2z-1) \, f(z) + 1 \, = 0.\]
It follows that $g(z) = f(1/z)/z$ is also D-finite, and satisfies the
differential equation:
\[
z\, \left( {z}-4 \right)  \, g'(z) -2\, g \left( z \right) +1=0.
\]
Recall that by the Stieltjes inversion formula, one has
\[ \mu(x) = -\frac{\psi(x)}{\pi} , \qquad \textrm{where \; } \psi(x) = \lim_{y \rightarrow 0^{+}} \Im \Big(  \, g(x+iy)\Big).\] 
It follows, by taking the limit at $y\to 0^{+}$, by taking the imaginary part, and by linearity, that $\psi(x)$, and thus also the density $\mu(x)$, satisfy the homogeneous part of the previous differential equation, that is:
\begin{equation}\label{eq:deqMuCatalan}
x\, \left( {x}-4 \right)  \, {\mu'  \left( x \right)} -2\,\mu \left( x \right)=0.
\end{equation}
Therefore, $\mu(x)$ is equal, up to a multiplicative constant $\lambda,$ to $\sqrt{ \frac{4-x}{x}}$.
The constant $\lambda$ can be determined using $1 = C_0 = \, \int_0^4 \mu(x) dx = \lambda \,  \int_0^4 \sqrt{ \frac{4-x}{x}} dx  = 2 \pi \,\lambda$.
In conclusion, 
\[ \mu(x) = \, \frac{1}{2\pi} \,  \sqrt{ \frac{4-x}{x} }.\]

\vskip .1cm 

\paragraph{The third way.} There is an alternative way to find $\mu(x)$,
without using the Stieltjes inversion formula. The method is also algorithmic,
close to the one presented in~\cite{BJL19}.

The starting point is, once again, that the generating function 
$f(t) = \sum_{n \geq 0} C_n t^n$ of the Catalan numbers satisfies a linear  differential equation,
 \[t\, (4t-1) \, f'(t) + (2t-1) \, f(t) + 1 = 0.\]
Write $ \, C_n = \, \int_0^4 \, x^n \, \mu(x)\,  dx$. Then,
summation and differentiation imply that
\[f(t) = \int_0^4\frac{1}{1-xt} \, \mu(x) \, dx, \quad \quad 
f'(t) = \int_0^4 \, \frac{x}{(1-xt)^2} \, \mu(x) \, dx.\]
Since $1 = C_0 = \int_0^4 \, \mu(x) \, dx,$
the differential equation for $f$ becomes
\[\int_0^4 \, \left( \frac{t\, (4t-1)\, x}{(1-xt)^2} + \frac{2t-1}{1-xt} + 1\right) \, \mu(x) \, dx \,  = \,  0,\]
which is, after dividing through by $t$ and reducing the pole order by Hermite reduction,
\[\int_0^4 \,  \left ( \frac{2(1-x)}{1-xt} \, \, \,
+ \left (\frac{1}{1-xt} \right)'
\, \,  x \, \, (4-x)  \right) \, \,  \mu(x) \,  dx  = 0.\]
Next, we integrate by parts, to reduce
the derivative $(1/(1-xt))'$ to $1/(1-xt)$:
\begin{align*}  
0  &= \, \left[ \frac{1}{1-xt}  \,  x \, \, (4-x) \,  \mu(x)\right]_{0}^{4} \,
  +\int_0^4 \, \frac{2\, (1-x)}{1-xt}\, \, \mu(x) \, \,
  - \frac{1}{1-xt} \,  \Bigl( x \, \, (4-x) \,  \mu(x)\Bigr)'\\
  &= \; 0 \, + \, 
\int_0^4 \frac{-2\, \mu(x) \,\,  -  (4x-x^2) \,  \mu'(x)}{1-xt} \, \, dx.
\end{align*}
This provides an alternative proof for the differential
equation~\eqref{eq:deqMuCatalan}, and the rest of the argument is similar. In
the next subsection we give a more complex example, corresponding to
pattern-avoiding permutations, which form the main objects studied in this
article.

\subsection{$Av(1342)$} \label{1342}

\begin{figure}[htb]
\centering
\includegraphics[scale =0.44] {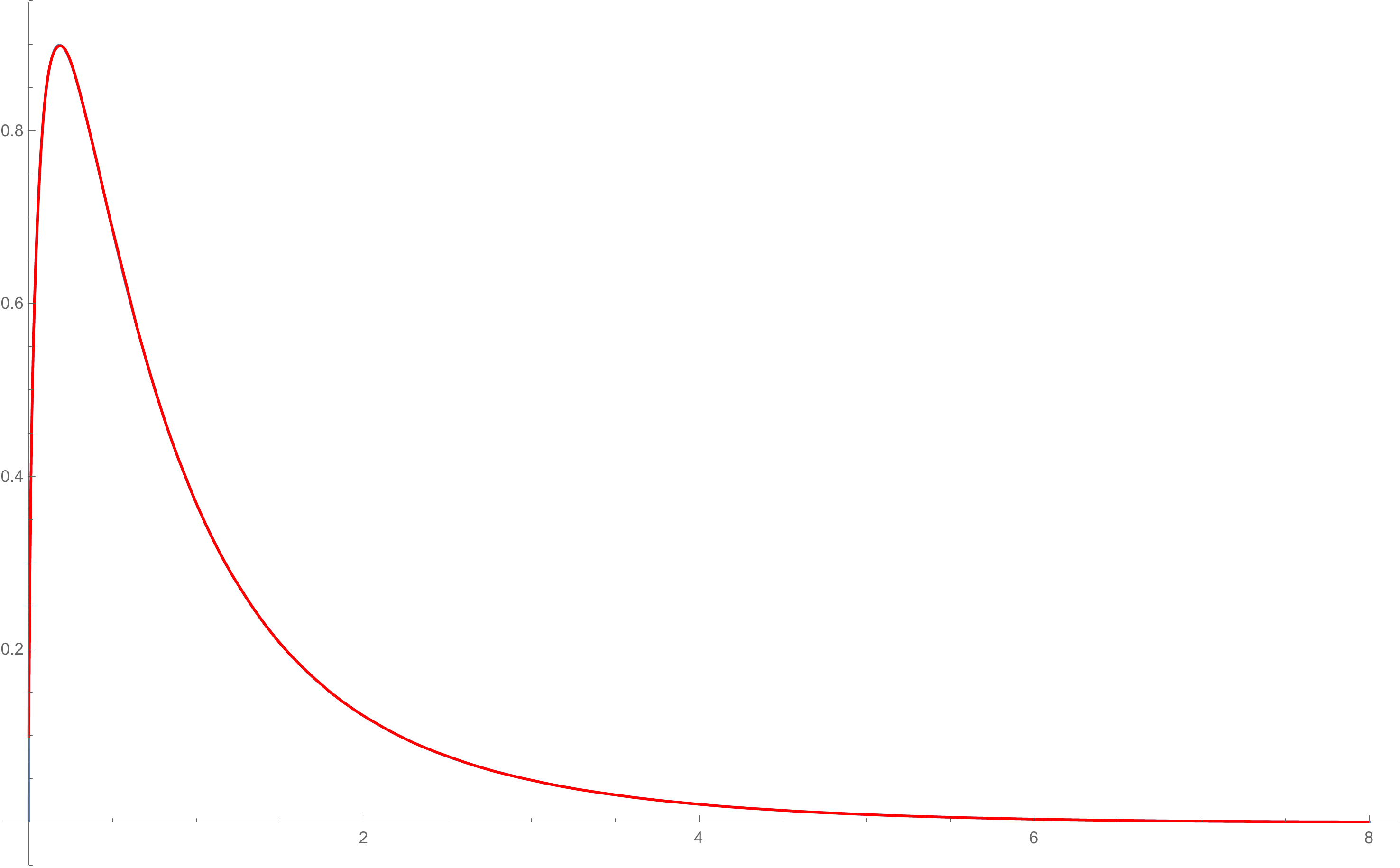}
 \caption{The density function for $Av(1342)$ constructed numerically, and a plot of the exact expression. }
 \label{fig:1342}
\end{figure}

\begin{figure}[htb]
\centering
\includegraphics[scale =0.44] {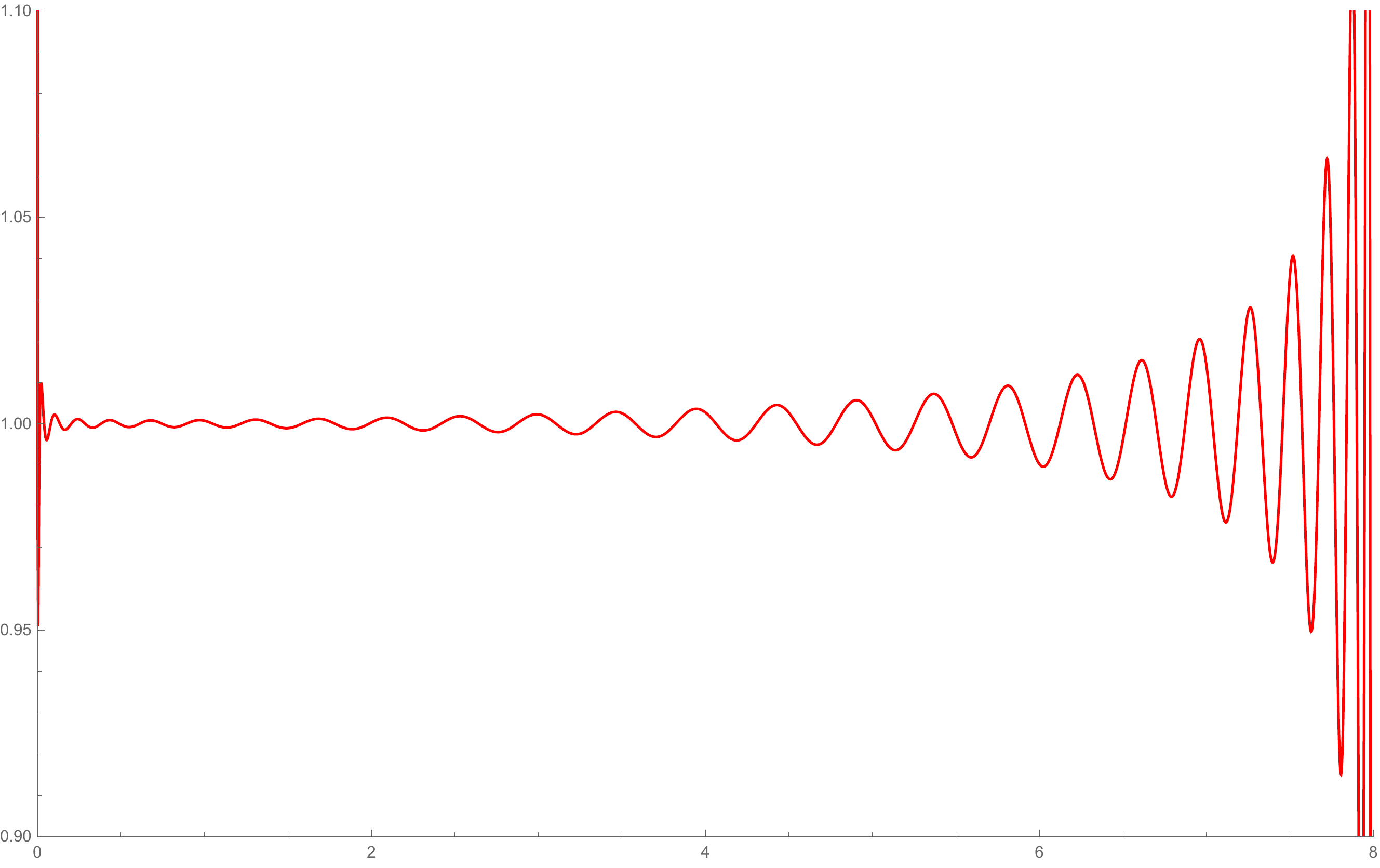}
 \caption{The ratio of the density function for $Av(1342)$ constructed numerically over the exact expression.}
 \label{fig:1342_ratio}
\end{figure}

The generating function of the sequence counting permutations of length $n$
that avoid the pattern $1342$ (\url{https://oeis.org/A022558}), starting
\[1+x+2 x^2+6 x^3+23 x^4+103 x^5+512 x^6+2740 x^7+15485 x^8+\ldots,\] is
algebraic and is equal~\cite{B97} to:
\begin{equation}
f(x) \, = \, \frac{(1-8x)^{3/2}}{2 \, (1+x)^3} \,  \,  + \frac{1+20x-8x^2}{2 \, (1+x)^3}.
\end{equation}
Making the substitution $x \,  \rightarrow \,  \frac{1}{x}$ in $f(x)$
and multiplying by $-\frac{1}{\pi x}$ yields
\[- \frac{1}{\pi x} \, f\left(\frac{1}{x}\right)
\, \, = \,\,  \, \frac{(x-8)^{3/2} \, \sqrt{x}}{2\pi \, (1+x)^3}
\, \,  + \frac{x^2+20x-8}{2\pi \, (1+x)^3}.\]
Since the imaginary part of the rational function is zero, the imaginary part
of the first term {\bf is} the density function. That is to say, the corresponding density is
\[\mu(x)=\frac{(8-x)^{\frac{3}{2}}\sqrt{x}}{2\pi (1+x)^3}.\]
{We plot $\mu(x)$ in Fig.~\ref{fig:1342}, as well as the numerically
constructed density function. To construct the density function numerically,
we make a polynomial approximation $\mu(x) \approx P(x)$ over the known range
$x \in [0,8]$ with the property that the polynomial reproduces the initial
moments of the density function -- that is, just the coefficients of the
sequence $Av(1342),$ -- constrained by $P(8)=P'(8)=0$, see \S\ref{sec:num} for
details. Graphically the two curves are indistinguishable. {In Fig.~\ref{fig:1342_ratio}, we {display} the ratio between these two curves, which shows that the approximation is very good in the bulk of the distribution, but that it becomes less accurate as a ratio as the density approaches 0.}}
                       
\vskip .2cm 
                       
\section{Increasing  $Av(123\cdots (k+1))$ sequences}
\label{Increasing}

We apply the approach from the previous section, notably the Stieltjes
inversion formula, to the explicit study of densities emerging from the
classical ``increasing subsequence problem", or equivalently, from the study
of permutations that avoid an increasing pattern of the form $(123\cdots
(k+1))$.

Later, in Section~\ref{1324}, we address one of the most interesting unsolved
problems in the area of pattern-avoiding permutations, that of enumerating
permutations avoiding the pattern $(1324)$, by using a purely numerical
approach to the Stieltjes inversion.

Let $\pi$ denote a permutation of $\{1,\ldots, n \}$. An increasing
subsequence is a sequence $i_1 < i_2 < \cdots < i_k$ such that $\pi(i_1) <
\pi(i_2)< \ldots <\pi(i_k).$ Let $f_{nk}$ denote the number of permutations
$\pi$ with longest increasing subsequence\footnote{Note that $f_{nk}$ is denoted $u_k(n)$ by
Gessel~\cite[p.~280]{G90} and Stanley~\cite[p.~11]{Stanley}, and $T_k(n)$ by
Bergeron and Gascon~\cite{BG00}.} of length at most $k$.

The number $\, f_{nk}$ is clearly identical to the number of permutations of
$\{1,\ldots, n \}$ that avoid the pattern $(123\ldots (k+1))$; the
set of these permutations is denoted by $Av_{n}(123\ldots (k+1))$, or simply
by $Av(123\ldots (k+1))$, when the size $n$ is implicit.

In random matrix theory, a classical result of Diaconis and
Shashahani~\cite{DS94} states that, if $U(k)$ represents the group of complex unitary matrices of size $k$, then
\[
	E_{M \in U(k)} \left( |\textrm{Tr}(M)|^{2n} \right ) =
n! \qquad \text{for} \;\; n \le k.
\]
Note that, in this context, the expectations are taken with respect to the
Haar measure $\E(f) := \int_{M\in U(k)} f(M) dM$, where $\int$ is the Haar
integral.

This result was extended to the case $n > k$ by Rains~\cite{R98}, who proved 
that 
\begin{equation} \label{eq:Rains} 
E_{M \in U(k)} \left( |\textrm{Tr}(M)|^{2n} \right ) = \, f_{nk}.
\end{equation}
From Rains'{s} result, we can deduce that for any $k$, the counting sequence
$(f_{nk})_{n\in\mathbb{N}}$ of $Av(123\ldots k+1)$ is a Stieltjes moment
sequence. To see this, let $X$ be a random variable with the same distribution
as $|\textrm{Tr}(M)^2|_{M\in U(k)}$. Clearly $X$ is supported on the subset
$[0, \,k^2]$ of $\mathbb{R}_{\geq0}$; indeed, for a $k\times k$ unitary
matrix~$M$, $0\leq|\textrm{Tr}(M)|^2\leq k^2$, as each entry of the matrix has
modulus~$1$.
Rains'{s} {result} implies that $f_{nk}$ is equal to the $n$-th moment $\E(X^{n})$
of the distribution of $X$. Hence $(f_{nk})_{n\in\mathbb{N}}$ is a Stieltjes
moment sequence for any $k$. Moreover, the corresponding density function,
that we will denote by~$\mu_k(x)$, is precisely the density function of
$|\textrm{Tr}(M)^2|_{M\in U(k)}$\footnote{Regev proved in~\cite{R81} that, for
any fixed $n$, the cardinality $f_{nk}$ of $Av_{n}(123\ldots k+1)$ grows
asymptotically like $c_{n,k}\,  k^{2n} \,  n^{\frac{1-k^2}{2}}$ for some
positive constant $c_{n,k}>0$, which in conjunction with
Theorem~\ref{thm:pick} yields another proof that the support of $\mu_k(x)$ is
$[0, \,k^2]$.}.

Equality~\eqref{eq:Rains} yields the following multiple integral representation for $f_{nk}$:
\begin{equation}\label{eq:Rains2}
f_{nk} = \frac{1}{(2 \pi)^k \, k!} \,  \int_{[0,2\,\pi]^k} | e^{i 
\theta_1} + \cdots + e^{i \theta_k}|^{2n} \, \prod_{1\leq p < q \leq k} | 
e^{i \theta_p} - e^{i \theta_q} |^2 \, d \theta_1\ldots d \theta_k.
\end{equation}

In particular, the generating function for $Av(123\ldots (k+1))$ is D-finite,
a result proved on the level of exponential generating functions and without
using~\eqref{eq:Rains2} by Gessel~\cite[p.~280]{G90}\footnote{Note that an
integral identity due to Heine allows to prove, via the exponential generating
function $\sum_n \frac{f_{nk}}{n!} \, t^n$ that~\eqref{eq:Rains2} is
equivalent to Gessel's result, see~\cite[p.~63--65]{Johansson98} and
\cite[p.~1122--1123]{BDJ99}.}. Moreover, Bergeron and Gascon~\cite{BG00}
explicitly computed the differential equations corresponding to the
exponential generating functions for $k\leq 11$ (solving earlier conjectures
made in~\cite{BFK95}). Therefore, the results from the previous section imply
that, for any $k$, the corresponding density~$\mu_k(x)$ is also a D-finite
function.

In subsequent subsections we explicitly find this density. In the simple case
$k = \, 2,$ this corresponds to the Catalan numbers, already treated in
\S\ref{sec:Catalan}. For $k= \, 3,\, 4$ we calculate the density function
$\mu_k(x)$ in terms of (pullbacks of) a rather special subset of $_2F_1$
hypergeometric functions corresponding to {\em classical modular forms.}

Note the formal resemblance of~\eqref{eq:Rains2} with the multiple integral
\begin{equation}\label{eq:BNSW1}
W_{k}(n) = \frac{1}{(2 \pi)^k} \, \int_{[0,2\,\pi]^k} | e^{i 
\theta_1} + \cdots + e^{i \theta_k}|^{2n} \, d \theta_1\ldots d \theta_k.
\end{equation}
which occurs in the theory of uniform random walk integrals in the plane,
where at each step a unit step is taken in a random
direction~\cite{BNSW11,BSWZ12,BSW13}. The integral~\eqref{eq:BNSW1} expresses
the $n$-th moment of the distance to the origin after $n$ steps. Of particular
interest is the evaluation $W_n(1)$, which is equal to the expected distance
after $n$ steps. We will actually see that there is a nontrivial link between
the generating functions of $f_{nk}$ and of $W_{k}(n)$, mirrored by an equally
nontrivial link between the corresponding densities $\mu_k(x)$ and $p_k(x)$.

The relationship  between classical modular forms and some particular $_2F_1$ hypergeometric functions 
is discussed in a number of places in the  literature (see for example \cite{S88}, \cite{M08}, \cite{Zagier08}).
One of the simplest illustrations 
of this intriguing relation can be seen in the identity~\cite[Eq.~(74)]{Zagier08} 
\begin{equation}\label{eq:E4}
E_4(q) =  {_2F_1} \left(\frac{1}{12}, \frac{5}{12}; 1; \frac{1728}{j(q)} \right)^4,
\end{equation}
where $E_4(q)=1 + 240 \, q + 2160 \, q^2 + \cdots $ denotes the weight-four classical Eisenstein series, and where
$j(q) = q^{-1} + 744 + 196884 \, q + 21493760 \, q^2+ \cdots$ is the $j$-invariant (of an elliptic curve). In our examples, the Hauptmodul 
$1728/j$  will be a rational function $\mathcal{H}(x)$ of a variable~$x$.

Modular forms (such as~$E_4$) satisfy infinite order symmetries corresponding
to the modular equations (see~\S\ref{Av1234as}, or Appendix~\ref{Av1234as-App}
and Appendix~\ref{classical2}). The pulled-back hypergeometric functions
representing the classical modular forms\footnote{Only a restricted finite set of $_2F_1$ hypergeometric
functions of the form $_2F_1(a,b; 1; \mathcal{H}(x))$ will yield nomes
corresponding to integer series.} (like ${_2F_1}
\left(\frac{1}{12}, \frac{5}{12}; 1; \mathcal{H}(x)
\right)$) are annihilated by order-two linear
differential operators. The corresponding nome $q(x)$ (explicitly, the
exponential of the ratio of two solutions of the order-two operator) is a
\emph{globally bounded} series, i.e., it can be recast into a series
\emph{with integer coefficients} after a rescaling $x \rightarrow N \, x$.

The emergence of classical modular forms (closely associated with elliptic
curves) in the density functions of $Av(1234)$ and $Av(12345)$ discussed below
(\S\ref{Av1234} and \S\ref{Av12345}), is a hint that the general $Av(123\ldots
(k+1))$ series could be related to integrable theory, an observation which is
in agreement with the relation with random matrix problems, discussed above.

For $k= 5,\,6, \,7$, we construct explicitly a $(k-1)$-th order linear ODE for
the density function $\mu_k(x)$ for $Av(123\ldots (k+1))$, and comment on its
properties (\S\ref{Av123456}, \S\ref{Av1234567} and \S\ref{Av12345678}). These
examples clearly demonstrate the procedure for the cases $k \geq 8.$

\vskip .1cm
\subsection{$Av(1234)$} \label{Av1234}
The sequence $Av(1234)$ admits nice closed forms (\cite[p.~281]{G90},~\cite[p.~12]{Stanley}),
\[
f_{n3} \, = \,
\sum _{k=0}^{n}{\frac {{2\,k\choose k}{n+1\choose k+1}{n+2\choose k+1}}{ \left( n+1 \right) ^{2} \left( n+2 \right) }}
 \, = \,
2\, \sum _{k=0}^{n}{\frac {{2\,k\choose k} {n\choose k} ^{2} \left( 3\,{k}^{2}-2\,kn+2\,k-n+1 \right) }{ \left( 
k+1 \right) ^{2} \left( k+2 \right)  \left( n+1-k \right) }}.
\]
{Its generating function $F_{3}(x)$ starts}\footnote{See 
  \url{https://oeis.org/A005802} (On-line Encyclopedia of Integer Sequences).}
\begin{eqnarray}
\label{A005802}
  \hspace{-0.98in}&&
F_3(x)\,=\,1 \, +x \, +2\,{x}^{2} \, +6\,{x}^{3} \, +23\,{x}^{4} \, +103\,{x}^{5}
\, +513\,{x}^{6} \, +2761\,{x}^{7} 
   \, + \, \, \cdots
\end{eqnarray}
and
 is a solution of an order-three operator   $\, L_3 \, = \, \, L_1 \oplus \, L_2$
 which is the direct sum\footnote{The least common left multiple,
\textsf{LCLM} in Maple's package \textsf{DEtools}.} of an order-one operator $\, L_1$
 with a simple rational solution $\, {\cal R}(x)$
 and an order-two linear differential operator  $\, L_2$ ({recall that} $D_x$ denotes
 the derivati{on} $\, d/dx$):
\begin{eqnarray}
\label{defL2}
\hspace{-0.98in}&&  \quad  \quad    \quad   
L_2 \, \, = \, \, \,\,
D_x^2 \, \,
+ \, \, {\frac {5 \, -27\,x}{ \left(1 \, -9\,x \right) \,  x}} \, \, D_x \, \,
+ \, {\frac {9\,{x}^{2}-9\,x+4}{ x^{2} \, \, (1 \, -9\,x)  \, \, (1 \, - x) }},
\end{eqnarray}
which can be shown to admit, as in~\cite{2DWalks},  a pulled-back $\, _2F_1$ hypergeometric function
$\, {\cal H}(x)$ as solution. Consequently, in a neighborhood of $x=0$, the series
 (\ref{A005802}) can be written as the sum
$F_3(x)\,=\,f_{1}\,:=\,  {\cal R}(x) \, + \, {\cal H}(x)$, where:
\begin{eqnarray}
\label{defL2-bis}
\hspace{-0.98in}&& 
       {\cal R}(x) \, = \, \, {{ 1 \, + \, 5 \, x } \over { 6 \, x^2}}, \\
\hspace{-0.98in}&&       
       {\cal H}(x)  = \,
       - \, {{(1\, -x)^{1/4} \, \, (1\, -9\, x)^{3/4}} \over { 6 \, x^2}}  \, \,
 _2F_1\left( -\, {{1} \over {4}}, \, {{3} \over {4}}; \, 1; 
 \, \, {{ - \, 64 \, x } \over { (1\, -x) \, \, (1 \, -9 \, x)^3}}  \right) \! \! .
\end{eqnarray}
Note that the pulled-back $\, _2F_1$ hypergeometric  $\, {\cal H}(x)$
is {\em not} a Taylor series but a Laurent series. 
The  Taylor series with positive integer coefficients (\ref{A005802})
amounts to getting rid of the poles of the  Laurent series
for  $\, {\cal H}(x)$. Note that the pulled-back $\, _2F_1$ hypergeometric
$\, {\cal H}(x)$ can also be rewritten using the identity
\begin{eqnarray}
\label{Identity1bis}
  \hspace{-0.98in}&&     \quad  
  _2F_1\left( -{{1} \over {4}}, \, {{3} \over {4}} ; \, 1;  \, \, x \right) 
 \, 
  \\
 \hspace{-0.98in}&& \quad   
 \,  \,=  \, \, \,  (1\, -3\, x) \, \,
 _2F_1\left( {{3} \over {4}}, \, {{3} \over {4}}; \, 1;  \, \, x \right)
\, \, \, + \, 4  x  \, (1 \, -x) \, \,
   {{ d \, } \over { dx}}
   \left(\,   _2F_1\left( {{3} \over {4}}, \, {{3} \over {4}} ; \, 1;  \, \, x \right) \right ),
\nonumber 
\end{eqnarray}
as an {\em order-one linear differential operator acting on a classical modular form}
since the $\, _2F_1$
pulled-back hypergeometric function
$ _2F_1\Bigl({{3} \over {4}}, \, {{3} \over {4}}; \, 1; \, \, x \Bigr)$
is actually a classical modular form, {see (B.1)
 in~\cite[Appendix B]{Heun} and Theorem~3 in~\cite{Roques13}, see also~\S\ref{Av1234as}.
Equation~\eqref{involution9x} in Appendix~\ref{Av1234as-App}
also underlies a  % intriguing    	
$\, x \, \leftrightarrow \, \frac{1}{9x}$
involutive symmetry for the pulled-back hypergeometric function~$\, {\cal H}(x)$
and thus also for the $Av(1234)$ series (\ref{A005802})}.

\subsection{Density of $Av(1234)$: Piecewise analysis}
\label{Av1234as}

We recall that, for the sequence $f_{n3}$ for $\, Av(1234)$, the growth rate is $9$~\cite[p.~280--281]{G90}.

The general theory previously displayed tells us that the density can be
obtained from formula (\ref{inversemu}), {\em no matter whether the series we
consider is an algebraic series, a D-finite series} (as in this case $\,
Av(1234)$), {\em a differentially algebraic series, or a
not-even-differentially algebraic series}. Let us rewrite the Stieltjes
inversion formula~(\ref{inversemu}) giving the density as:
\begin{eqnarray}
\label{inversemu2}
  \hspace{-0.98in}&& \quad \quad \quad  \quad 
  \mu_3(y) \, = \,\, - \frac{1}{\pi} \, \,  \Im(\theta(y))
  \quad  \quad \hbox{where:} \quad \quad  \quad 
 \theta(y)  \, = \,\,  {{1} \over {y}} \, \, F_3\Bigl({{1} \over {y}}\Bigr).
\end{eqnarray}

\subsubsection{Naive approach}

In the case of the $\, Av(1234)$ series (\ref{A005802}), 
the function $\, \theta(y)$
is the sum of 
$\, 1/y \cdot \, {\cal R}(1/y)$ and $\, \, 1/y \cdot \,{\cal H}(1/y)$,
which  reads:
\begin{eqnarray}
\label{inversemu2rho}
\hspace{-1.2in}&&  \quad \quad
\theta(y)= {{y+5} \over {6}} \, 
 - \, \, {{ (y-1)^{\frac14} \,  (y-9)^{\frac34}} \over {6}} \,  \,
 _2F_1\left(  -{{1} \over {4}}, \, {{3} \over {4}} ; \, 1;
  \, \, -\, {{64 \, y^3 } \over {(y-1) \, (y-9)^3 }} \right).
\end{eqnarray}
Since {\em  the imaginary part of a polynomial}
(such as $\, {{y+5} \over {6}}$) {\em evaluated at real values is zero},
it is clear that finding the imaginary part of (\ref{inversemu2rho})
 {gives} just the imaginary part of the pulled-back $\, _2F_1$ hypergeometric
function $\, \, 1/y \cdot \,{\cal H}(1/y)$. However evaluating the
imaginary part of a pulled-back $ _2F_1$ hypergeometric is
{\em not} straightforward. This compels us to examine the function as a piecewise
function.

The density $\, \mu_3(y)$ is defined on the interval $\, [0, \, 9]$. We first
note that the pullback in the pulled-back $\, _2F_1$ hypergeometric function
in (\ref{inversemu2rho}), is small when $\, y$ is small. Since we are dealing
with $\, _2F_1$ hypergeometric functions we need to know, for $\, y \in \,[0,
\, 9]$ when the pullback $\, H \, = \, -\, {{64 \, y^3 } \over {(y-1) \,
(y-9)^3 }}$ satisfies $\, H \in \,(-\infty, \, 0]$ or $\, H \in \,[0, \, 1]$
or $\, H \in \,[1, \, \infty)$. This leads one to introduce the following
intervals inside $\, [0, \, 9]$: $[0, \, 1]$, $[1, \, r]$, $[r, \, 9]$, where
$\, r$ is the algebraic number $\, r \, = \, 6\sqrt{3}-9 \, = \,
1.39230485\ldots,$ and is a solution of the quadratic equation $\, y^2 \,
+18\,y \, -27 \, = \,0$ corresponding to $\, H \, = \, 1.$

\vskip .1cm

The detailed piecewise analysis of this function and its imaginary part is
given in Appendix~\ref{app:maple}. The conclusion of this naive approach is
that we find an expression which is the correct density, but divided by a
factor of 3.

\subsubsection{The magic of modularity}
Actually, the origin of the factor 3 comes from the modular equation~\eqref{Mod13}. More precisely,
the following identity $\mathcal{H}(x) = \mathcal{H}_2(x)$ holds:
\begin{eqnarray}
\label{pullbackL21over9xidcalHbis}
\hspace{-0.48in}
 - \, {{(1\, -x)^{1/4} \, \, (1\, -9\, x)^{3/4}} \over { 6 \, x^2}}  \, \,
 _2F_1\left(-\, {{1} \over {4}}, \, {{3} \over {4}};  \, 1; 
 \, \, {{ - \, 64 \, x } \over { (1\, -x) \, \, (1 \, -9 \, x)^3}}  \right)
 \nonumber  \\
\hspace{-0.48in}
= -\frac{(1-9x)^{1/4}\, (1+3x)^2}{6\, (1-x)^{5/4}\, x^2} \, \, \,  _2F_1\left(-\, {{1} \over {4}}, \, {{3} \over {4}};  \, 1; 
 \, \, {{ - \, 64 \, x^3 } \over { (1\, -x)^3 \, \, (1 \, -9 \, x)}}  \right)\nonumber \\
\hspace{-0.38in}
 +\frac{32\, x^3\, (1-6x-3x^2)}{27}\, \, \, _2F_1\left(-\, {{3} \over {4}}, \, {{7} \over {4}};  \, 1; 
 \, \, {{ - \, 64 \, x^3 } \over { (1\, -x)^3 \, \, (1 \, -9 \, x)}}  \right),
\end{eqnarray}
where $\mathcal{H}_2$ is the right-hand side of~\eqref{pullbackL21over9xidcalHbis}.

This expresses the fact that {the} principal part of the generating function of
$Av(1234)$ corresponds to an order-one differential operator acting on a
classical modular form.

Let us introduce $\mathcal{G}(y):=\frac{1}{y} \mathcal{H}(1/y)$ and
$\mathcal{G}_2(y):=\frac{1}{y} \mathcal{H}_2(1/y)$. Naively, one could imagine
that $\mathcal{G}(y)$ and $\mathcal{G}_2(y)$ are equal for small~$y$ and
possibly in some interval like $[0,9]$. In fact, this is not true: their
series expansions around $y=0$ satisfy $3 \, \mathcal{G}(y) =
\mathcal{G}_2(y)$. In the Stieltjes inversion formula, the density candidate
gets multiplied by a factor of 3. Thus, the corresponding density becomes the
correct one.

\subsubsection{Less naive approach} 
\begin{figure}[htb]
\centering
\includegraphics[scale =0.44] {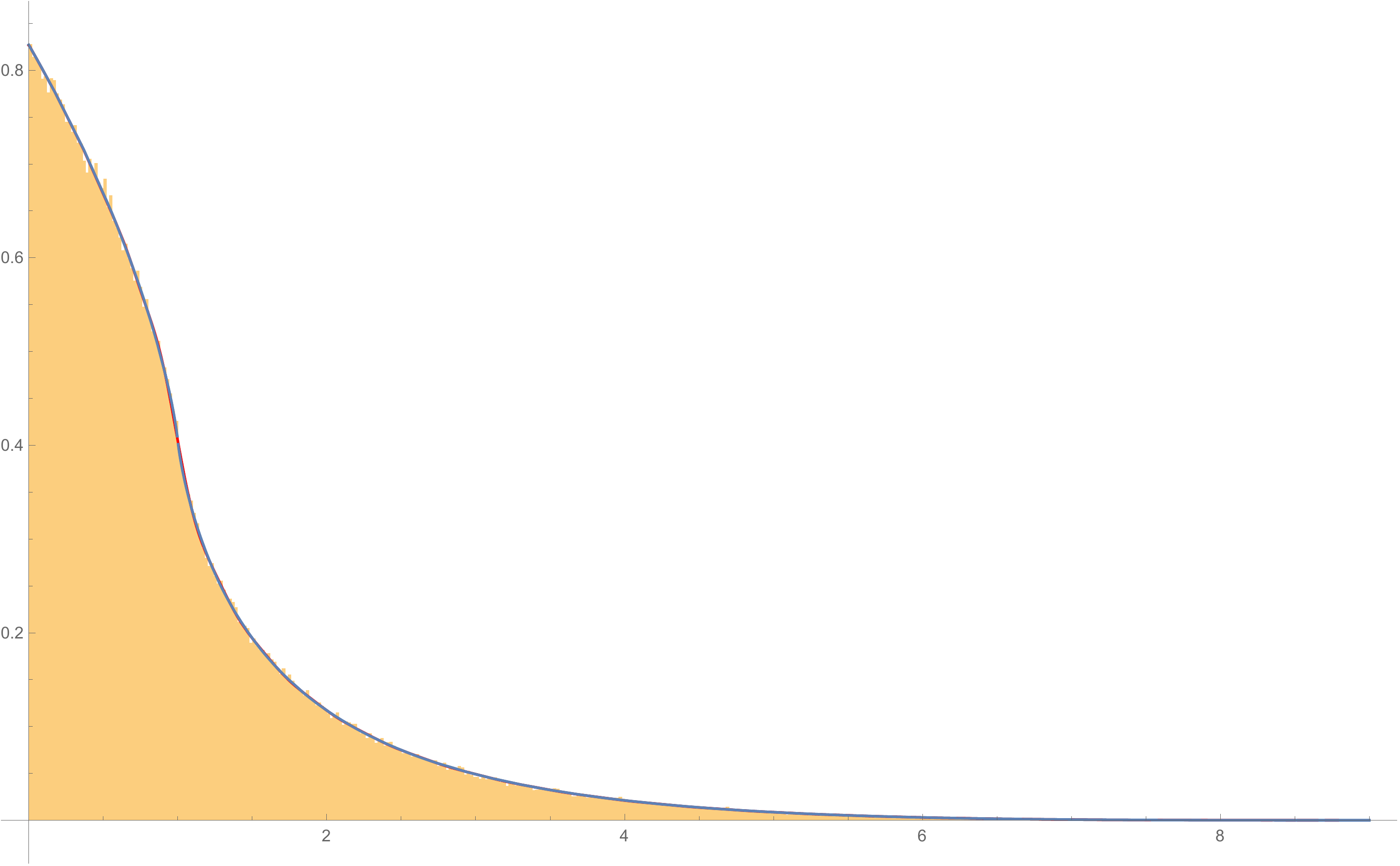}
 \caption{The density function $\mu_3(x)$ for $Av(1234)$ constructed numerically by a histogram (orange bars) and the plot created from the exact expression (blue curve). The plot also contains a polynomial approximation (red curve), but this is perfectly covered by the exact plot.}
 \label{fig:1234}
\end{figure}

Again, we need to perform a piecewise
analysis of the slightly more complicated expression $\mathcal{H}_2(y)$, as we
did with $\mathcal{H}(y)$ in Appendix~\ref{app:maple}. Not surprisingly, after
performing this analysis, one gets the correct density. In
Appendix~\ref{app:H2} we sketch some arguments showing that this alternative
expression~$\mathcal{H}_2(y)$ behaves, as expected, better than the simpler
one using $\mathcal{H}(y)$. {In Fig.~\ref{fig:1234} we plot this correct
density function, as well as plotting the density function, numerically
constructed as described in Sec. \ref{sec:num}. Again, they are graphically
indistinguishable. }

 \vskip .3cm 

\subsection{$Av(12345)$}
\label{Av12345}

The positive integer coefficient series $F_4(x)$ for $Av(12345)$ reads\footnote{See
  \url{https://oeis.org/A047889}.}:
\begin{eqnarray}
\label{A047889}
  \hspace{-0.98in}&& 
  1 \, +x \, +2\,{x}^{2} \, +6\,{x}^{3} \, +24\,{x}^{4}
  \, +119\,{x}^{5} \, +694\,{x}^{6}
+4582\,{x}^{7}+33324\,{x}^{8} \, + \, \cdots 
\end{eqnarray}
It is {a} solution of an order-four linear differential operator
$\, L_4 \, = \, \, L_1 \, \oplus \, L_3$
 which is the direct sum 
of an order-one linear differential operator~$\, L_1$
 with a simple rational solution $\, {\cal R}(x)=\frac{ 18\,{x}^{2}+10\,x+1}{12 \,  {x}^{3}}$
 and an order-three linear differential operator~$\, L_3$
\begin{eqnarray}
  \label{L3def}
  \hspace{-0.98in}&& \quad  \quad 
L_3 \, \, = \, \, D_x^3 \,\,
+ {\frac { \left( 448\,{x}^{2}-182\,x+13 \right)}{
x \,  \, (1 \, - 16\,x) \,  \, (1 \, -4\,x) }} \,  \, D_x^2
 \nonumber \\
  \hspace{-0.98in}&& \quad  \quad  \quad   \, \, \,\,\,\,
+\,{\frac { 4 \,  \, \left( 160\,{x}^{2}-87\,x+11 \right) }{{x}^{2}
\,  \, (1 \, - 16\,x) \,  \, (1 \, -4\,x)}} \,  \, D_x
\, \,\, +\,{\frac { 4 \,  \, (32\,{x}^{2}-20\,x+9)}{
    (1 \, - 16\,x) \,  \, (1 \, -4\,x)  \,  \, {x}^{3}}}, 
\end{eqnarray}
 which can be seen, using van Hoeij's algorithm in~\cite{Hoeij}\footnote{The implementation of this algorithm (\textsf{ReduceOrder}) is available at \url{https://www.math.fsu.edu/~hoeij/files/ReduceOrder/ReduceOrder}}, to be homomorphic
 (with order-two intertwiners) to the
 symmetric square of an order-two  linear differential operator:
\begin{eqnarray}
  \label{U2def}
  \hspace{-0.98in}&& \, 
  U_2 \, \, = \, \, \,
D_x^2  \, \,  \, +  \,{\frac { 2 \,  \, (32\,x \, -5) }{
(1 \,  -16\,x)  \,  \, (1 \, -4\,x) }} \,  \, D_x
\, \,  \, + \,{\frac {1\, -8\,x}{
    4 \,  \, (1 \, - 16\,x) \,  \, (1 \, -4\,x) \,  \,  x^{2}}}, 
\end{eqnarray}
which has a classical modular form solution that can be written as a
pulled-back hypergeometric function. In Appendix~\ref{classical2} we show that this
pulled-back hypergeometric function is a classical modular form.
Appendix~\ref{classical2} also underlines another involutive symmetry $\, x \,
\leftrightarrow \, \frac{1}{64x}$ on this new classical modular form and thus
on the $Av(12345)$ series (\ref{A047889}). 

Introducing $\, {\cal F}(x) \, = \, {\cal H}(x)^2$, the square of this classical
modular form 
\begin{eqnarray} \label{classmodformsquare} 
	\hspace{-0.98in}&&
\quad \quad \quad \, \, \, {\cal F}(x) \, \, = \, \, \, \, x \,  \, (1\, -16 \,
x)^{-1/3} \,  \, (1\, +2 \, x)^{-2/3} \nonumber \\ \hspace{-0.98in}&& \quad
\quad \quad \quad \quad \quad \quad \times \, \, _2F_1\left( {{1} \over {6}},
\, {{2} \over {3}};  \, 1;  \, \, - \,{\frac { 108\,{x}^{2}}{ (1\, - 16\,x)
\,  \, (1 \, +2\,x)^{2}}} \right)^2, 
\end{eqnarray} 
one finds that the series $F_4(x)$ for $\,Av(12345)$ can be written as the sum of the
rational function $\, {\cal R}(x)$ (the solution of $\, L_1$) and an order-two
linear differential operator acting on $\, {\cal F}(x)$ given by
(\ref{classmodformsquare}). Thus, the $\, Av(12345)$ series
(\ref{A047889}) is actually the expansion of 
\begin{eqnarray} \label{ordertwoinvolved} 
 	\hspace{-0.88in}&&  \quad -{\frac { \left(
4\,x-1 \right) \left( 16\,x-1 \right) \left( 1376\,{x}^{2}+590\,x-1 \right)
}{864\,{x}^{3}}} \,  \, {{ d^2 {\cal F}(x)} \over {dx^2}} \nonumber \\
\hspace{-0.98in}&& \quad \quad \quad \quad -{\frac { \left(
63488\,{x}^{4}+10304\,{x}^{3}+396\,{x}^{2}-344\,x+1 \right)}{ 864\,{x}^{4}}}
\,  \, {{ d {\cal F}(x)} \over {d{x}}} \\ \hspace{-0.98in}&& \quad \quad \quad \quad
\quad \, \, -{\frac { \left( 2176\,{x}^{3}-3372\,{x}^{2}+420\,x-1 \right) }{
864\,{x}^{5}}} \,  \, {\cal F}(x) \, \, \, \, \, + {\frac {18\,{x}^{2}+10\,x+1}{12
\,\, {x}^{3}}}. \nonumber 
\end{eqnarray}

Alternatively, $L_3$ can also be seen as being homomorphic (with order-two
intertwiners) to the \emph{symmetric square} of another order two linear differential
operator, with a pulled-back hypergeometric solution:
\[
_2F_1\left({{1} \over {6}},
\, {{1} \over {3}};  \, 1;  \, \,  \,{\frac { 108\,{x}^{2}}{ (1\, - 4\,x)^3}} \right),
\]
which emerges in the study of Domb numbers
 (counting $2n$-step polygons on the diamond lattice, and also the
moments of a 4-step random walk in two dimensions,
cf. \url{https://oeis.org/A002895}). In other words,
the following identity holds:
\begin{eqnarray}
 \left( 1-4\,x \right) ^{\frac12}
\,  \,
_2F_1\left({{1} \over {6}},\, {{2} \over {3}};  \, 1;  \, \,  \,{\frac { -108\,{x}^{2}}{ (1\, - 16\,x) \,  (1+2\,x)^2}} \right)
\nonumber \\ =
{ \left( 1-16\,x\right)^{\frac16}  \,  \left( 1+2\,x \right) ^{\frac13} }
\,  \,
_2F_1\left({{1} \over {6}},\, {{1} \over {3}};  \, 1;  \, \,  \,{\frac { 108\,{x}^{2}}{ (1\, - 4\,x)^3}} \right).
\end{eqnarray}

As we will see in the next subsection, the relation between the generating
function~$F_4(x)$ of $Av(12345)$ and the Domb generating function is not
fortuitous: it mirrors what happens on the level of the density functions,
namely that the density function~$\mu_4(x)$ corresponding to $Av(12345)$ is
intimately related to the density function~$p_4(x)$ of the distance travelled
in 4 steps by a uniform random walk in the plane~\cite{BSWZ12}.

\vskip .2cm

\subsection{Density of $Av(12345)$}
\label{Av12345as}

The $Av(12345)$ density function $\mu_4(x)$ is also the probability density function of
the squared norm of the trace of a random 4 by 4 unitary matrix. We recall
that the support of the measure is $\, [0, \, 16]$.
The sequence (\ref{A047889}) is a Stieltjes sequence, the density being given
by the general formula~\eqref{inversemu2}. Again one has to change $\, x \,
\rightarrow \, 1/x$ in~\eqref{ordertwoinvolved} and take the imaginary part of
the slightly involved hypergeometric result.

\begin{figure}[htb]
\centering
\includegraphics[scale =0.44] {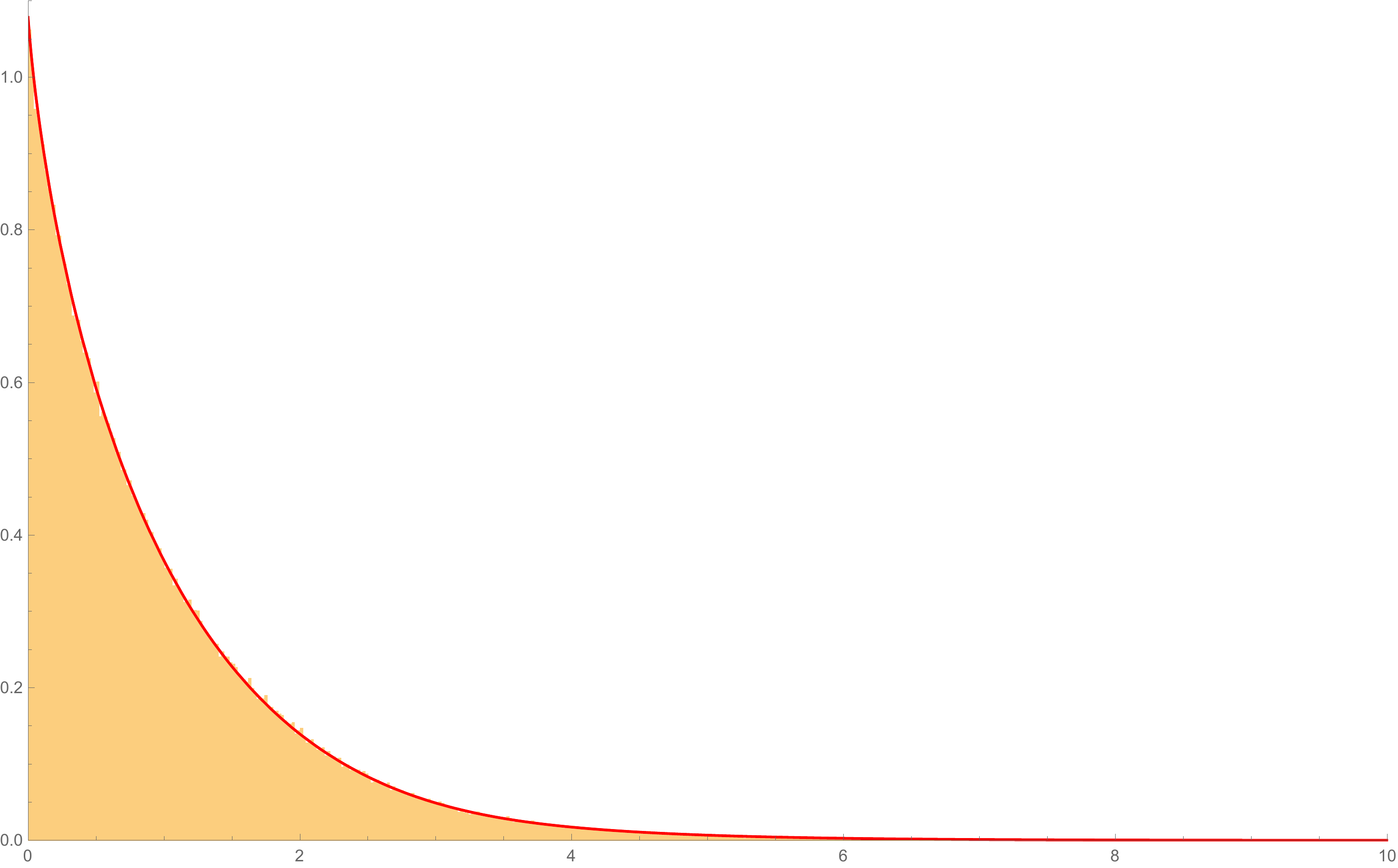}
 \caption{The density function $\mu_4(x)$ for $Av(12345)$ constructed numerically by a histogram (orange bars) and the polynomial approximation (red line).}
 \label{fig:12345}
\end{figure}

As far as finding the density $\mu_4(x)$ is concerned (before undertaking the
delicate piecewise analysis of taking the imaginary part) one finds, using~\eqref{ordertwoinvolved} and Stieltjes inversion formula, that
\[\theta(x) = \frac{1}{x} \,  \left( F_4 \left(\frac{1}{x} \right) - \frac{1}{12}\, \left( {x}^{2}+10\,x+18 \right) x\right)\]
is equal to
\begin{eqnarray}
\label{density}
\hspace{-0.98in}&& \quad \quad \quad
    \theta(x)=   x^{2} \,  \, (4-x)  \,  \, (16-x)
       \,  \, (x^{2}-590\,x-1376) \,  \, {{d^2 M(x)} \over { dx^2}}
\nonumber \\
\hspace{-0.98in}&& \quad \quad \quad \quad \quad
\, \, + \, x \,  \, (3\,x^{4}-1564\,x^{3}
+21372\,x^{2}-10176\,x-112640) \,  \, {{dM(x)} \over { dx}}
\nonumber \\
\hspace{-0.98in}&& \quad \quad \quad \quad \quad
\, \,  +x  \,  \, \left( x^{3}-420\,x^{2}+3372\,x-2176 \right) \,  \,  M(x), 
\end{eqnarray}
where:
\begin{eqnarray}
\label{density2}
\hspace{-0.98in}&& \quad
M(x) \, \, = \, \, \, (16\, -x)^{-1/3} \,  \,  (2\, +x)^{-2/3} \,   \,
 _2F_1\left({{1} \over {6}}, \, {{2} \over {3}}; \, 1; \, \,
\,\,{\frac {  108\, x}{ (16\, - \,x) \,  \, (2 \, +\,x)^{2} }}   \right)^2.
\nonumber 
\end{eqnarray} 
 In this case the support of the measure $\mu_4(x)$ is $\, [0,\, 16]$
 and the piecewise analysis will need to consider the intervals $\, [0,\, 4]$
 and $\, [4,\, 16]$. The value $\, 4$ corresponds to the pullback
 $\, {\frac {  108\, y}{ (16\, - \,y) \,  \, (2 \, +\,y)^{2} }}$
 being equal to $\, 1$. 
We will not pursue these calculations in detail here, as they are similar to
those in the previous section for $Av(1234)$. We simply show in
Fig.~\ref{fig:12345} a plot of the density function $\mu_4(x)$ for $Av(12345)$
obtained numerically, using the approach described in \S\ref{sec:num}.

\vskip .1cm

The linear differential operator annihilating the density $\, \mu_4(x)$ is
easily obtained from the pullback of the order-three operator $L_3$ in~\eqref{L3def}. It reads\footnote{Note that, starting from this operator, and from the initial conditions that uniquely
specify the density in the space of its solutions, Mezzarobba's
algorithms~\cite{Mezzarobba} could in principle obtain (numerically) this
density, by effective (numerical) analytic continuation.}:
\[
\left( x-4 \right)  \left( x-16 \right) \, {x}^{2}D_x^{3} 
-2\,x \left( 2\,{x}^{2}-x-64
 \right) D_x^{2}+10\,x \left( x+2 \right) D_x-12\,(x-2).
\]
Similarly, $\mu_4(x^2)$ is a solution of the linear differential operator
\begin{align}\label{eq:mux2}
\left( x^2-4 \right) \,  \left( x^2-16 \right) 
{x}^{2} D_x^{3}- \left( 11\,{x}^{4}-64\,{x}^{2}-64 \right) x D_x^{2}\nonumber
\\+ \left( 51\,{x}^{4}+16\,{x
}^{2}-64 \right) D_x-96\,x\,  \left( {x}^{2}+2 \right).
\end{align}

\subsection{Relation between densities of $Av(12345)$ and of $4$-steps uniform random walks} \label{Av12345as:rel}

An unexpectedly simple relation connects the density $\mu_4(x)$ for
$Av(12345)$ and the density $p_4(x)$ of short uniform random walks studied
in~\cite{BSWZ12}.

Recall that $p_4(x)$ denotes the density function of the distance travelled in~$4$ steps by a uniform random walk in $\mathbb{Z}^2$
that starts at $(0,0)$ and consists of 
unit steps, each step being taken into a uniformly random direction.
In~\cite[Thm.~4.7]{BSWZ12} it was show{n} that for $x\in (2,4)$, 
\[
p_4(x) = \frac{2}{\pi^2} \, \frac{\sqrt{4-x^2}}{x}
\,  
\,
_3F_2\left(\frac12,\, \frac12, \, \frac12;, \,  \frac56, \, \frac76 ; \, \,  \,{\frac { (16-x^2)^3}{ 108 \,  x^4}} \right),
\]
which satisfies the order-3 operator~\cite[Eq.~(2.7)]{BSWZ12}
\begin{equation*}\label{eq:p4}
A_4 = 	 \left( x^2-4 \right)  \left( x^2-16 \right) {x}^{3}  D_x^{3}+6\,{x}^{4} \left( {x}^{2}-10 \right) D_x^{2}+x \left( 7\,{x}^{4}-32\,{x}^{2}+
64 \right) D_x+ \left( {x}^{4}-64 \right) .
\end{equation*}

Interestingly, the operator~\eqref{eq:mux2} satisfied by $\mu(x^2)$ 
is  gauge equivalent (homomorphic) with $A_4$, 
with a nontrivial order-2 intertwiner
\begin{align}\label{eq:p4-to-mu}
x \left( x+4 \right)  \left( x-2 \right)  \left( x+2 \right)  \left( x-4 \right)  \left( {x}^{
4}-590\,{x}^{2}-1376 \right) D_x^{2} \nonumber \\
+ \left( 3\,{x}^{8}-1298\,{x}^{6}+11280\,{x}^{4}+
10368\,{x}^{2}+38912 \right) D_x \nonumber  \\
+{\frac {{x}^{8}-382\,{x}^{6}+2208\,{x}^{4}-19072\,{x}^{2
}-38912}{x}}.
\end{align}
This shows that the density $p_4(x)$ and our density $\mu_4(x^2)$, both with
support $[0,4]$, are actually \emph{differentially related in a non-trivial
way}. More precisely, the operator~\eqref{eq:p4-to-mu} sends~$p_4(x)$ to
a solution of the operator~\eqref{eq:mux2} satisfied by~$\mu_4(x^2)$. 
Conversely, $\mu_4(x^2)$ is mapped by a similar order-2 operator to a solution
of the operator satisfied by~$p_4(x)$.

We can push a bit further the analysis made in~\cite{BSWZ12}.
The operator $A_4$ annihilating $p_4(x)$ is the symmetric square of
\[
4   \left( x^2-4 \right)  \left( x^2-
16 \right) \,{x}^2 D_x^{2} +8\,{x}^{3} \left( {x}^{2}-10 \right) D_x+({x}^{4}-12\,{x}^{2}+64),
\]
which has a basis of solutions $\left\{ y_1(x) = A(x)\,  H_1(x), \; y_2(x) =  A(x)\,  H_2(x)\right\}$ with
\[
H_1(x) = \;
_2F_1\left({{1} \over {6}},\, {{2} \over {3}};  \, 1;  \, \, 
\,\phi(x)\right),\quad
H_2(x)= \;
_2F_1\left({{1} \over {6}},\, {{2} \over {3}};  \, \frac56 ;  \, \, 
1 - \,\phi(x)\right),
\]
where 
\[A(x) = \left(\frac { x^3}{ (x^2+2)^2 \,  (16-x^2)} \right)^{\frac16} \qquad
\text{and} \qquad \phi(x) = {\frac {108\, {x}^{2}}{ \left( 16-x^2 \right)    \left( {x}^{2}+2 \right) ^{2}}}.\] 
Thus, $p_4(x)$ is equal to
$c_1 \, y_1(x)^2 + c_2 \, y_2(x)^2 + c_3 \, y_1(x) \,  y_2(x),$
where
\[
c_1 = {\frac {\sqrt {6} }{2\, \pi}},
\quad\;
c_2=
\frac23 \,{\frac {\sqrt [6]{2}\, \pi\, \left(i- \sqrt {3} \right) }{  \Gamma \left( 2/3 \right) ^{6}}},
\quad\;
c_3= - \frac{2 \,i}{3} {\frac {{2}^{5/6}\sqrt {3}}{ \Gamma \left( 2/3 \right)  ^{3}}}.
\]
Similarly, an explicit differential expression of order 2 in $\mu_4(x^2)$ could also  be written  as a linear combination with constant coefficients of $y_1^2, y_2^2$ and $y_1 \,  y_2$.

\subsection{$Av(123456)$: order-four operator homomorphic to its adjoint}
\label{Av123456}

The positive integer coefficient series $F_5(x)$ for $Av(123456)$ reads\footnote{See 
  \url{https://oeis.org/A047890}.}:
\begin{eqnarray}
\label{??}
  \hspace{-0.98in}&& 
  1 \, +x \, +2\,{x}^{2} \, +6\,{x}^{3} \,
  +24\,{x}^{4} \, +120\,{x}^{5} \, +719\, x^6 \, +5003 \, x^7
  \, +39429 \, x^8 \, + \, \cdots.
\end{eqnarray}
It is {a} solution of an order-five linear differential operator
$\, L_5 \, = \, \, L_1 \, \oplus \, L_4$
 which is the direct sum of an order-one operator $\, L_1$
 with a simple rational solution $\, {\cal R}(x)$
 \begin{eqnarray}
\label{L4rat}
  \hspace{-0.98in}&& \quad \quad  \quad \quad  \quad \quad 
         {\cal R}(x) \, \, = \, \, \,  \, \,
         {\frac {  1\, +17\,x \, +71\,{x}^{2}  + 63\,{x}^{3}   }{{20 \, \, x}^{4}}}, 
\end{eqnarray}
 and an order-four linear differential operator  $\, L_4.$
This is
 \begin{eqnarray} \label{operL4}
  \hspace{-0.98in}&& 
 L_4= (x-1)\,  (9 \,x-1)\, (25 \, x-1)\,  x^4 D_x^4+2x^3(1350x^3-1865x^2+336x-13)D_x^3 \nonumber \\
   \hspace{-0.98in}&& 
+ x^2(8550x^3-14241x^2+3650x-215)D_x^2+x(7200x^3-14643x^2+5708x-649)D_x  \nonumber \\
  \hspace{-0.98in}&& 
 +(900x^3-2325x^2+1225x-576).
 \end{eqnarray}
 The linear differential operator annihilating the corresponding density $\,
\mu_5(x)$ can easily be obtained from the pullback of the order-four
operator~\eqref{operL4}. Since we do not have any closed exact formula for its
solutions\footnote{No $ _4F_3$ expressions, no reduction to one variable
of an Appel-Lauricella function.} we cannot perform the kind of piecewise
analysis we performed previously. Mezzarobba's algorithms~\cite{Mezzarobba}
could in principle obtain (numerically) this density, by effective (numerical)
analytic continuation.

The order-four operator $\, L_4$ is MUM (``maximal unipotent monodromy'') and
is non-trivially homomorphic to its adjoint with an order-two intertwiner. Its
\emph{exterior square} has a rational solution of the form $\, P(x)/(301 \,
x^{11})$, where $\, P(x)$ is a polynomial with integer coefficients. The
differential Galois group of the order-four operator $\, L_4$ is thus $\,
Sp(4, \, \mathbb{C})$, or a special subgroup of it~\cite{BHMW15}\footnote{Note that in contrast to the two previous cases, where intriguing involutions
occur, we do not have a homomorphism between the order-four linear
differential operator $\, L_4$ and its pulled-back operator by an involution
$\, x \, \leftrightarrow \, A/x$.}. 

\smallskip 
{\bf Remark.} A priori, we cannot exclude the fact that $L_4$ could be
homomorphic to the symmetric cube of a second-order linear differential
operator, or to a symmetric product of two second-order operators.
Furthermore, it could also be, in principle, that these second-order operators
admit classical modular forms as solutions (pullbacks of special $_2F_1$
hypergeometric functions). However, these options can both be excluded by
using some results from differential Galois theory~\cite{Singer}, specifically
from \cite[Prop.~7, p.~50]{Person} for the symmetric cube case, and from
\cite[Prop.~10, p.~69]{Person} for the symmetric product case, see
also~\cite[\S3]{Hoeij}. Indeed, if $L_4$ were either a symmetric cube or a
symmetric product of order-two operators, then its symmetric square would
contain a (direct) factor of order~3 or~1. This is ruled out by a
factorization procedure which shows that the symmetric square of~$L_4$ is
(LCLM-)irreducible.

Still, we cannot exclude the fact that the solutions of $L_4$ (and in
particular the generating function for $Av(123456)$, minus the rational
part~\eqref{L4rat}) could be written as an algebraic pullback of a $_4F_3$
hypergeometric function.

\begin{figure}[htb]
\centering
\includegraphics[scale =0.44] {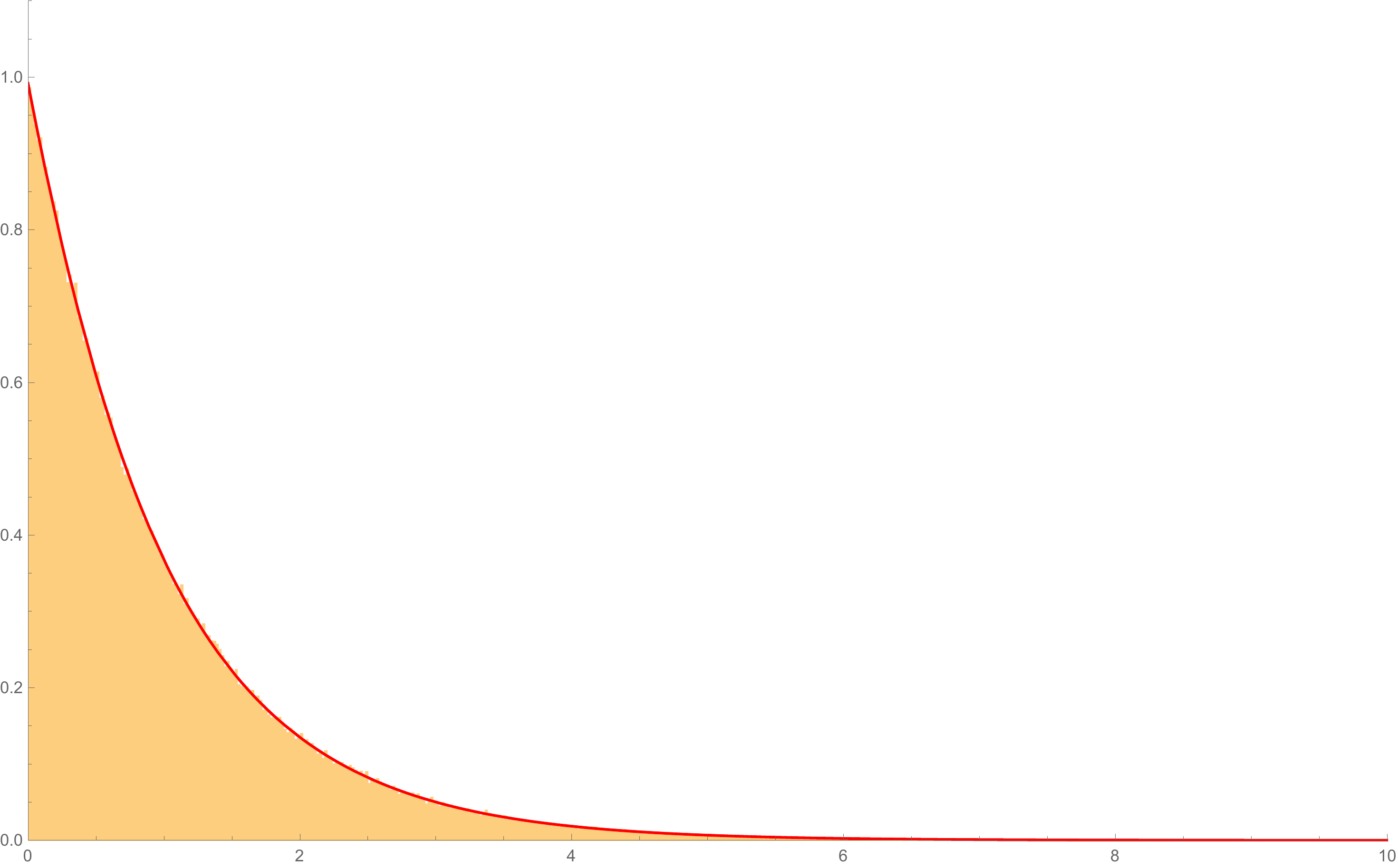}
 \caption{The density function $\mu_5(x)$ for $Av(123456)$ constructed numerically by a histogram (orange bars) and the polynomial approximation (red line).}
 \label{fig:123456}
\end{figure}

{We show in Fig.~\ref{fig:123456} a plot of the density function $\mu_5(x)$
for $Av(123456)$ obtained numerically, using the approaches described in
\S\ref{sec:num}. For $k>6$ we find that our approximations are visually
indistinguishable from the graph of $e^{-x}$, so we do not display these plots
separately.}

\vskip .1cm

\subsection{Relation between densities of $Av(123456)$ and of 5-steps uniform random walks} \label{Av123456as:rel}

An unexpectedly simple relation connects the density $\mu_5(x)$ for
$Av(123456)$ and the density $p_5(x)$ of short uniform random walks studied
in~\cite{BSWZ12}.

Recall that $p_5(x)$ denotes the density function of the distance travelled in~$5$ steps by a uniform random walk in $\mathbb{Z}^2$
that starts at $(0,0)$ and consists of unit steps, each step being taken into a uniformly random direction.
In~\cite[Thm.~5.2]{BSWZ12} it was show{n} that for $x\in (0,5) \setminus \{ 1, 3 \}$, the density $p_5(x)$ satisfies
the order-4 operator~\cite[Eq.~(5.4)]{BSWZ12}
\begin{align*}
A_5 = 	
\left( x-1 \right)  \left( x-3 \right)  \left( x-5 \right)  \left( x+5 \right)  \left( x+3 \right)  \left( x+1
 \right)  {x}^{4}  D_x^{4}
\\+2\, \left( 5\,{x}^{6}-105\,{x}^{4}+259\,{x}^{2}+225 \right) \, {x}^{3}  D_x^{3}
\\+ \left( 25\,{
x}^{6}-287\,{x}^{4}+363\,{x}^{2}-1125 \right) \, {x}^{2} D_x^{2}
\\+  \left( 15\,{x}^{6}-77\,{x}^{4}-363\,{x}^{2}+2025 \right) \, x D_x\\+ \left({x}^{6}-3\,{x}^{4}+363\,{x}^{2}-2025 \right).
\end{align*}
Interestingly, the operator satisfied by our density $\mu_5(x^2)$,
\begin{align} 
-{x}^{3} \left( x-1 \right)  \left( x-3 \right)  \left( x-5 \right)  \left( x+5 \right)  \left( x+3 \right)  \left( x+1
 \right) D_x^{4} \nonumber
\\+2\,{x}^{2} \left( 13\,{x}^{6}-217\,{x}^{4}+363\,{x}^{2}+225 \right) D_x^{3}
\nonumber
\\-x \left( 287\,{x}^{
6}-1685\,{x}^{4}-1323\,{x}^{2}+225 \right) D_x^{2}\nonumber
\\+ \left( 1575\,{x}^{6}-437\,{x}^{4}-1683\,{x}^{2}+225 \right) D_x\nonumber
\\
-1200\,{x}^{3} \,  \left( 3\,{x}^{2}+5 \right), \label{eqmux2ter}
\end{align}
is  gauge equivalent (homomorphic) with $A_5$, 
with a nontrivial order-3 intertwiner.
This shows that the density $p_5(x)$ and our density $\mu_5(x^2)$, both with
support $[0,5]$, are actually \emph{differentially related in a non-trivial
way}. 
More precisely, the intertwiner sends~$p_5(x)$ to
a solution of the operator~\eqref{eqmux2ter} satisfied by~$\mu_5(x^2)$. 
Conversely, $\mu_5(x^2)$ is mapped by a similar order-2 operator to a solution
of the operator satisfied by~$p_5(x)$.

%%%%
\vskip .1cm

\subsection{$Av(1234567)$: order-five operator homomorphic to its adjoint}
\label{Av1234567}

The positive integer coefficient series $F_6(x)$ for $Av(1234567)$ reads\footnote{See \url{https://oeis.org/A052399}.}:
\begin{eqnarray}
\label{A052399}
  \hspace{-0.98in}&& 
  1 \, +x \, +2\,{x}^{2} \, +6\,{x}^{3} \, +24\,{x}^{4} \, +120\,{x}^{5}
  \, +720\,{x}^{6} \,\, +5039 \,{x}^{7}  \, + \,40270 \, x^8 \, + \, \,  \cdots 
\end{eqnarray}
It is {a} solution of an order-six linear differential operator
$\, L_6 \, = \, \, L_1 \oplus \, L_5$
 which is the direct sum
 of an order-one operator $\, L_1$
 with a simple rational solution $\, {\cal R}(x)$
\begin{eqnarray}
\label{L5rat}
  \hspace{-0.98in}&& \quad \quad  \quad \quad  \quad \quad 
         {\cal R}(x) \, \, = \, \,  \,
         \, {\frac { 1\,  +26\, x\, +198\, x^2\, + 476\, x^3\, + \, 247\, x^4
           }{{ 30 \, \, x}^{5}}}.
\end{eqnarray}

The linear differential operator $\, L_5$ is MUM and 
non-trivially homomorphic to its adjoint. The intertwiners have order 4. 
The \emph{symmetric square} of  $\, L_5$ has a
rational solution of the form $\, P(x)/x^N$, where
$\, P(x)$ is a polynomial. Its differential Galois group is thus 
a special subgroup of 
the orthogonal group $\, SO(5, \, {\mathbb C})$:
\begin{eqnarray}
\label{L5op}
 \hspace{-0.98in} 
\qquad \qquad L_5=\, \left( 4\,x -1 \right)\,  \left( 16\,x-1 \right)\,  \left( 36\,x-1 \right) \, {x}^{5} \, D_x^5 \\
+ 2\, \left( 29376\,{x}^{3}-12232\,{x}^{2}+1106\,x-25 \right) {x}^{4} \, D_x^4\nonumber \\
+ 2\, \left( 241920\,{x}^{3}-124430\,{x}^{2}+14824\,x-449 \right) {x}^{3}
\, D_x^3 \nonumber \\
+ 8\, \left( 190944\,{x}^{3}-122397\,{x}^{2}+20161\,x-893 \right) {x}^{2}
\, D_x^2 \nonumber \\
+ 4\, \left( 425088\,{x}^{3}-343620\,{x}^{2}+82202\,x-6185 \right) x
\, D_x \nonumber \\
+ \, 456192\,{x}^{3}-475344\,{x}^{2}+166912\,x-29160 \nonumber. 
\end{eqnarray}
Once again, the operator satisfied by $\mu_6(x^2)$ 
is  gauge equivalent (homomorphic) with the operator~$A_6$ 
satisfied by the density $p_6(x)$ in~\cite{BSWZ12}.
Therefore,  $p_6(x)$ and our density $\mu_6(x^2)$, both with
support $[0,6]$, are actually \emph{differentially related in a non-trivial
way}. 

\vskip .1cm

\subsection{$Av(12345678)$: order-six operator homomorphic to its adjoint}
\label{Av12345678}

The positive integer coefficient series $F_7(x)$ for $Av(12345678)$
reads\footnote{See \url{https://oeis.org/A072131}; note that the sequence
\url{https://oeis.org/A230051} shares the same first 8 terms,
but it counts permutations avoiding \emph{consecutively} the pattern
$(12345678)$.}:
\begin{eqnarray*}
\label{A072131}
  \hspace{-0.98in}&& 
  1  +x  +2\,{x}^{2}  +6\,{x}^{3} +24\,{x}^{4}
  +120\,{x}^{5}  +720\,{x}^{6}  + 5040\,{x}^{7} + 40319 \,{x}^{8}  + 362815 \,{x}^{9} + \cdots 
\end{eqnarray*}
It is a solution of an order-seven linear differential operator
$\, L_7 \, = \, \, L_1 \, \oplus \, L_6$
 which is the direct sum
 of an order-one operator $\, L_1$
 with a simple rational solution $\, {\cal R}(x)$
 \begin{eqnarray}
\label{L6rat}
  \hspace{-0.98in}&& \quad  \quad 
\mathcal{R}(x) \, \, = \, \,  \,
\, {\frac {1 +37\,x  +447\,{x}^{2}  +2079\,{x}^{3} +3348\,{x}^{4} +1160\,{x}^{5}
  }{ 42 \, \, {x}^{6}}}, 
\end{eqnarray}
 and an order-six operator $L_6$,
\begin{align*}
L_6 = \quad  \left( x-1 \right) \,  \left( 9\,x-1 \right)  \,  \left( 25\,x-1 \right)   \,  \left( 49\,x-1 \right)  {x}^{6} D_x^{6} \\
+  \left( 297675\,{x}^{4}-434024\,{x}^{3}+85686\,{x}^{2}-4704\,x+71 \right)  {x}^{5} D_x^{5} \\
+  \left( 
2679075\,{x}^{4}-4858886\,{x}^{3}+1273932\,{x}^{2}-94818\,x+1913 \right)  {x}^{4} D_x^{4}\\
+ \, 2\,  \left( 4862025
\,{x}^{4}-10985516\,{x}^{3}+3929062\,{x}^{2}-424260\,x+12385 \right)  {x}^{3} D_x^{3}\\
+ \left( 13693050\,{x}^
{4}-38766996\,{x}^{3}+19376513\,{x}^{2}-3329230\,x+160367 \right)  {x}^{2}D_x^{2}
\\+ \left( 5953500\,{x}^{4}-
21425184\,{x}^{3}+15390947\,{x}^{2}-4684008\,x+483601 \right)  x D_x \\+ \, 396900\,{x}^{4}-1887480\,{x}^{3}+2073337\,
{x}^{2}-1002001\,x+518400.
 \end{align*}
The linear differential operator $\, L_6$ is
non-trivially homomorphic to its adjoint. The intertwiners have order~4. 
Its differential Galois group is thus
a special subgroup of
the symplectic group $\, Sp(6, \, {\mathbb C})$.

As before, the operator satisfied by $\mu_7(x^2)$ is gauge equivalent
(homomorphic) with the operator~$A_7$ satisfied by the density $p_7(x)$
in~\cite{BSWZ12}. Therefore, $p_7(x)$ and our density $\mu_7(x^2)$, both with
support $[0,7]$, are actually \emph{differentially related in a non-trivial
way}.

 \vskip .1cm 

\subsection{$Av(1234{\ldots}(k+1))$: order-$(k-1)$ operator homomorphic to its adjoint}
\label{Av123k}

Recall that, by aforementioned results of Rains~\cite{R98} and of
Regev~\cite{R81}, the sequence $ \, Av(1234{\ldots}{k+1})$ is a Stieltjes
moment sequence with support $[0,k^2]$. Further, the generating function
$F_k(x)$ of $Av(123\ldots (k+1))$ is D-finite, as proved by Gessel~\cite{G90},
see also~\cite{GWW98,Novak,Xin10}.
More particularly, it belongs to the interesting subset of D-finite functions
which can be written as diagonals of a rational function, as proved by
Bousquet-M\'elou~\cite[Prop.~13, p.~597]{BM}. It seems likely that many of the
various remarkable properties enjoyed by the sequence $ \,
Av(1234{\ldots}(k+1))$ and by its generating function are consequences of
these diagonal representations.

\vskip .1cm

From all the examples of the previous subsections, it is legitimate to
expect that for any~$k$, the generating function $F_k(x)$ of $ \, Av(1234{\ldots}(k+1))$ is
a series solution of a linear differential operator $\, L_{k}(x,D_x) \,$ of order
$\, k$. The upper bound $k$ can be proven in the spirit
of~\cite[Prop.~1]{BG00}. It also appears that $\, L_{k} \,$ is moreover the
direct sum $\, L_{k} \, = \, \, L_1^{(k)}\, \oplus \, L_{k-1}^{(k)}$ of an
irreducible operator $L_{k-1}^{(k)}(x,D_x)$ of order $k-1$ and of an order-one
operator $L_1^{(k)}$ with a rational solution $\, {\cal R}_{k-2}(x)$ with a
unique pole at $x=0$,
 \begin{eqnarray}
\label{Lkrat-gen}
  \hspace{-0.98in}&& \quad  \quad  \quad  \quad 
{\cal R}_{k-2}(x) \, \, = \, \,  \,
\, {\frac { P_{k-2}(x)}{ k (k-1) \, \,  {x}^{k\, -1}}}, 
\end{eqnarray}
  where the polynomial $\,  P_{k-2}(x)$ has degree $\, k \, -2$ and is
 of the form  
\[P_{k-2}(x) \, = \, \, 1 \, + \, \cdots\, + \,c_{k-2} \,  x^{k-2}. \]
Their coefficients are given in Table~\ref{tab:coeffsP} for ~$k\leq 10$.
Note that they seem to have a very nice combinatorial interpretation: 
for fixed $r$, the sequence of coefficients $\left([x^r] \, P_\ell(x)\right)_{\ell \geq r+1}$ coincides with the sequence
$\left( Q_r(\ell)\right)_{\ell \geq r+1}$
for some polynomial $Q_r(x)$ of degree $2r$ in~$x$.
For instance, when $r=1$, this sequence is $(\ell^2+1)_{\ell \geq 2} = (5, 10, 17, \ldots)$, and when $r=2$, this is $\left(\frac12\,\ell^{4}-\ell^{3}+\frac12\,\ell^{2}-\ell+3\right)_{\ell \geq
3} = (18, 71, 198, \ldots)$. 
Even more remarkably, the evaluation $Q_r(-d)$ of $Q_r(x)$ at \emph{negative integers} $x=-d$ appears to 
count permutations of
length $r + d$ that \emph{do contain} an increasing subsequence of length~$d$. In other words, the polynomials $Q_r(-d)$ are exactly the polynomials $B_d(d + r)$ of~\cite{ESZ}. This surprising coincidence deserves a combinatorial explanation\footnote{Another intriguing property is that $Q_r(r)=r+1$ and $Q_r(r-1)=r$.
}.

\begin{table}[t]
\centering
\begin{tabular}{c|c|c|c|c|c|c|c|c}
$\ell $  $\backslash$  $r$ & 0 & 1 & 2 & 3 & 4 & 5 & 6 & 7\\
\hline
2& 1 & 5 &  &  &  &  &  &\\
3& 1 & 10 & 18 &  & &  &  & \\
4& 1 & 17 & 71 & 63 & &  &  &\\
5& 1 & 26 & 198 & 476 &  247 &  &  &\\
6& 1 & 37 & 447 & 2079 &  3348 & 1160  &  &\\
7& 1 & 50 & 878 & 6668 &  21726 & 25740  & 6588  & \\
8& 1 & 65 & 1563 & 17539 &  95339 & 235755  & 218844 & 44352\\
\end{tabular}
\caption{Coefficients of $x^d$ in the numerators $P_\ell(x)$ of 
$\mathcal{R}_\ell(x)$,
for $\ell \leq 8$.}
\label{tab:coeffsP}
\end{table}

The order-$(k-1)$ linear differential operator $\, L_{k-1}^{(k)}$ appears to
be Fuchsian\footnote{Its finite singularities are at 0, and at $1, \frac19,
\ldots, \frac{1}{k^2}$ if $k$ is odd and at $\frac14, \ldots, \frac{1}{k^2}$
if $k$ is even.}, MUM and homomorphic to its adjoint. As a consequence of
this, it has a differential Galois group which is a subgroup of the orthogonal
group $\, SO(k-1, \, {\mathbb C})$ when $k$ is even, and of the symplectic
group $\, Sp(k-1, \, {\mathbb C})$ when $k$ is odd.

The linear differential operator $\, L_{k-1}^{(k)}$ has a Laurent series
solution $\, \ell_{k-1}(x)$ whose sum with the Laurent series expansion of the
rational function $\, {\cal R}_{k-2}(x)$ is exactly the generating function 
$F_k(x)$ of
$Av(1234{\ldots}(k+1))$:
 \begin{eqnarray}
\label{Lkrat}
\hspace{-0.98in}&& \quad  \quad   \quad  \quad  
 \ell_{k-1}(x) \, +\,   {\cal R}_{k-2}(x)  \, \, = \, \, \, \,     
       1 \, +x \, +2\,{x}^{2} \, +6\,{x}^{3} \, +24\,{x}^{4} \, + \, \, \cdots .
\end{eqnarray}
Therefore, the explicit calculation of the density function $\mu_k(x)$ for the Stieltjes sequence 
$Av(1234{\ldots}(k+1))$ amounts to performing a piecewise analysis of the sum
\begin{eqnarray}
\label{Lkrat1surx}
\hspace{-0.98in}&& \quad  \quad \quad  \quad \quad  \quad \quad
-{{1} \over {\pi y}} \, \, \ell_{k-1}\Bigl(  {{1} \over { y}}  \Bigr)
 \, \, - {{1} \over { \pi y}} \, \, {\cal R}_{k-2}\Bigl(  {{1} \over { y}}  \Bigr).
 \end{eqnarray}
Because of the special form of ${\cal R}_{k-2}(x)$, the function $\, 1/y \cdot {\cal R}_{k-2}(1/y)$ is a {\em polynomial}. The imaginary part
of the evaluation of this polynomial for real values of $\, y$ is zero.
The density function for the Stieltjes sequence 
$Av(1234{\ldots}(k+1))$ is thus a piecewise D-finite function, 
and a solution of the order-$(k-1)$ linear differential operator
\begin{eqnarray}
\label{Lkrat1surx-bis}
\hspace{-0.98in}&& \quad  \quad \quad  \quad \quad  \quad \quad  \quad
  pullback\Bigl( L_{k-1}^{(k)}, \,  {{1} \over {y}} \Bigr) \, \, y.  
\end{eqnarray}
Of course, this  linear differential operator is again (non-trivially) homomorphic to its adjoint, and it has the same differential Galois group as $\, L_{k-1}^{(k)}$.

\smallskip {\bf Remark.} For $k\geq 5$, we exclude the fact that
$L_{k-1}^{(k)}$ could be homomorphic to the $(k-2)$-th power of a second-order
linear differential operator. However, we cannot exclude the fact that the
solutions of $L_{k-1}^{(k)}$, and in particular the generating function
$F_k(x)$ of $Av(12\ldots (k+1))$, minus the rational part $\, {\cal
R}_{k-2}(x)$ in~\eqref{Lkrat-gen}, could be written as an algebraic pullback
of a $_{k-1}F_{k-2}$ hypergeometric function. If we were forced to make a bet,
our guess would be that the solutions might belong to the world of
(specializations of) multivariate hypergeometric functions.

\subsection{The double Borel transform of the  $Av(123\ldots (k+1))$ sequences}
\label{Bergeron}

Another explicit study of the generating functions for the $Av(123\ldots
(k+1))$ sequences was undertaken by Bergeron and
Gascon~\cite{BG00}\footnote{One hundred terms of the generating functions of
$Av(123\ldots (k+1))$ are given in~\cite{ESZ} for $k\leq 60$.}.

In 1990 Gessel \cite{G90} had show{n} that the double Borel transform\footnote{Since the exponential generating function of a power series $F(x) = \sum
a_nx^n$ is also known as the {\em Borel transform} of $F$, we will similarly
call $\sum \frac{a_n}{n!^2} x^n$ the {\em double Borel transform} of $F$.
} of the generating function of $Av(123\ldots (k+1))$ 
could be expressed as the determinant of a $k \times k$ matrix  whose coefficients are Bessel functions.
More precisely, recalling that the generating function of $Av(123\ldots (k+1))$  is $\sum_n
f_{n k} \, x^n,$ with $f_{n k}$ as in~\eqref{eq:Rains} and~\eqref{eq:Rains2},
then the double Borel transform
\[Y_k(x)=\sum_{n\geq 0} \frac{f_{n k}}{n!^2} \,  x^n 
\] 
is equal to the determinant of the Toeplitz matrix 
$\displaystyle{[I_{i-j}(x)]_{1\le i,j\le k}},$
where
\[I_k(x)=\sum_{n=0}^\infty \frac{x^{n+k/2}}{n!(n+k)!}.\]
Clearly, for fixed $k$, the power series $I_k(x)$ is D-finite, so an immediate consequence of Gessel's result is that  $Y_k(x)$ is also D-finite.
Bergeron and Gascon~\cite{BG00} showed that, moreover,
$Y_k(x)$ can  be expressed as homogeneous polynomials in the D-finite functions $I_0(x)$ and $I_1(x)$\footnote{In Maple jargon, these are 
$I_0=\text{BesselI}(0, 2\,  \sqrt{x})$ and $I_1=\text{BesselI}(1, 2\,  \sqrt{x})$.}
\[I_0(x) = 1+x+{\frac{1}{4}}{x}^{2}+{\frac{1}{36}}{x}^{3}+{\frac{1}{576}}{x}^{4}+\cdots,
\quad 
I_1(x) = \sqrt{x} \, \left( 1+{\frac{1}{2}}x+{\frac{1}{12}}{x}^{2}+{\frac{1}{144}}{x}^{3}
+ \cdots \right) \]
satisfying the order-two linear differential equations
\begin{equation}\label{eq:I0I1}
	x  I_0''(x) +  I_0'(x) - I_0(x)  =0 \quad \text{and} \quad
4\,{x}^{2} I_1''(x) 
+4\,xI_1'(x)
-\left( 4\,x+1 \right) I_1(x) 
= 0. 
\end{equation}
For instance, 
\[Y_2 = I_0(x)^{2}-I_1(x)^{2}
\quad \text{and} \quad
Y_3 = 
{\frac {2\,\sqrt {x}I_0(x)^{2}I_1(x)-I_0(x)\,I_1(x)^{2}-2\,\sqrt {x}I_1(x)^{3}}{x}}.
\]
{The underlying differential
operators $\mathcal{L}_2(x), \, \mathcal{L}_3(x), \,\ldots, \,\mathcal{L}_7(x)$
for $Y_2(x), \, Y_3(x),\, \ldots,\, Y_7(x)$
were derived Bergeron and Gascon in~\cite{BG00}.}
 For instance,
\[
\mathcal{L}_2 = {x}^{2}D_x^{3}+4\,xD_x^{2}-2\, \left( 2\,x-1 \right) {
\it D_x}-2
\]
and
\[
\mathcal{L}_3 =
{x}^{3}D_x^{4}+10\,{x}^{2}D_x^{3}-x \left( 10\,x-23
 \right) D_x^{2}- \left( 32\,x-9 \right) D_x+9\,x-9.
\]

The solutions of the differential equations $\mathcal{L}_k(y(x)) =0$ are not
given by regular hypergeometric functions or solutions of Fuchsian linear
differential operators, but rather by solutions of linear differential
operators with an irregular singularity at $\, \infty.$ The differential
operators $\mathcal{L}_k$ are again homomorphic to their adjoint. Moreover,
since the operators~\eqref{eq:I0I1} are homomorphic and the $Y_k$ are
homogeneous polynomials in $I_0$ and $I_1$, one deduces immediately that the
$\mathcal{L}_k$'s are homomorphic to the $k$-th symmetric power of the order
two operator $\mathcal{L}_1 = xD_x^2+D_x-1$.

It is known~\cite{BMW97} that symmetric powers of a fixed second-order
operator are related by a recurrence of order two with simple order-1
operators as coefficients: if $\mathcal{M} = D_x^2 + a(x) D_x + b(x)$, and
then its $n$-th symmetric power is given by the $n$-th term of the sequence of
operators
\begin{eqnarray}
\hspace{-0.58in}\mathcal{M}_0 =1, \quad \mathcal{L}_1 = D_x, \quad \text{and for all } \;  \; k\geq 1, \nonumber \\  \hspace{-0.18in} \mathcal{M}_{k+1} = (D_x + k \cdot a(x)) \cdot \mathcal{M}_k + k (n - k + 1) \cdot b(x) \cdot \mathcal{M}_{k-1}.
\end{eqnarray}
Our operators $\mathcal{L}_k$ are not symmetric powers, but are homomorphic to such
symmetric powers. Consequently, we do not expect a recurrence as simple as the
one on the $\mathcal{M}_k$'s. However, we remark that the $\mathcal{L}_k$'s
still possess some remarkable features. For instance,
\begin{eqnarray}
\hspace{-0.58in}\mathcal{L}_{k+1} = \left(xD_x + \binom{k}{2}+1 \right) \cdot \mathcal{L}_k \nonumber \\  \hspace{-0.18in} + \left(\binom{k+1}{4} - \binom{k+1}{2}x\right) \cdot \mathcal{L}_{k-1} \, + \, \left( \text{terms of ord} \; \leq k-2 \right).
\end{eqnarray}
By analogy with the emergence of these symmetric powers for $\mathcal{L}_k$,
one could imagine that the $L_k$'s associated with $Av(12\ldots (k+2))$ might
also be homomorphic to the symmetric $(k-1)$-th power of an order 2
differential operator, which would hopefully correspond to a classical modular
form.

\section{A difficult case: $Av(1324)$}\label{1324}

In this section, we investigate the density function for the most difficult
Wilf class of length-4 pattern-avoiding permutations, namely~$Av(1324),$ whose
generating function is unknown. The study of this (conjectured) density
function is entirely numerical, based on the extensive exact enumerations of
Conway, Guttmann and Zinn-Justin~\cite{CGZ18}.

The series with positive integer coefficients for $Av(1324)$ is \[1+x+2\,
x^2+6\, x^3+23\, x^4+103\, x^5+513\, x^6+2762\, x^7+15793 \, x^8+\ldots,\] and
is known \cite{CGZ18} up to order $x^{50}.$ The generating function is not
known, and there are compelling arguments \cite{CGZ18} that \emph{it is not
D-finite}. While it is known that the coefficients grow exponentially as
$\lambda^n,$ the value of $\lambda$ is not known. The best estimate is in
\cite{CGZ18} and is $\lambda=11.60 \pm 0.01.$ The best rigorous bounds
\cite{BBEP17} are $10.27< \lambda < 13.5$, {while the paper~\cite{CJS12}
gives the improved upper bound $e^{\pi \, \sqrt{2/3}} \approx 13.001954$, 
but only
under a certain conjecture about
the number of inversions in a 1324-avoiding permutation.}

The Hankel determinants are all positive and are monotonically increasing, which provides strong evidence for the conjecture
that the series is a Stieltjes moment sequence. 

 Recall that if the sequence $\seq$ is positive, then its
log-convexity implies that the coefficient ratios $\frac {a_n}{a_{n-1}}$ are lower bounds
on the growth rate $\mu$ of the sequence. In this way we obtain $\lambda > 9.03.$

In the case that $\seq
$ is a
Stieltjes moment sequence, 
stronger lower bounds for $\mu$ can be calculated using a method first given by Haagerup,
Haagerup and Ramirez-Solano~\cite{HHR15}, that we now summarize.

Using the coefficients $a_0, a_1, \ldots \, a_n,$ one calculates the terms $\,
\alpha_0, \alpha_1, \ldots \alpha_n$ in the continued fraction representation given in Theorem~\ref{thm:Stieltjes}.
 It is easy to see that the coefficients of $A(x)$ are nondecreasing in
each $\alpha_j.$ Hence $A(x)$ is (coefficient-wise) bounded below by the
generating function $A_n(x)$, defined by setting $ \, \alpha_n,
\alpha_{n+1},\, \alpha_{n+2} \, \ldots,$ to 0. Therefore, the growth rate
$\mu_n$ of $A_n(x)$ is not greater than the growth rate $\mu$ of $A(x).$ The
growth rates $ \, \mu_1,\, \mu_2, \, \ldots $ clearly form a non-decreasing
sequence, and, since the coefficients of $A_n(x)$ are log-convex, $\mu_n \ge
a_n/a_{n-1}$. It follows that this sequence $ \, \mu_1, \, \mu_2, \, \ldots $
of lower bounds converges to the exponential growth rate $\mu$ of \seq.

If we assume further that the sequences $ \, \alpha_0, \, \alpha_2, \,
\alpha_4 \, \ldots $ and $\alpha_1, \alpha_3, \alpha_5 \ldots $ are
non-decreasing, as we find empirically in many of the cases we consider, we
can get stronger lower bounds for the growth rate by setting $ \,
\alpha_{n+1}, \, \alpha_{n+3} \, \ldots$ to $ \, \alpha_{n-1}$ and $
\,\alpha_{n+2}, \, \alpha_{n+4} \, \ldots$ to $\alpha_{n}$. For this sequence
the exponential growth rate of the corresponding sequence $\seq
$ is
$(\sqrt{\alpha_n} \, +\sqrt{\alpha_{n-1}})^2.$ By the method with which we
constructed this bound, it is clear that $(\sqrt{\alpha_n}
\,+\sqrt{\alpha_{n-1}})^2 \ge \mu_n.$ Hence, the lower bounds
$(\sqrt{\alpha_n} \,+\sqrt{\alpha_{n-1}})^2$ converge to the growth rate
$\mu.$

\begin{figure}[htb]
 \begin{center}\includegraphics[scale =1.]{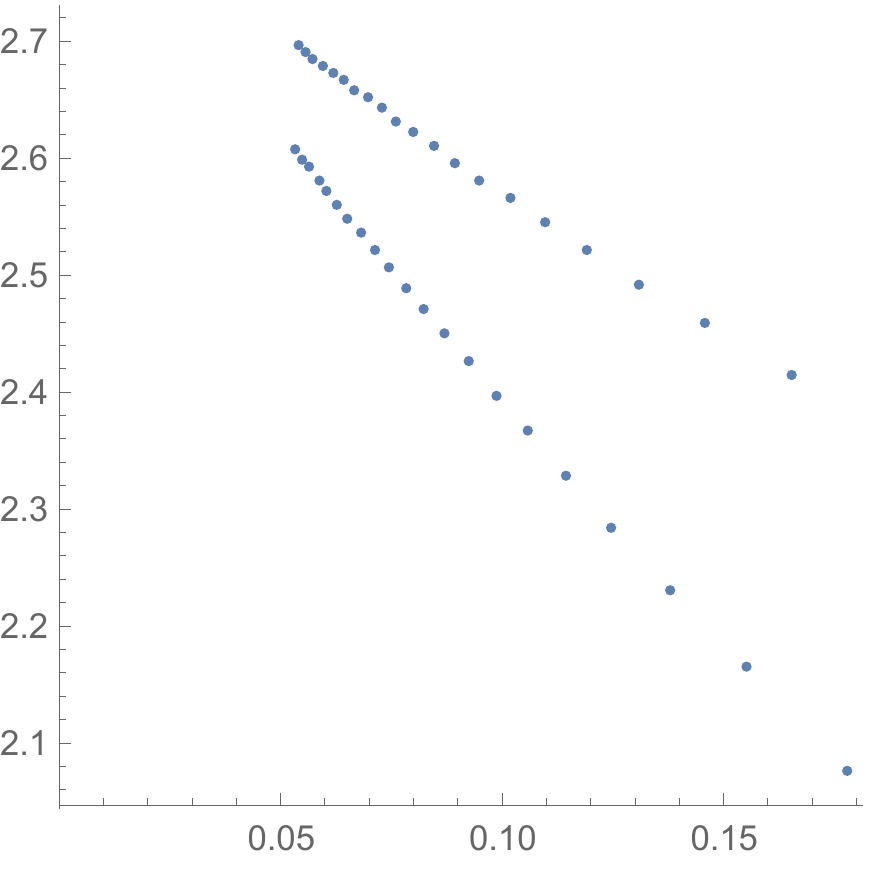}\end{center}
 \caption{Plot of $\alpha_n$ (vertical axis) against $n^{-2/3}$ (horizontal axis) for the sequence $Av(1324)$ using coefficients $10 \ldots 50$.}
 \label{fig:1324-alpha}
\end{figure}

  In particular, if $\alpha_n \le \alpha_{n+2}$ for each $n$ and the
limit $\lim_{n \to \infty} \alpha_n$ exists, then it is equal to $\mu/4.$
Extensive use of this result in studying the cogrowth sequences of various
groups has been made in~\cite{HHR15} and~\cite{EG19}.
{A plot of the known values of $\alpha_{n}$ for $Av(1324)$ is shown in Fig.~\ref{fig:1324-alpha}.}

 It appears that the coefficients $\alpha_{2k}$ and the coefficients
$\alpha_{2k+1}$ each form an increasing sequence which is roughly linear when
plotted against $n^{-2/3}.$ 
If this trend continues, the terms $\alpha_n$ will
certainly all be positive. If it could be proved that this is a Stieltjes
moment sequence, these $\alpha$ values would imply a lower bound of $\lambda
>10.302,$ which is an improvement on the best lower bound cited above. If we
assume further that the sequences $(\alpha_{2k})_{k \ge 0}$ and
$(\alpha_{2k+1})_{k\geq 0}$ are both increasing, we get an improved lower
bound $\lambda > 10.607.$

Extrapolating the two subsequences $(\alpha_{2k})_{k \ge 0}$ and
$(\alpha_{2k+1})_{k\geq 0}$, it seems plausible that these both converge to
the same constant $c \approx 2.9,$ which is consistent with the growth rate
$\lambda \approx 11.60$ predicted by Conway, Guttmann and Zinn-Justin
\cite{CGZ18}.

\begin{figure}[t]
\centering
\includegraphics[scale =0.44] {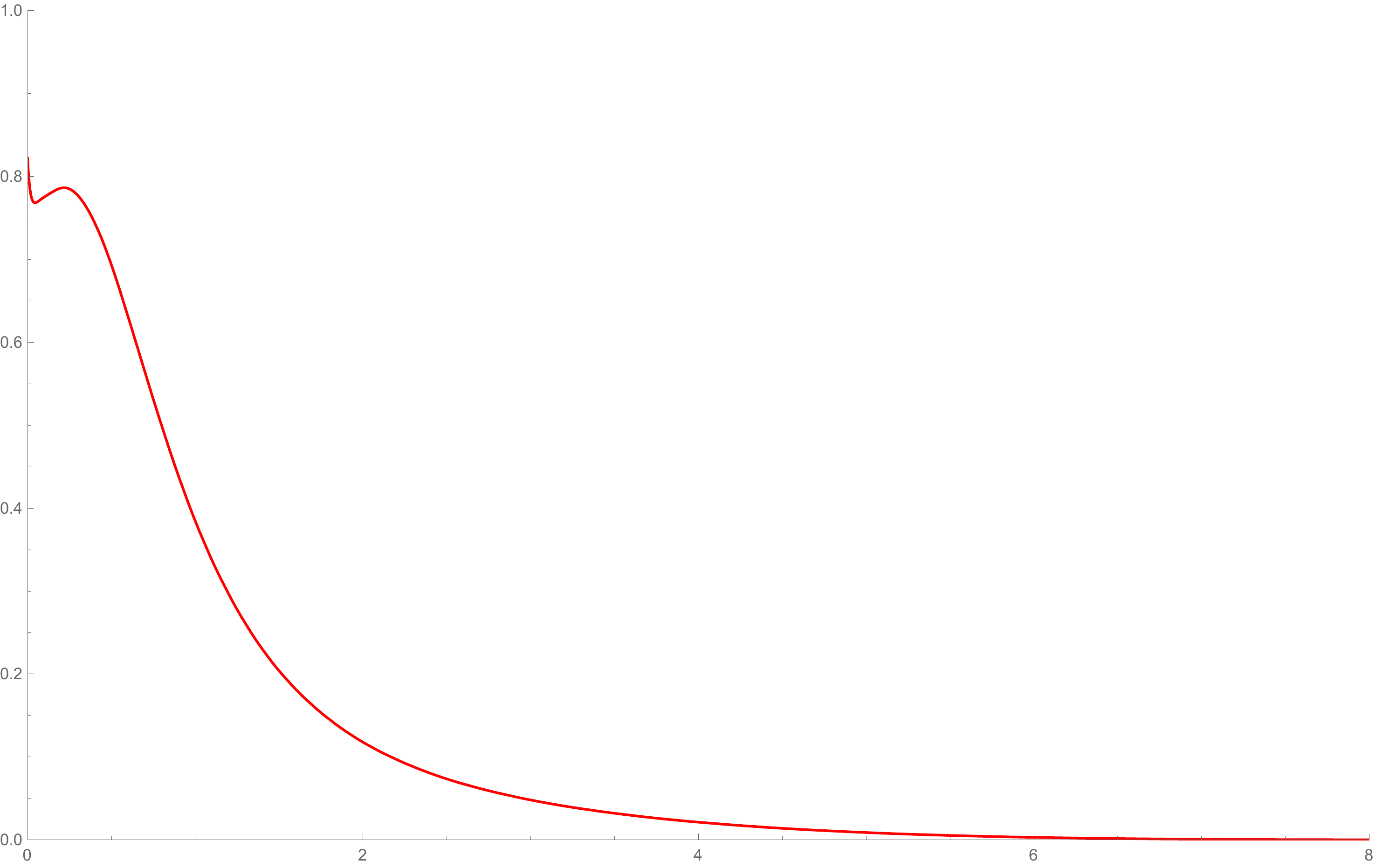}
 \caption{The density function for $Av(1324)$ constructed numerically.}
 \label{fig:1324}
\end{figure}

From the numerical data, a histogram of the density function can be
constructed, and this is shown in Fig.~\ref{fig:1324}. We have no
explanation for the little wiggle near the origin.

\subsection{Numerics}\label{sec:num}

We have used two different types of numerical construction of the density functions.

\paragraph{}
The {\em histogram} for {$Av(12\ldots k)$} (D-finite, discussed in \S\ref{Av1234}) is constructed as follows:
\begin{itemize}

\item take a random {$(k-1)\times(k-1)$} matrix with entries which independently have standard
complex normal distributions

\item orthogonalize it to create a random unitary matrix

\item compute the norm of the trace squared (by the Rains result~\cite{R98}
the $n$-th moment of this distribution is exactly the number of
$12\ldots k$-avoiding permutations of size~$n$)

\item repeat the above steps $500\,000$ times and draw a histogram of the
results, where each bar has width $0.02$

\end{itemize}
Note that this procedure does not use the exact coefficients.

\paragraph{}

The {\em polynomial approximations} for the densities are constructed by
computing the moments $a_0, a_1,\ldots,a_n$ (that is, just the initial terms
of the sequence), then approximating the density function by the unique
polynomial $P(t)$ of degree $n+2$ over the known range $[0,\mu)$ with the
following properties:

\begin{itemize}
\item the initial moments of the distribution with density $P(t)$ are $a_0, a_1,\ldots, a_n$
\item $P(\mu)=P'(\mu)=0$
\end{itemize}
For $Av(1342)$, we have $\mu=8$. For $Av(1324)$ we chose $\mu=12$, though for any value of $\mu$ between 11 and 13 the figure looks almost identical.
 
\subsection{Final remarks}\label{sec:rem}
                        
\smallskip 
{
\quad {\bf \ Remark 1.} For each permutation $\pi$ of length at most $4$, we have proved or given strong evidence that the sequence $\big(Av_{n}(\pi)\big)_{n \geq 0}$,
{whose general term is the number of permutations of
$\{1,\ldots, n \}$ that avoid the pattern $\pi$},
 is a Stieltjes moment sequence. 
Indeed, there is only one Wilf class of length 3, namely $(123)$, considered in \S\ref{sec:Catalan} and there are 3 Wilf classes of length 4, namely $(1342), (1234)$ and $(1324)$. We proved that the first three are Stieltjes moment sequences (in \S\ref{sec:Catalan}, \S\ref{1342} and \S\ref{Av1234}) and presented numerical evidence (in \S\ref{sec:num}) for the last one. Moreover, the same equally holds for any $\pi$ of the form $(12\ldots k)$, as showed in \S\ref{Increasing}. 
This naturally suggests the further question:
\vspace{-0.1cm}
\begin{quote}
	{\bf Open question:} Is it the case that, for every permutation $\pi$, the sequence $\big(Av_{n}(\pi)\big)_{n \geq 0}$ is a Stieltjes moment sequence? 
\end{quote}	
\vspace{-0.1cm}
We tested all known terms of these sequences and found they are consistent with being Stieltjes moment sequences. However,
to our knowledge, substantial computations of the initial terms of these sequences have only been done in the cases discussed in our paper.
}

\medskip 
{
{\bf Remark 2.} One of the referees suggested that the classes of
triples and pairs of length-4 pattern-avoiding permutations enumerated by Albert et al.~\cite{AHPSV18},
and conjectured to be non-D-finite, should be similarly studied. These sequences are
$Av(4123, 4231, 4312)$, $Av(4123, 4231)$, $Av(4123, 4312)$ and $Av(4231,
4321)$, and their first terms are given in the OEIS as 
\href{https://oeis.org/A257562}{A257562},   
\href{https://oeis.org/A165542}{A165542}, \href{https://oeis.org/A165545}{A165545}, \href{https://oeis.org/A053617}{A053617},
respectively. Unfortunately, they all have Hankel determinants that become
negative after a certain order, and so cannot be described as Stieltjes moment
sequences.
}

\section{Conclusion}
\label{Conclusion}

We have shown how considering coefficients of combinatorial sequences as
moments of a density function can be useful. Quite often in the literature the
analysis of Stieltjes moment sequences with explicit densities is performed
only for algebraic series which are slightly over-simplified degenerate cases.

We have studied here in some detail examples that are not algebraic series.
For instance we have show{n} that $\, Av(1234)$ is a Stieltjes sequence whose
generating function corresponds, up to a simple rational function, to an
order-one linear differential operator acting on a classical modular form
represented as a pulled-back $\, _2F_1$ hypergeometric function, and $\,
Av(12345)$ is a Stieltjes sequence whose generating function, up to a simple
rational function, corresponds to an order-two linear differential operator
acting on the square of a classical modular form represented as a pulled-back
$\, _2F_1$ hypergeometric function. The corresponding densities are of the
same type, and do not involve the rational function. This scheme generalizes
to all the $\, Av(12345\dots k)$ which are all Stieltjes moment sequences.

The linear differential operators annihilating such series are direct sums of
an order-one operator with a simple rational solution where the denominator is
just a power of $\, x$, and a linear differential operator homomorphic to its
adjoint (thus corresponding to selected differential Galois groups). This last
operator has a Laurent series solution, the generating function of the
Stieltjes sequence corresponding to getting rid of the finite number of poles.
The density is the solution of the $\, x \, \rightarrow \, 1/x$ pullback of
this last operator, but finding the actual density requires delicate piecewise
analysis which is difficult to perform when one does not have exact
expressions for the generating function of the Stieltjes moment sequence.

We have show{n} that the density function for the Stieltjes moment sequence
$Av(12\ldots k)$ is closely, but non-trivially, related to the density
attached to the distance traveled by a walk in the plane with $k-1$ unit steps
in random directions.

Finally, we have considered the challenging (and still unsolved!) case of
$Av(1324)$, and provided compelling numerical evidence that the corresponding
counting sequence is a Stieltjes sequence. Assuming this, we obtained lower
bounds to the growth constant that are better than existing rigorous bounds.

\bigskip

\noindent {\bf Acknowledgments.} {We thank the referees for their very
careful reading and for providing many useful and constructive suggestions.}
Jean-Marie Maillard thanks The School of Mathematics and Statistics, The
University of Melbourne, Australia, where part of this work has been
performed. Alin Bostan and Jean-Marie Maillard address their warm thanks to
the Caf\'e~du~Nord (Paris, France) where they had the chance to work daily
during the massive French pension reform strike (December 2019). Alin Bostan
has been supported by
\textcolor{magenta}{\href{https://specfun.inria.fr/chyzak/DeRerumNatura/}{DeRerumNatura}}
ANR-19-CE40-0018. Andrew Elvey Price has been supported by European Research
Council (ERC) under the European Union's Horizon~2020 research and innovation
programme under the Grant Agreement No.~759702. Anthony~John~Guttmann wishes
to thank the ARC Centre of Excellence for Mathematical and Statistical
Frontiers (ACEMS) for support.

\newpage

\appendix

\section{Density of $Av(1234)$: Piecewise analysis}
\subsubsection*{Discussion of the evaluation of $\, _2F_1$ in Maple}
\label{app:maple}

As far as the evaluation of a Gauss hypergeometric
series like $\,_2F_1(-1/4,3/4;\, 1;  \, H)$ goes, 
the hypergeometric series converges for $\, |H| < 1$, and 
 $\,_2F_1(a,b;\, c, \, H)$
is then defined for $\, |H| \ge \, 1$ {\em by analytic 
  continuation}. The point $z = 1$ is a branch point, and
  the interval $ [1, \, \infty]$  is the branch cut, with the evaluation on this cut taking the limit from below the cut.
We need to look at~\cite{Johansson1,Johansson2,Becken}. The  hypergeometric function
$\,_2F_1(a,b;\, c; \, H)$  can be rewritten\footnote{See section 7.1
  page 15 in ~\cite{Johansson1}.} using the Pfaff
connection formulas: 
\begin{eqnarray}
\label{Johans}
\hspace{-0.98in}&& \quad  \quad  \quad 
_2F_1 (a, \, b; \, c; \, H)
\, \, \, \, = \, \, \, \,
(1\, -H)^{-a} \,  \,
_2F_1\Bigl(a, \, c \, -b; \, c; \,\, - \, {{H} \over {1\, -H}} \Bigr)
\nonumber \\
\hspace{-0.98in}&& \quad \quad  \quad  \quad  \quad  \quad  \quad  \quad  \quad 
\, \, = \, \, \,
(1\, -H)^{-b} \,  \,
_2F_1\Bigl(c\, -a, \, b; \, c; \,\, - \, {{H} \over {1\, -H}}  \Bigr). 
\end{eqnarray}
These formulas enable one to evaluate $\, _2F_1(-1/4, 3/4;\, 1;  \, H)$ for
{\em negative real values of} $\, H$, changing them  {to the} evaluation of  
$\, _2F_1(3/4, 5/4;\, 1;  \, \tilde{H})$ or
$\, _2F_1(-1/4, 1/4;\, 1;  \, \tilde{H})$
with $\, \tilde{H} \, \in \, [0, \, 1]$.
In our case this gives for $\, H \, < \, 0$:
\begin{eqnarray}
\label{Johans0}
\hspace{-0.98in}&& \quad \quad   \,
_2F_1\Bigl(-{{1} \over {4}}, {{3} \over {4}}; \, 1;  \, H\Bigr)
\, \, \, = \, \, \, \,  (1 \, -H)^{1/4} \,  \,
_2F_1\Bigl(-{{1} \over {4}}, \, {{1} \over {4}};\, 1; 
\,  \, -\, {{H} \over {1 \, -H}} \Bigr),
\end{eqnarray}
which corresponds to  real positive values.  

To evaluate   $\, _2F_1(-1/4, 3/4;\, 1;  \, H)$ for
positive real values of $\, H$ larger than~$\ 1$ with ($H  \in  (1, \, \infty] $),
we must use the connection formula  15.10.21 in \url{https://dlmf.nist.gov/15.10\#E21}
(see~\cite{Johansson1,Johansson2,Becken}): 
\begin{eqnarray}
\label{Johans2}
\hspace{-0.98in}&& \quad \quad 
 _2F_1\Bigl(a, \, b; \, c; \, H \Bigr)
\, \, = \, \, \,
   {{ \Gamma(a) \,\Gamma(b)  } \over {\Gamma(c-a) \,\Gamma(c-b) }}  \, \, w_3
   \, \, \, +   {{ \Gamma(c) \,\Gamma(a+b-c)  } \over {\Gamma(a) \,\Gamma(b) }}  \, \, w_4, 
\end{eqnarray}
where $\, w_3$ is alternatively one of the two formulas: 
\begin{eqnarray}
\label{Johans33}
\hspace{-0.98in}&&  \quad  \quad  \quad  \, 
w_3 \, \, = \, \, \,
_2F_1\Bigl(a, \, b; \, a \, +b \, -c +1; \,  1 \, -\,  H \Bigr), 
\nonumber \\
\hspace{-0.98in}&&  \quad  \quad  \quad  \quad  \quad  \, \, \, = \, \, \,
H^{-a} \,  \,
_2F_1\Bigl(a, \, a \, -c \, +1; \, a \, +b \, -c +1;
\,  1 \, -\,  {{1} \over {H}} \Bigr),
\end{eqnarray}
and  $\, w_4$ is alternatively one of the two formulas:
\begin{eqnarray}
\label{Johans34}
\hspace{-0.98in}&& \quad   \, 
w_4 \, \, = \, \, \,\, \,(1\, -H)^{c-a-b}  \,  \,
_2F_1\Bigl(c \, -a, \, c \, -b; \, c \, -a \, -b +1; \,  1 \, -\,  H \Bigr)
 \\
\hspace{-0.98in}&& \quad \quad   \, 
\, \, = \, \, \,\, 
H^{a-c} \,  \, (1\, -H)^{c-a-b}\,  \,
_2F_1\Bigl(c \,-a, \, 1 \, -a; \, c \, -a \, -b \, +1;
\,  1 \, -\,  {{1} \over {H}} \Bigr).
\nonumber
\end{eqnarray}
In our case this gives for $\, H \, > \, 1$:
\begin{eqnarray}
\label{Johans345}
\hspace{-0.98in}&& \quad \quad   \,
_2F_1\Bigl(-{{1} \over {4}}, {{3} \over {4}}; \, 1;  \, H\Bigr)
\, \, \, = \, \, \, \, 
 {{ 2\, \Gamma(3/4)^2 } \over {\pi^{3/2} }}
 \, \,  H^{1/4} \,  \,
 _2F_1\Bigl(-{{1} \over {4}}, \, -{{1} \over {4}}; \,
       {{1} \over {2}}; \, 1 \, - \, {{1} \over {H}} \Bigr) 
 \nonumber \\
 \hspace{-0.98in}&& \quad \quad \quad  \quad  \quad   \,
 \, \, +  {{ \pi^{1/2} } \over { 2 \, \Gamma(3/4)^2  }} \, \,
 H^{-5/4} \,  \, (1\, -H)^{1/2}  \,  \,
 _2F_1\Bigl({{5} \over {4}}, \, {{5} \over {4}};\, {{3} \over {2}};
 \, 1 \, - \, {{1} \over {H}} \Bigr).
\end{eqnarray}
For $\, H \, > \, 1$, that the first  term on the RHS
of (\ref{Johans345}) corresponds  to {\em positive real values}
while the second  term on the RHS
of (\ref{Johans345}) corresponds to purely imaginary\footnote{Coming from
the $\,  (1\, -H)^{1/2}$ term when $\, H \, > \, 1$.}
values:
\begin{eqnarray}
\label{Johans345Im}
\hspace{-0.98in}&& \quad \quad   \,
\Im \Bigl(\,   _2F_1\Bigl(-{{1} \over {4}}, {{3} \over {4}}; \, 1;  \, H\Bigr) \Bigr)
\nonumber \\
\hspace{-0.98in}&& \quad \quad  \quad \quad   \,
\, \, \, = \, \, \, \, 
   {{ \pi^{1/2} } \over { 2 \, \Gamma(3/4)^2  }} \, \,
 H^{-5/4} \,  \, (H\, -1)^{1/2}  \,  \,
 _2F_1\Bigl({{5} \over {4}}, \, {{5} \over {4}}; \, {{3} \over {2}};
 \, 1 \, - \, {{1} \over {H}} \Bigr).
\end{eqnarray}
See equation 15.10.21 in \url{https://dlmf.nist.gov/15.10\#E21}

Note that both Maple and Mathematica compute these by defining $\,
_2F_1(a,b;c;z)$ to be analytic in~$z$ except on the cut $\, [1,\infty]$,
and, on this cut, they both take {\em the limit from below the cut}.

\vskip .3cm 

\subsubsection*{The piecewise analysis}
\label{piecewise}

For   $\, y \, \in \, \, [0, \, 1]$
the pullback gives $\, H  \, \in \, \, (- \infty, \, 0]$, so $\,  _2F_1\Bigl(-{{1} \over {4}}, \, {{3} \over {4}};
\, 1;\, \, H \Bigr) \, \in \, \,[1, \, \infty]$.

For   $\, y \, \in \, \, [1, \,  6\sqrt{3}-9]$ the  pullback
gives  $\, H  \, \in \, \, [1, \, \infty]$.
For $\, y \, \in \, \, [6\sqrt{3}-9, \, 9]$
the   pullback
 also gives  $\, H  \, \in \, \, [1, \, \infty]$. 
Since {\em  the imaginary part of a polynomial
  (such as $\, {{y+5} \over {6}}$) evaluated at real values is zero},
the  imaginary part of (\ref{inversemu2rho}) is the
imaginary part of the pulled-back hypergeometric term in (\ref{inversemu2rho})
which reads:

\vskip .3cm

$\, \bullet$ For $\, y \, \in \, [0, \, 1]$ (i.e.  $\, H \, \in \, [-\infty, \, 0]$),
and $ _2F_1(-1/4, 3/4; \, 1;  \, H)$
is given by (\ref{Johans0}) and is real and positive), 
the imaginary part  of (\ref{inversemu2rho})
\begin{eqnarray}
\label{inversemu2rhoy010}
\hspace{-0.98in}&&  \quad \quad \quad \quad \quad
\frac {{ (y-1)^{1/4} \,  (y-9)^{3/4}} } {6}
\,   \,  _2F_1\Bigl(-{{1} \over {4}}, {{1} \over {4}}; \, 1;  \, H\Bigr).
\end{eqnarray}
can be rewritten, using  (\ref{Johans0}), as the imaginary part of
\begin{eqnarray}
\label{inversemu2rhoy01-bis}
\hspace{-0.98in}&&  \quad \quad \quad \quad \quad
 \frac{  \Bigl( y^2+18\, y \, -27 \Bigr)^{1/2}}{6}
 \,  \,  _2F_1\Bigl(-{{1} \over {4}}, \, {{1} \over {4}};  \, 1; 
  \, \, \, {{64 \, y^3 } \over {  (y^2+18\, y \, -27)^2  }} \Bigr),  
\end{eqnarray}
which reads, since $\, y^2+18\, y \, -27 \, < 0$:
\begin{eqnarray}
\label{inversemu2rhoy01-ter}
\hspace{-0.98in}&&  \quad \quad   \Im\Bigl(\theta(y)\Bigr)  \, = \, \, 
\frac{  (27 \, -18 \, y -y^2)^{1/2}}{6}
 \,  \,  _2F_1\Bigl(-{{1} \over {4}}, \, {{1} \over {4}};  \, 1; 
  \, \, \, {{64 \, y^3 } \over {  (y^2+18\, y \, -27)^2  }} \Bigr).
\end{eqnarray}

$\, \bullet$ For $\, y \, \in \, [1, \, 9]$ the imaginary part
of (\ref{inversemu2rho}) is thus the imaginary part  of
\begin{eqnarray}
\label{inversemu2rhoy010bis}
\hspace{-0.98in}&&  \quad \quad \quad \quad \quad
 {{ (y-1)^{1/4} \,  (y-9)^{3/4}} \over {6}} 
\,   \,  _2F_1\Bigl(-{{1} \over {4}}, {{3} \over {4}}; \, 1; \, H\Bigr).
\end{eqnarray}
This can be rewritten, using  (\ref{Johans345}), as the imaginary part of
\begin{eqnarray}
\label{Johans345as}
\hspace{-0.98in}&& 
 {{ 2\, \Gamma(3/4)^2 } \over {\pi^{3/2} }} \,  \,  {{ (y-1)^{1/4} \,  (y-9)^{3/4}} \over {6}} 
 \,  \,  H^{1/4} \,  \,
 _2F_1\Bigl(-{{1} \over {4}}, \, -{{1} \over {4}}; \, {{1} \over {2}}; 
 \, 1 \, - \, {{1} \over {H}} \Bigr) 
  \\
 \hspace{-0.98in}&& 
 +  {{ \pi^{1/2} } \over { 2 \, \Gamma(3/4)^2  }}\,  \,  {{ (y-1)^{1/4} \,  (y-9)^{3/4}} \over {6}}  \,  \,
 H^{-5/4} \,  \, (1\, -H)^{1/2}  \,  \,
 _2F_1\Bigl({{5} \over {4}}, \, {{5} \over {4}}; \, {{3} \over {2}}; 
 \, 1 \, - \, {{1} \over {H}} \Bigr).
 \nonumber
\end{eqnarray}
Note that  the first term  in (\ref{Johans345as}) is just
a complex number because of the factor
$\, (y-9)^{3/4}$ (all the other factors are real and positive), while
the second term
in (\ref{Johans345as}) is a complex number because of the factor $\, (y-9)^{3/4}$ 
combined with the {\em  pure imaginary factor} $\, (1 \, -H)^{1/2}$
(all the other  factors are real and positive).

After simplification, this can be written as 
\begin{eqnarray}
\label{inversemu2rhoy19}
\hspace{-0.98in}&&  \quad\quad \quad  \quad
       {{ \Gamma(3/4)^2 } \over { 3 \, \pi^{3/2} }}
 \, \,  (-64 \, y^3)^{1/4} \,  \,
 _2F_1\Bigl(-{{1} \over {4}}, \, -{{1} \over {4}}; \, {{1} \over {2}}; 
 \,  {{ (y^2+18\, y \, -27)^2 } \over {64 \, y^3 }} \Bigr) 
 \nonumber \\
 \hspace{-0.98in}&& \quad \quad \quad 
  +  {{ \pi^{1/2} \, i } \over { 12 \, \Gamma(3/4)^2  }} \, \,
  (y-1) \,  \, (y-9)^3  \,  \,  (-64 \, y^3)^{-5/4} \,  \, ((y^2+18\, y \, -27)^2)^{1/2}
\nonumber \\
 \hspace{-0.98in}&& \quad \quad \quad \quad  \quad \quad \quad  \quad  \quad  \,  \,
 _2F_1\Bigl({{5} \over {4}}, \, {{5} \over {4}}; \, {{3} \over {2}}; 
 \,  {{ (y^2+18\, y \, -27)^2 } \over {64 \, y^3 }} \Bigr).
\end{eqnarray}

Since $\, (-1)^{1/4}$ is $\,\, (1+i)\, 2^{-1/2},$  $\, (-1)^{-5/4}$ is
$\,\, (i-1) \, 2^{-1/2}$,  $\, 64^{1/4} \, = \, 2^{6/4} \, = \,  2^{3/2}$,
and $\, 64^{-5/4} \, = \, 2^{-30/4} \, = \,  2^{-15/2}$, 
one can rewrite (\ref{inversemu2rhoy19}) as:
\begin{eqnarray}
\label{inversemu2rhoy20}
\hspace{-0.98in}&&  \quad\quad \quad  \quad
     {{ 2} \over {3   }}    {{  \Gamma(3/4)^2 } \over { \pi^{3/2} }}
 \, \, (1\, +i) \, \, y^{3/4} \,  \,
 _2F_1\Bigl(-{{1} \over {4}}, \, -{{1} \over {4}}; \, {{1} \over {2}}; 
 \,  {{ (y^2+18\, y \, -27)^2 } \over {64 \, y^3 }} \Bigr) 
 \nonumber \\
 \hspace{-0.98in}&& \quad \quad \quad 
  +  {{ \pi^{1/2} } \over { 3072 \, \Gamma(3/4)^2  }} \, \, (i-1)  \, \, i \cdot \,
  (y-1) \,  \, (9 \, -y)^3  \,   \,  y^{-15/4} \,  \, |y^2+18\, y \, -27|
\nonumber \\
 \hspace{-0.98in}&& \quad \quad \quad \quad  \quad \quad \quad  \quad  \quad  \, \,  
 _2F_1\Bigl({{5} \over {4}}, \, {{5} \over {4}}; \, {{3} \over {2}}; 
 \,  {{ (y^2+18\, y \, -27)^2 } \over {64 \, y^3 }} \Bigr).
\end{eqnarray}

This should give an imaginary part which reads, for $\, y \in [1, \, 9]$, 
\begin{eqnarray}
\label{inversemu2rhoy21}
\hspace{-0.98in}&&  \quad \quad \quad  \quad
\Im\Bigl(\theta(y)\Bigr) \,  \, = \, \,  \, 
     {{ 2} \over {3   }}    {{  \Gamma(3/4)^2 } \over { \pi^{3/2} }}
  \, \, y^{3/4} \,  \,
 _2F_1\Bigl(-{{1} \over {4}}, \, -{{1} \over {4}}; \, {{1} \over {2}}; 
 \,  {{ (y^2+18\, y \, -27)^2 } \over {64 \, y^3 }} \Bigr) 
 \nonumber \\
 \hspace{-0.98in}&& \quad \quad \quad  \quad \quad  \quad 
  +   {{ \pi^{1/2} } \over { 3072 \, \Gamma(3/4)^2  }} \, \,
  (y-1) \,  \, (9 \, -y)^3  \,   \,  y^{-15/4} \,  \, (27  - 18\, y \, - y^2)
\nonumber \\
 \hspace{-0.98in}&& \quad \quad  \quad \quad  \quad \quad  \quad \quad \quad  \quad  \quad  \, \, 
 _2F_1\Bigl({{5} \over {4}}, \, {{5} \over {4}}; \, {{3} \over {2}};
 \,  {{ (y^2+18\, y \, -27)^2 } \over {64 \, y^3 }} \Bigr).
\end{eqnarray}

\section{$Av(1234)$ as a derivative of a classical modular form}
\label{Av1234as-App}

First note that 
\begin{eqnarray}  \label{involution9x}
 \hspace{-0.98in}&& 
(1-6x-3x^2)\, (1-9x)^{3/2}\, {_2F_1}\left ( \frac{3}{4},\frac{3}{4} ;  1  ; \frac{-64x^3}{(1-x)^3(1-9x)} \right )=\\ \nonumber
\hspace{-0.98in} &&
(1-x)^{3/2}\, (1+18x-27x^2)\, {_2F_1}\left ( \frac{3}{4},\frac{3}{4} ;  1  ; \frac{-64x}{(1-x)(1-9x)^3} \right ).
\end{eqnarray}

The fact that ${_2F_1}\left ( \frac{3}{4},\frac{3}{4} ;  1  ; \star \right )$ is a classical modular form gives rise to this simple but infinite-order transformation.

For the non-classical modular form emerging naturally for $Av(1234)$ a more complicated, but analogous, transformation exists. It is:
\begin{eqnarray}
\label{pullbackL21over9xidcalHter}
  \hspace{-0.98in}&&  \quad
  - \, {{1} \over {2}}   \, \,  x \,  \,(1\, -x)^{3/4} \,  \, (1\, -9\, x)^{1/4}  \,  \,
 _2F_1\Bigl(-\, {{1} \over {4}}, \, {{3} \over {4}}; \, 1; 
 \, \, {{ - \, 64 \, x^3 } \over { (1\, -x)^3 \,  \, (1 \, -9 \, x)}}  \Bigr)
 \nonumber  \\
  \hspace{-0.98in}&&     \quad \quad \quad \quad
\, \, = \, \,  \,  \,   {{ x^3} \over {1 \, - \, 9 \, x}}
 \,  \, \Bigl( 8 \,  \, (x-1) \,  \, x \,  \, {{ d {\cal H}(x)} \over { d x}}
 \, \, +(5\, x \,-13) \,  \, {\cal H}(x) \Bigr),             
\end{eqnarray}
where:
\begin{eqnarray}
\label{hypL2recall}
\hspace{-0.99in}&&
 {\cal H}(x) \, \, = \, \,
       - \, {{(1\, -x)^{1/4} \,  \, (1\, -9\, x)^{3/4}} \over { 6 \, x^2}}  \,  \,
 _2F_1\Bigl(-\, {{1} \over {4}}, \, {{3} \over {4}}; \, 1; 
 \, \, {{ - \, 64 \, x } \over { (1\, -x) \,  \, (1 \, -9 \, x)^3}}  \Bigr).
\end{eqnarray}

The relation between these two Hauptmoduls
\[\, A \, = \, \,  {{ - \, 64 \, x } \over { (1\, -x) \,  \, (1 \, -9 \, x)^3}} \qquad  \text{and} \qquad
\, B \, = \, \,  {{ - \, 64 \, x^3 } \over { (1\, -x)^3 \,  \, (1 \, -9 \, x)}} \, \, \]
is the modular equation:
\begin{eqnarray}
\label{Mod13}
\hspace{-0.18in}&&   
4096\,{A}^{3}{B}^{3} \, \, -4608\,{A}^{2}{B}^{2} \, \, (A+B)\,
\, \, -{A}^{4}+900\,{A}^{3}B-28422\,{A}^{2}{B}^{2}+900\,A{B}^{3}-{B}^{4}
\nonumber \\
  \hspace{-0.98in}&&   \quad  \quad  \quad  \quad \quad 
\, \, -4608 \,AB  \, \, (A+B) \, +4096\,AB    \,  \, = \,\, \,  0,              
\end{eqnarray}
which is a representation of $\, \tau \, \rightarrow \, \, 3  \, \tau,$ where $\tau$ is the ratio of the two periods of the underlying elliptic function.

\section{Cut analysis for $\mathcal{H}_2(y)$} \label{app:H2}
\begin{figure}[htb]
\centering
\includegraphics[scale =0.54]{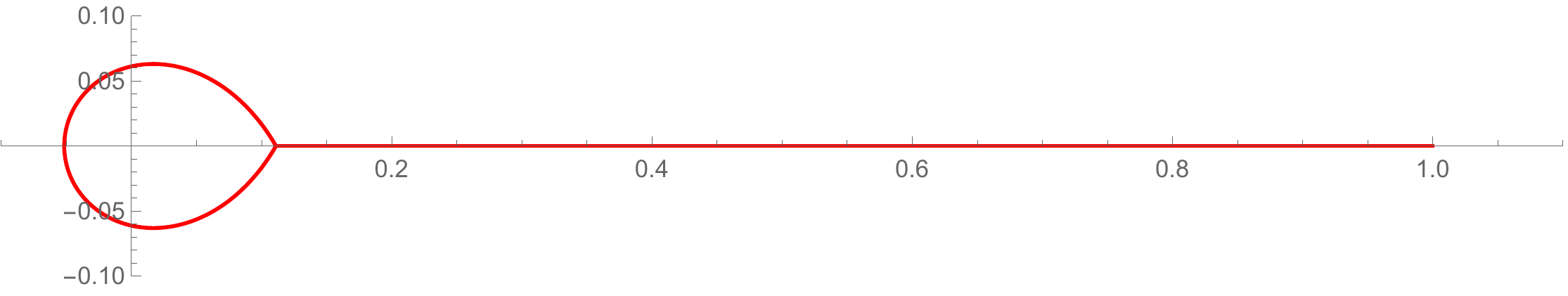} 
 \caption{The cuts of the functions $f_{1}(x)$ and ${\cal H}(x)$.}
 \label{fig:F1}
\end{figure}

We have the formula $f(x)=f_{1}(x):={\cal R}(x)+{\cal H}(x)$ near $x=0$. We cannot directly apply the Stieltjes inversion formula to this, because
\[{\cal H}(x)     \, = \, \,
       - \, {{(1\, -x)^{1/4} \,  \, (1\, -9\, x)^{3/4}} \over { 6 \, x^2}}  \,  \,
 _2F_1\Bigl(-\, {{1} \over {4}}, \, {{3} \over {4}}; \, 1; 
 \, \, {{ - \, 64 \, x } \over { (1\, -x) \,  \, (1 \, -9 \, x)^3}}  \Bigr)\]
is not analytic on $\mathbb{C}\setminus \mathbb{R}_{>0}$. Indeed ${\cal H}(x)$ has a cut wherever ${{ - \, 64 \, x } \over { (1\, -x) \,  \, (1 \, -9 \, x)^3}}\in[1,\infty)$. These cuts are shown in Fig.~\ref{fig:F1}. We can therefore only deduce that $f_{1}(x)$ coincides with the desired analytic extension $f(x)$ inside the ring in Fig.~\ref{fig:F1}.

\begin{figure}[htb]
\centering
\includegraphics[scale =0.44]{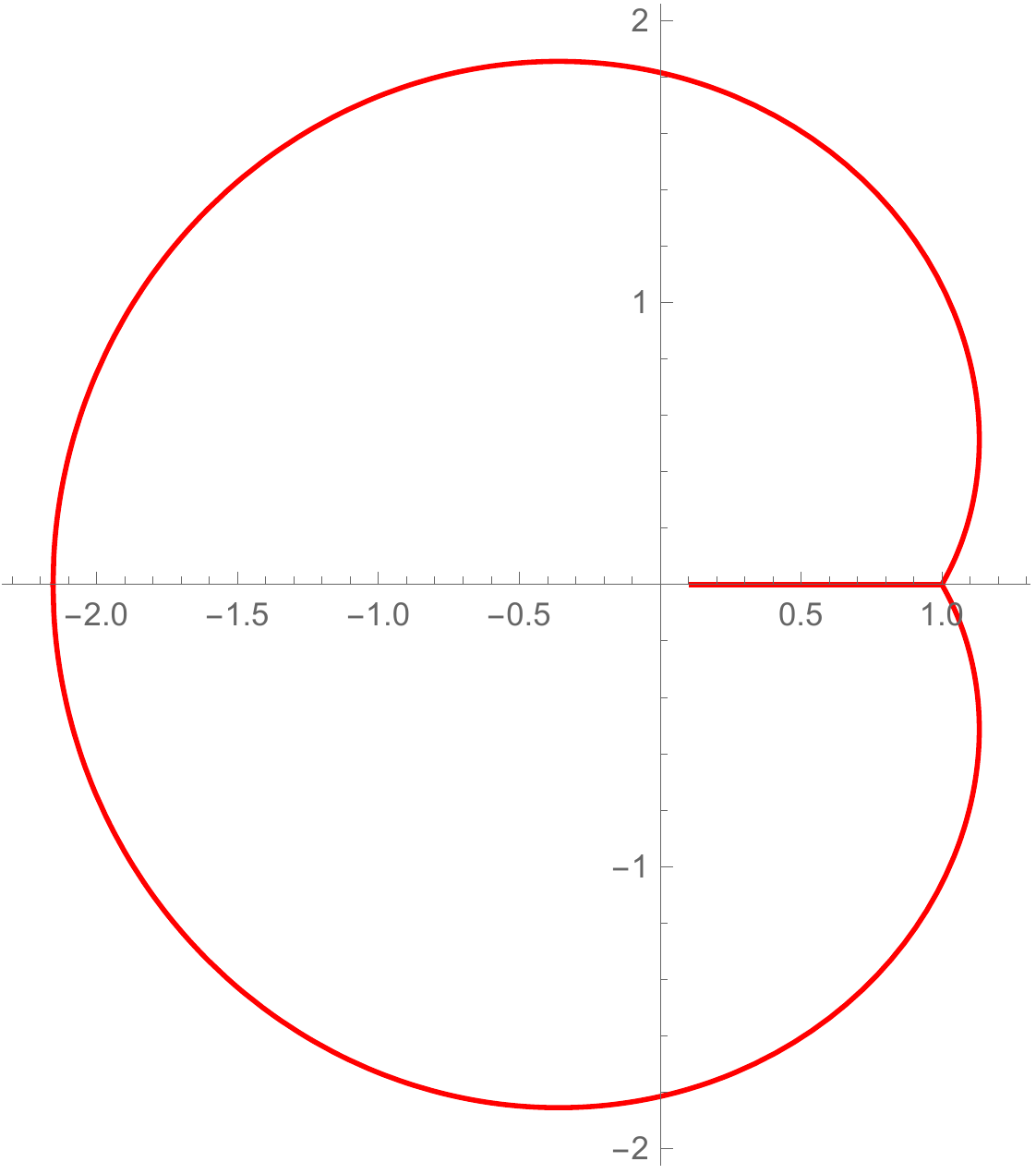} 
 \caption{The cuts of the functions $f_{2}(x)$ and ${\cal H}_{2}(x)$.}
 \label{fig:F2}
\end{figure}

Alternatively we can use the following formula which can be derived from \eqref{pullbackL21over9xidcalHter}: $f(x)=f_{2}(x):={\cal R}(x)+{\cal H}_{2}(x)$ near $x=0$, where 
\begin{align*}
{\cal H}_{2}(x)&=-\frac{(1-9x)^{1/4}(1+3x)^2}{6(1-x)^{5/4}x^2} 
\,  \, _2F_1\Bigl(-\, {{1} \over {4}}, \, {{3} \over {4}}; \, 1; 
 \, \, {{ - \, 64 \, x^3 } \over { (1\, -x)^3 \,  \, (1 \, -9 \, x)}}  \Bigr)\\
 &~+ {\frac{16x(1-6x-3x^2)}{(1-9x)^{3/4}(1-x)^{17/4}}} \,  \, _2F_1\Bigl(-\, {{3} \over {4}}, \, {{7} \over {4}}; \, 1; 
 \, \, {{ - \, 64 \, x^3 } \over { (1\, -x)^3 \,  \, (1 \, -9 \, x)}}  \Bigr).
\end{align*}
This function is also not analytic on $\mathbb{C}\setminus
\mathbb{R}_{>0}$, as it has a cut wherever ${{ - \, 64 \, x^3 } \over { (1\,
-x)^3 \,  \, (1 \, -9 \, x)}}\in[1,\infty)$. These cuts are shown in Fig.~\ref{fig:F2}. We can therefore only deduce that $f_{2}(x)$ coincides
with the desired analytic extension $f(x)$ inside the ring in
Fig.~\ref{fig:F2}. Nonetheless, this is an improvement on $f_{1}(x)$, as the
ring in Fig.~\ref{fig:F2} contains the cut $[1/9,1]$. As a consequence,
applying the Stieltjes inversion formula directly to $f(x)={\cal R}(x)+{\cal
H}_{2}(x)$ yields the correct expression for the density on the interval~$[1,9]$.

\section{The classical modular form emerging in $\, Av(12345)$}
\label{classical2}
Introducing $\, U_2^{(p)}$ the linear operator $\, U_2$ given by (\ref{U2def}) pulled-back
by $\, 64/x$:
\begin{eqnarray}
 \label{pU2}
 \hspace{-0.98in}&& \quad  \quad \quad
 U_2^{(p)}   \, \, = \, \,  \,  pullback\Bigl(U_2, \, {{64} \over {x}}\Bigr)
 \\
 \hspace{-0.98in}&& \quad  \quad  \quad \quad
\, \, \, \, = \, \,  \, \,
 D_x^2  \, \, \,
 +{\frac { \left( 128\,{x}^{2}-30\,x+1 \right)}{
     x \,  \, (1 \, - 4\,x) \,   \, (1 \, - 16\,x) }} \,  \,D_x
 \, \, \,  +\,{\frac {2 \,  \, (8\,x-1)}{
     x \,  \, (1 \, - 4\,x) \,   \, (1 \, - 16\,x) }}.
\nonumber       
\end{eqnarray}
One has the relations in (\ref{pU2}) and in $G(x)=\frac{ F(x)}{\sqrt{x}}$ given by (\ref{classmodformsquare}):
\begin{eqnarray}
 \label{pU2Homo}
\hspace{-0.98in}&&  \quad \quad \, \, 
U_2^{(p)}   \, \, = \, \,  \,   \, x^{-1/2}  \,  \,  U_2  \,  \, x^{1/2},
\quad \quad  \quad \quad
 G\Bigl( {{1} \over {64 \, x}}\Bigr)
\, \,  = \, \, \,  {{2(-1)^{1/3}}} \,  \, G(x). 
\end{eqnarray}
Consequently the classical modular form (the square of which is $\, F(x)$) is
\begin{eqnarray}
  \label{classmodform}
  \hspace{-0.98in}&& \quad  \quad  \quad  \quad   \,  \,  \,  \quad 
   x^{1/2} \,  \, (1\, -16 \, x)^{-1/6} \,  \, (1\, +2 \, x)^{-1/3}
  \nonumber \\
 \hspace{-0.98in}&& \quad  \quad \quad  \,  \,  \,   \quad  \quad   \quad  \quad   \quad  
 \,   \,  \, _2F_1\Bigl({{1} \over {6}}, \, {{2} \over {3}}; \, 1;  \, \,
  - \,{\frac { 108 \,  \, {x}^{2}}{ (1\, - 16\,x) \,  \, (1 \, +2\,x)^{2}}}   \Bigr),
\end{eqnarray}
{\rm which can also be written}:
\begin{eqnarray}
  \label{classmodformbis}
  \hspace{-0.98in}&& \quad \quad \quad \quad  
x^{1/2} \,  \, (1\, -4 \, x)^{-1/6} \,  \, (1\, +32 \, x)^{-1/3}
\nonumber \\
 \hspace{-0.98in}&& 
 \quad  \quad  \quad     \quad  \quad  \quad       \quad  \quad       
 \,  \,  \, _2F_1\Bigl({{1} \over {6}}, \, {{2} \over {3}}; \, 1;  \, \,
  \, \,{\frac {108\, x}{ (1\, - 4\,x) \,  \, (1 \, +32\,x)^{2}}}   \Bigr).  
\end{eqnarray}
The relation between the two (Hauptmodul) pullbacks
\[\, \, A \, = \,  - \,{\frac { 108\, {x}^{2}}{ (1\, - 16\,x) \,  \, (1 \, +2\,x)^{2}}}  \qquad
\text{and} \qquad \, \, B\, = \,  \,{\frac {108\, x}{ (1\, - 4\,x) \,  \, (1 \, +32\,x)^{2}}} \, \]
is the (fundamental) modular equation
\begin{eqnarray}
  \label{modular}
  \hspace{-0.9in}&& \quad \quad \quad  
  625\, \, \, {A}^{3}{B}^{3} \, \, -525\, {A}^{2}{B}^{2}\, \, (A+B) \,  \,  \, \, \,
  -3\,AB \left( 32\,{A}^{2}+AB+32\,{B}^{2} \right)
 \nonumber \\
\hspace{-0.98in}&& \quad \quad \quad \quad \quad \quad  
\,  \,  \, -4\, \, \, (A+B)  \, ({A}^{2}-133\,AB+{B}^{2})
\,  \, \,  -432\,AB \,  \, \, = \, \,  \, \, 0.
\end{eqnarray}
which is a representation of $\,\, \tau \, \rightarrow \, \, 2 \, \tau$.
Note that one could have seen directly that~$\, L_3$, the order-three linear
differential operator (\ref{L3def}), is non-trivially homomorphic to its
pullback by $\, x \, \rightarrow \, \frac{1}{64x}$.

\section{Densities for (Hamburger) moment sequences of walks}
\label{sec:hamburger} 

Fig.~\ref{fig:walks} displays Hamburger moment sequences (not Stieltjes!),
coming from the enumeration of walks restricted to the quarter plane. The
context is the following: given a certain model (i.e., set of allowed steps), one
forms the sequence $(a_n)_{n \geq 0}$ whose $n$-th term counts the walks
confined to the quarter plane $\mathbb{N}^2$, starting at the origin, and
consisting of exactly $n$ steps. In general, the generating function $F(z) = 
\sum_n a_n z^n$ is not algebraic, and not even D-finite. However, in
Fig.~\ref{fig:walks} we consider only models leading to algebraic generating
functions. These are basically of two types: either the support of the model
is contained in a half-plane, e.g. for the first model (\oeis{A001405}), in
which case the generating function is known to be algebraic~\cite{BaFl02}; or
the support has a more general shape, and then algebraicity is not obvious,
such as for 11th model (\oeis{A001006})~\cite[Prop.~9]{BMM}, or for the last
model (\oeis{A151323})~\cite[Prop.~15]{BMM}, see also~\cite{BoKa09}.

In each case we have computed the density function using the Stieltjes inversion formula. We describe below this process for the sixth model \begin{tiny}{\TwoDstepset00100011}\end{tiny}, as in this case the density involves a Dirac delta function, which has not appeared in our previous examples. Starting with the generating function
\[F(z)=\frac{4z-1+\sqrt{1-8z^2}}{4z(1-3z)}\]
for the model, we write $G_{\mu}(z)=\frac{1}{z}F(\frac{1}{z})$, as then the density $\mu(x)$ should be given by
\[\mu(x)=\lim_{\epsilon\to0^{+}}-\frac{1}{\pi}\, \Im(G_{\mu}(x+i\epsilon)).\]
For $x\neq3$, this yields the correct formula, however, $G_{\mu}(x)$ has a pole at $x=3$, as the limit does not exist at this point. There are several ways to resolve this problem (all of which, of course, yield the same density), the simplest being to remove the pole from $F(z)$ before applying the Stieltjes inversion formula. To do this we write 
\[F(z)=\tilde{F}(z)+\frac{1}{2(1-3z)},\]
as then $\tilde{F}(z)$ has no poles. Applying the Stieljes inversion yields
\[\tilde{F}(z)=\int_{-2\sqrt{2}}^{2\sqrt{2}}\frac{1}{1-xz}\tilde{\mu}(x)dx,
\qquad \text{where} \quad \tilde{\mu}(x)=\frac{1}{4\pi}\, \frac{\sqrt{8-x^2}}{3-x}.\]
We then simply observe that
\[\frac{1}{2(1-3z)}=\int_{0}^{\infty}\frac{1}{1-xz}\cdot \frac{1}{2}\delta(x-3)dx,\]
where $\delta$ is the Dirac delta function. It follows that the density function $\mu(x)$ is
\[\mu(x)=\tilde{\mu}(x)+\frac{1}{2}\delta(x-3).\]

Actually, much more can be said about the 13 sequences displayed in Fig.~\ref{fig:walks}.
For any sequence $\seq$ among them, the leading principal minors
$\Delta_0^n(\seq)$ and $\Delta_1^n(\seq)$ have very nice expressions:
\begin{itemize}
 \item in all cases except the last one, $\Delta_0^n(\seq)$ is equal to $q^{\binom{n}{2}}$, where $q$ is 1, 2, 3, or 4;
 \item 	in all cases except the last one, the quotient sequence $u_n:=\Delta_1^n(\seq)/\Delta_0^n(\seq)$ is a linearly recurrent sequence with constant coefficients.
\end{itemize}	
For instance, for the tenth model \begin{tiny}\TwoDstepset11101001\end{tiny}, 
with sequence \oeis{A151292} = $(1, 2, 7, 23, 85, 314, \ldots)$:
\begin{itemize}
\item the sequence $\Delta_0^n(\seq)$ is 
$3^{\binom{n}{2}} = (1, 3, 27, 729, 59049, 14348907,\ldots)$,
\item the sequence $\Delta_1^n(\seq)$ is 
$(2, -3, -189, -2916, 1003833, 416118303, \ldots)$
with quotient $u_n:=\Delta_1^n(\seq)/\Delta_0^n(\seq)$
equal to
$(2, -1, -7, -4, 17, 29, -22, -109, -43, \ldots)$
which satisfies the recurrence $u_{n+2} =  u_{n+1} - 3 u_n$
with initial conditions $u_0 = 2, u_1 = -1$.
\end{itemize}
For the last model \begin{tiny}\TwoDstepset11011011\end{tiny},
with sequence \oeis{A151323} = $(1, 3, 14, 67, 342, 1790, \ldots)$, 
the situation is even more interesting: 
\begin{itemize}
\item the sequence $\Delta_0^n(\seq)$ is 
$(1, 5, 105, 9009, 3128697, 4379132901,\ldots)$
and it appears to coincide with sequence \oeis{A059490},
which emerges in the context of \emph{alternating symmetric matrices} (ASMs).
Precisely, it is Kuperberg's sequence $\text{AQT}^{(1)}(4n;2)$ 
from  \cite[Theorem~5]{Kuperberg} 
(counting quarter-turn symmetric ASMs),
which equals
\begin{equation}\label{eq:AQT}
 \left( -4 \right) ^{{n\choose 2}}\,  
\prod _{i=1}^{n}\prod _{j=1}^{n}
{\frac {4\,(j-i)+1}{j-i+n}}.
\end{equation}
\item the sequence $\Delta_1^n(\seq)$ is 
$(3, 5, -1113, -227799, -23986677, 379126748429, \ldots)$
and has quotient $u_n:=\Delta_1^n(\seq)/\Delta_0^n(\seq)$ equal to
\[ \left(3, 1, -\frac{53}{5}, -\frac{177}{7}, -\frac{23}{3}, \frac{2857}{33}, \frac{29169}{143}, \frac{3921}{65},\ldots\right),\]
which no longer satisfies a recurrence with constant coefficients, but still
(conjecturally) satisfies a linear recurrence with \emph{polynomial} coefficients
\begin{equation}\label{eq:recAQT}
	 u_{n+2} = 2 \, u_{n+1} - \frac{(4n+9)(4n+7)}{(2n+5)(2n+3)} u_n.
\end{equation}
\end{itemize}
In cases 1--12, all these assertions can be proved using the explicit form of
the corresponding (Stieltjes, or Jacobi) continued fractions, which can be
determined since the generating function of $\seq$ is an algebraic function
of degree~2. This approach is very classical, and is for instance
explained by Krattenthaler in~\cite[\S2.7]{Krat99} and \cite[\S5.4]{Krat05}.
It was applied in similar walk enumeration contexts by Tamm~\cite{Tamm01}, by
Brualdi and Kirkland~\cite{Brualdi}, and by Chang, Hu, Lei and
Yeh~\cite{CHLY}, to name just a few.

Going beyond algebraicity degree 2 appears to be quite a challenging task.
For degree 3, namely for some sequences in the orbit of the Catalan-Fuss
sequence $\frac{1}{3n+1}\binom{3n+1}{n}$, explicit evaluations of Hankel
determinants have been possible due to an approach introduced by
Tamm~\cite{Tamm01} and extended by Gessel and Xin~\cite{GesselXin06}, see
also~\cite{ETR} and Theorem~31 in~\cite{Krat05}. Some of these evaluations are
already connected to the enumeration of alternating sign matrices presenting
(vertical) symmetries.

{For the evaluations~\eqref{eq:AQT}
and~\eqref{eq:recAQT} of the Hankel determinants occurring in case~13, we can again use the J-fraction
\begin{equation}\label{Jfrac}A(x) \, = \,  \,
    \cfrac{\beta_0}{1 \, -\gamma_0 x\, -\cfrac{\beta_1 x^2}{1- \gamma_1 x-\, \cfrac{\beta_2x^2}{
      \begin{array}{@{}c@{}c@{}c@{}}
        1 \, - \cdots 
      \end{array}
    }}} , 
\end{equation}
with $\gamma_{0}=3$, $\gamma_{j}=2$ for $j>0$ and
$\beta_{j}=\frac{(4j-1)(4j+1)}{(2j-1)(2j+1)}$. Similarly to cases 1--12, the
Hankel determinants can easily be determined from these continued fraction
coefficients. To prove that this is indeed the correct continued fraction, we
adapt a method applied implicitly by Euler \cite[\S21]{Euler} to S-fractions,
which is described in detail in \cite[\S2.1]{PSZ} and \cite{Sokalsimple}. The
first step in our case is to set $A_{1}(x):=A(x)/\beta_{0}$, and define
$A_{j+1}(x)$ recursively by
\begin{equation}\label{Arecursion}A_{j}(x)=\frac{1}{1-\gamma_{j-1}x-\beta_{j}x^2 A_{j+1}(x)}.\end{equation}
Then it suffices to show that each $A_{j}(x)\in\mathbb{R}[x]$, as this immediately yields the continued fraction form \eqref{Jfrac} for $A(x)=\beta_{0}A_{1}(x)$. The next step is to define another sequence $\{B_{n}(x)\}_{n\in\mathbb{N}}$ by $B_{0}(x)=1$ and $B_{j}(x)=A_{j}(x)B_{j-1}(x)$. Writing $A_{j}(x)=B_{j}(x)/B_{j-1}(x)$, the recursion \eqref{Arecursion} simplifies to
\begin{equation}\label{Brecursion}
	B_{j+1}(x)=\frac{(2j-1)(2j+1)}{(4j-1)(4j+1)x^2} \, \Big( (1-2x)B_{j}(x)-B_{j-1}(x) \Big).
\end{equation}
Finally, we guess and prove the exact form of each series $B_{j}(x)$. First, the guess is that $B_{j}(x)$ is the unique series in $\mathbb{R}[x]$ with constant term $1$ satisfying
\begin{equation}\label{Bdiff}-\frac{2j(2j+1+6(j+1)x)}{x(1-6x)(1+2x)}B_{j}(x)+\frac{2j-(8j+6)x-24(j+1)x^2}{x(1-6x)(1+2x)}B_{j}'(x)+B_{j}''(x)=0.
\end{equation}
To see that this equation uniquely defines $B_{j}(x)$, we observe that it is equivalent to the coefficients $b_{j,0},b_{j,1},\ldots$ of $B_{j}(x)$ being determined by the recurrence
\[b_{j,n+1}=\frac{2(2n+2j+1)(n+j)b_{j,n}
+12(n+j)(n+j-1)b_{j,n-1}}{(2j+n)(n+1)},\]
with the initial conditions $b_{j,0}=1$ and $b_{j,-1}=0$.}

Finally, we prove by induction that \eqref{Bdiff} holds for all $j$. The base cases are trivial to prove using the exact, algebraic forms of $B_{0}(x)$ and $B_{1}(x)$, so we proceed to the inductive step. We assume that \eqref{Bdiff} holds for $j$ and $j-1$, and prove that it holds for $j+1$. Using \eqref{Brecursion}, this reduces to proving
\[B_{j-1}(x)=\frac{x(1+2x)(1-6x)B_{j}'(x)-(1-2j-x+6jx+12jx^2)B_{j}(x)}{2j-1},\]
which is true as the right hand side satisfies the differential equation \eqref{Bdiff} for $B_{j-1}(x)$ and has constant term $1$.

\newcommand{\fnote}[1]{{\footnotesize #1}}
\newcommand{\seqU}[2]{\renewcommand{\arraystretch}{0.3}
\renewcommand{\arraystretch}{1.2}
\begin{tabular}{l} \qquad \qquad #1 \\[-1.0truemm] \fnote{#2} \end{tabular}}
\begin{figure}\small

\renewcommand{\tabcolsep}{0.1pt}

\vspace{1.5cm}

\begin{center}
\begin{tabular}{ll|c|cc}
\hline\hline
\multicolumn{2}{c|}{{Model, sequence and OEIS tag}} & Generating function $F(z)$ & \multicolumn{1}{c}{{Density $\mu(x)$}} & 
\multicolumn{1}{c}{Support $\Gamma$}
\\
\hline\hline
\oeis{A001405} & \seqU{\TwoDstepset00100001}{$(1, 1, 2, 3, 6, 10, 20, 35, 
	% 70, 
\ldots)$}
 &
$\ds \frac{2z-1+\sqrt{1-4z^2}}{2z(1-2z)}$ &
$\ds \frac{1}{2\pi}\sqrt{\frac{2+x}{2-x}}$ & $[-2,2]$ 
\\
\oeis{A005773} & \seqU{\TwoDstepset00101001}{$(1, 2, 5, 13, 35, 96, 267, 
	% 750, 
\ldots)$}
 &
$\ds \frac{3z-1+\sqrt{1-2z-3z^2}}{2z(1-3z)}$ &
$\ds \frac{1}{2\pi}\sqrt{\frac{1+x}{3-x}}$ & $[-1,3]$ 
\\
\oeis{A151318} & \seqU{\TwoDstepset01101011}{$(1, 3, 13, 55, 249, 1131, 
	% 5253, 
\ldots)$}
 &
$\ds \frac{5z-1+\sqrt{1-2z-15z^2}}{4z(1-5z)}$ &
$\ds \frac{1}{4\pi}\sqrt{\frac{3+x}{5-x}}$ & $[-3,5]$ 
\\
\oeis{A060899} & \seqU{\TwoDstepset01100011}{$(1, 2, 8, 24, 96, 320, 
	% 1280, 
	% 4480, 
\ldots)$}
 &
$\ds \frac{4z-1+\sqrt{1-16z^2}}{4z(1-4z)}$ &
$\ds \frac{1}{4\pi}\sqrt{\frac{4+x}{4-x}}$ & $[-4,4]$ 
\\
\oeis{A126087} & \seqU{\TwoDstepset01100001}{$(1, 1, 3, 5, 15, 29, 87, 181, 
	%543, 
\ldots)$}
 &
$\ds \frac{2z-1+\sqrt{1-8z^2}}{2z(1-3z)}$ &
$\ds \frac{1}{2\pi}\frac{\sqrt{8-x^2}}{3-x}$ & $[-2\sqrt{2},2\sqrt{2}]$ 
\\
\oeis{A151281} & \seqU{\TwoDstepset00100011}{$(1, 2, 6, 16, 48, 136, 408, 	
	% 1184, 
\ldots)$}
 &
$\ds \frac{4z-1+\sqrt{1-8z^2}}{4z(1-3z)}$ &
$\ds {\frac{1}{4\pi}\frac{\sqrt{8-x^2}}{3-x} + \frac12 \delta(x-3)}$ & {$[-2\sqrt{2},2\sqrt{2}]$}
\\ 
\oeis{A128386} & \seqU{\TwoDstepset11100001}{$(1, 1, 4, 7, 28, 58, 232, 
	% 523, 
\ldots)$}
 &
$\ds \frac{2z-1+\sqrt{1-12z^2}}{2z(1-4z)}$ &
$\ds \frac{1}{2\pi}\frac{\sqrt{12-x^2}}{4-x}$ & $[-2\sqrt{3},2\sqrt{3}]$ 
\\
\oeis{A151282} & \seqU{\TwoDstepset01101001}{$(1, 2, 6, 18, 58, 190, 638, 
	% 2170, 
\ldots)$}
 &
$\ds \frac{3z-1+\sqrt{1-2z-7z^2}}{2z(1-4z)}$ &
$\ds \frac{1}{2\pi}\frac{\sqrt{7+2x-x^2}}{4-x}$ & $[1-2\sqrt{2},1+2\sqrt{2}]$ 
\\ 
\oeis{A129637} & \seqU{\TwoDstepset00101011}{$(1, 3, 11, 41, 157, 607, 
	% 2367,  
\ldots)$}
 &
$\ds \frac{5z-1+\sqrt{1-2z-7z^2}}{4z(1-4z)}$ &
$\ds {\frac{1}{4\pi}\frac{\sqrt{7+2x-x^2}}{4-x} }$ & 
% {$[1-2\sqrt{2},1+2\sqrt{2}] \cup \{4\}$}
\\
 & 
 &
 &
$\,\,\,  {+ \frac12 \delta(x-4)}$ & {\hspace{-0.2cm}$[1-2\sqrt{2},1+2\sqrt{2}] \cup \{4\}$}
\\[+0.8ex] 
\oeis{A151292} & \seqU{\TwoDstepset11101001}{$(1, 2, 7, 23, 85, 314, 
	% 1207, 
	% 4682, 
\ldots)$}
 &
$\ds \frac{3z-1+\sqrt{1-2z-11z^2}}{2z(1-5z)}$ &
$\ds \frac{1}{2\pi}\frac{\sqrt{11+2x-x^2}}{5-x}$ & $[1-2\sqrt{3},1+2\sqrt{3}]$
\\
\oeis{A001006} & \seqU{\TwoDstepset00110010}{$(1, 1, 2, 4, 9, 21, 51, 127, 
	% 323, 
\ldots)$}
& $\ds\frac{1-z-\sqrt{1-2z-3z^2}}{2z^2} $ 
& $\ds \frac{1}{2\pi}\sqrt{(3-x)(1+x)} $ & $[-1,3]$ 
\\
\oeis{A129400} & \seqU{\TwoDstepset01111110}{$(1, 2, 8, 32, 144, 672, 
	% 3264, 
\ldots)$}
 &
$\ds \frac{1-2z-\sqrt{1-4z-12z^2}}{8z^2}$ &
$\ds \frac{1}{8\pi}{\sqrt{(x+2)(6-x)}}$ & $[-2,6]$ 
\\ 
\oeis{A151323} & \seqU{\TwoDstepset11011011}{$(1, 3, 14, 67, 342, 1790, 
	% 9580, 
\ldots)$}
 &
$\ds  \frac{\sqrt[4]{\frac{1+2z}{1-6z}}-1}{2z}$ &
$\ds \frac{1}{2\sqrt{2}\pi}\sqrt[4]{\frac{2+x}{6-x}}$ & $[-2,6]$ 
\\ 
\hline\hline
\end{tabular} 
\end{center}
\caption{\label{fig:walks}\small
{
Some sequences of walks in $\N^2$, counted by length, and their (Hamburger)
moment representations: walk model, initial terms, tag in the On-Line Encyclopedia
of Integer Sequences (\href{https://oeis.org}{OEIS}), (algebraic) generating
function, associated (algebraic) density function, and its support.}
}
\end{figure}

\end{document}